\documentclass[11pt]{article}
\usepackage[utf8]{inputenc}
\usepackage{amsmath}
\usepackage{authblk}
\usepackage{amssymb}
\usepackage[english]{babel}
\usepackage{hyperref}
\IfFileExists{apacite.sty}%
  {\usepackage{apacite}\def\bibstylename{apacite}}%
  {\usepackage[round]{natbib}\def\bibstylename{apalike}}
\providecommand{\citeA}{\citet}
\usepackage{ragged2e}
\usepackage{enumerate}
\usepackage{graphicx}
\usepackage{geometry}
\usepackage{pifont}
\usepackage{setspace}
\usepackage{xcolor}
\hypersetup{colorlinks=true,linkcolor=blue,anchorcolor=blue,citecolor=blue}
\newtheorem{example}{Example}
\newtheorem{definition}{Definition}

\newtheorem{theorem}{Theorem}
\newtheorem{lemma}{Lemma}
\newtheorem{proposition}{Proposition}

\newtheorem{corollary}{Corollary}

\newcommand{\M}{\mathbf{M}}

\renewcommand{\P}{\mathbb{P}}

\newcommand{\<}{\langle}
\renewcommand{\>}{\rangle}
\newcommand{\E}{\mathbb{E}}

\begin{document}
\title{Self-Normalized Inequalities for Time Series with Sample Variance Only and Their Application to Forecasting}
\author[1]{Zihao Yuan}
\affil[1]{University of Hamburg}
\date{}
\maketitle

\begin{abstract}
Evaluating a forecasting rule on a time series requires a risk bound that can be computed from the observed losses. Such generalization bounds rest on empirical Bernstein (EB) and self-normalized inequalities, which replace the unknown variance by a sample counterpart. These tools, however, have been developed almost exclusively for independent observations and martingale differences, and, for dependent data, the natural setting of forecasting, direct results are scarce. They require $\beta$-mixing and knowledge of the mixing coefficients, and their leading terms are not sharp. This paper fills this gap. We first establish, for martingale differences, an EB inequality in closed form that uses only the realized quadratic variation, holds with probability at least $1-\alpha$ for every sample size, and is sharp under heteroscedasticity, whereas the sharpness of existing EB inequalities is established for a constant conditional variance. Building on this, we derive EB and self-normalized inequalities for the marginal mean of $\tilde{\phi}$- and $\tau$-mixing processes, a class that contains processes that are not even strongly mixing, such as autoregressions with Bernoulli innovations. The radii use nothing but a block sample variance: no variance proxy has to be designed or tuned, and, for $\tau$-mixing processes, no dependence coefficient enters the radius. Under summability-type conditions on the dependence coefficients, their leading terms carry the constant of the oracle Bernstein inequality with the long-run variance, so, to our knowledge, these are the first sharp self-normalized inequalities under non-strong mixing conditions. As an application, we obtain a sharp upper bound on the risk of a time-series forecasting rule that is computed from its observed losses with a single block length and contains no dependence coefficient.
\end{abstract}

\section{Introduction}
\label{sec 1}
How far a sample mean can deviate from its expectation is the basic question behind every generalization bound. Classical answers, such as the inequalities of Bernstein and Bennett, are stated in terms of the variance, which is unknown in practice. Replacing it by a worst-case bound, such as the squared range of the data, keeps the bound valid but ignores how concentrated the data actually are. Self-normalized and empirical Bernstein (EB) inequalities avoid this difficulty by measuring the deviation against the realized quadratic variation or the sample variance, so that the bound adapts to the data and can be computed from them. This makes EB bounds a natural tool for certifying the performance of learned rules, for instance through the sample-variance penalization of \citeA{maurer2009empirical}.

Most of this theory, however, is built for independent observations or martingale differences, see, for instance, \citeA{howard2021time}, \citeA{waudby2024estimating}, \citeA{wang2024sharp} and \citeA{chugg2025variational}. Time series do not fit this frame. The losses of a forecasting rule evaluated along a time series are autocorrelated, and the variance that governs their average is the long-run variance rather than the marginal one. The EB inequality of \citeA{mirzaei2025empirical} for dependent data requires $\beta$-mixing and knowledge of the mixing coefficients, and therefore excludes processes that are not strongly mixing, such as the autoregression with Bernoulli innovations of \citeA{andrews1984non}. This paper fills this gap. We derive self-normalized and EB inequalities for the marginal mean of $\tilde{\phi}$- and $\tau$-mixing processes, a class that contains such processes, with leading terms that only involve a block sample variance and are sharp with respect to the long-run variance. We then use them to evaluate time-series forecasting rules directly from their observed losses.

Sections \ref{sec 1.1} and \ref{sec 1.2} review the two lines of work on which we build, and Section \ref{sec 1.3} summarizes our contributions.
\subsection{Self-Normalized Inequality}
\label{sec 1.1}
Suppose $(\Omega,\mathcal{F}, (\mathcal{A}_i)_{i\in \mathbb{Z}}, \mathbb{P})$ is a filtered probability space and $(Z_{i})_{i\in \mathbb{Z}}=:\mathbf{Z}$ is a real-valued stochastic process adapted to $(\mathcal{A}_i)_{i\in\mathbb{Z}}$. Denote
\begin{align*}
    \mathbf{M}_n=\sum_{i=1}^n(Z_i-\E[Z_i]),\quad [\mathbf{M}]_n=\sum_{i=1}^n(Z_i-\E[Z_i])^2,\quad
    \<\mathbf{M}\>_n=\sum_{i=1}^n\E[(Z_i-\E[Z_i])^2|\mathcal{A}_{i-1}].
\end{align*}
Generally speaking, the central goal of scalar self-normalized inequality is to show some concentration phenomenon of mean based on $\M_n$, $[\M]_n$ and $\<\M\>_n$. When $(Z_i,\mathcal{A}_i)_{i\in \mathbb{Z}}$ is a martingale difference sequence (MDS) such that $|Z_i|\leq b$, \citeA{freedman1975tail} established an exponential upper bound on $\P(\M_n\geq x, \<\M\>_n\leq v)$ for any $x,v\geq 0$. More specifically, when $\E[Z_i^2|\mathcal{A}_{i-1}]=\sigma^2>0$ holds, we have the corollary
\begin{align}
    \label{eq 1}
    \P\left(\frac{1}{n}\M_n\leq \sqrt{\frac{2\sigma^2\log(\frac{1}{\delta})}{n}}+\frac{b\log(\frac{1}{\delta})}{3n}\right)\geq 1-\delta.
\end{align}
Then, \citeA{de1999general} and \citeA{dzhaparidze2001bernstein} extended Freedman's work to unbounded cases. These results provide powerful concentration and limit theorems for Studentized partial sums and self-normalized martingales. However, the variance term typically appears as a predictable variance process or as an abstract conditional variance, which is rarely observable and often lives only as a process adapted to a given filtration. Recall that, when $(Z_i,\mathcal{A}_i)_{i\in \mathbb{Z}}$ is an MDS, $\E[Z_i]=0$ holds for every $i$, which implies that $[\M]_n=\sum_{i=1}^nZ_i^2$ is fully observable. Hence, deriving a concentration inequality purely based on $\M_n$ and $[\M]_n$ becomes a meaningful task. By strengthening the assumption of MDS to conditionally symmetric MDS, \citeA{de1999general} proved that, for all $x>0$,
\begin{align}
    \label{eq 2}
    \P\Big(\frac{\M_n}{[\M]_n}\geq x\Big)\leq \sqrt{\E\Big[\exp\Big(-\frac{1}{2}[\M]_nx^2\Big)\Big]}
\end{align}
and, for all $x,y\geq 0$,
\begin{align}
    \label{eq 3}
    \P\Big(\frac{\M_n}{[\M]_n}\geq x,[\M]_n\geq y\Big)\leq \exp\Big(-\frac{1}{2}x^2y\Big).
\end{align}
Since a conditionally symmetric MDS satisfies the definition of ``heavy-on-left'' in \citeA{bercu2008exponential}, they obtained the following important extension of \citeauthor{de1999general}'s \eqref{eq 2},
\begin{align}
    \label{eq 4}
     \P\Big(\frac{\M_n}{[\M]_n}\geq x\Big)\leq \inf_{p\geq 1}\Big(\E\Big[\exp\Big(-\frac{1}{2}(p-1)[\M]_nx^2\Big)\Big]\Big)^{\frac{1}{p}}.
\end{align}
They also proved that \eqref{eq 3} holds for ``heavy-on-left'' MDS as well. Compared with \eqref{eq 1}, \eqref{eq 2}-\eqref{eq 4} hold under much weaker moment conditions, but their statistical applications are not as convenient as \eqref{eq 1}. Meanwhile, the dependence on $[\M]_n$ appears inside the Laplace transform or in the event $\{[\M]_n\geq y\}$, and there is no closed-form expression of the type $\sqrt{[\M]_n\log(\frac{1}{\delta})}$ with an explicit confidence parameter $\delta$. However, the variance ($\sigma$) and boundedness ($b$) are usually unknown to us, which prevents the direct application of \eqref{eq 1}. More recently, \citeA{howard2020time} derived a class of exponential bounds for the probability that a martingale sequence crosses a time-dependent linear threshold. Generally speaking, their work strengthens nearly all of the aforementioned results to a sharper level under more general conditions. For one-parameter self-normalized inequalities involving a martingale $\M_n$ and $\<\M\>_n$, the key contribution of \citeauthor{howard2020time}'s paper is to give a single time-uniform (infinite-horizon) Chernoff/line-crossing framework that systematically produces bounds where the deviation of $\M_n$ is controlled by a variance process $V_n$ built from $\<\M\>_n$. The connection between \citeauthor{howard2020time}'s framework and our results is discussed after Proposition \ref{prop 1} and Theorem \ref{th 3}.

It is obvious that, under some circumstances, the value of boundedness can still be known from prior information. But, compared with boundedness, it is more challenging and risky to assume the variance to be known. Thus, it would be very ideal to derive an inequality having the form \eqref{eq 1} but in which $\sigma$ is observable (or data-dependent), which is actually the core task of a parallel but tightly connected area, empirical Bernstein inequality of marginal mean.
\subsection{Empirical Bernstein Inequality of Marginal Mean}
\label{sec 1.2}
Based on a variance estimator, \citeA{maurer2009empirical} derived an EB inequality for IID random variables taking values in $[0,1]$. However, compared with \eqref{eq 1}, their inequality is not sharp in the following sense, which is used throughout the paper. Let $L_n$ denote the length of a one-sided confidence interval for the mean of $n$ observations with marginal variance $\sigma^2$. The interval is called sharp if
\begin{align}
\label{eq 5}
    \sqrt{n}L_n\xrightarrow[]{\P} \sigma\sqrt{2\log(1/\alpha)},
\end{align}
and, for two-sided intervals of level $1-\alpha$ with half-length $L_n$, the constant $\sqrt{2\log(1/\alpha)}$ is replaced by $\sqrt{2\log(2/\alpha)}$. The interval of \citeA{maurer2009empirical} does \textbf{not} satisfy \eqref{eq 5}. Based on self-normalized martingale techniques, Theorem 4 of \citeA{howard2021time} establishes an anytime EB bound for the difference between sample mean and averaged conditional mean, that is, based on the $\mathbf{Z}$ introduced in Section \ref{sec 1.1}, the following inequality holds for any given $\delta>0$,
\begin{align}
\label{eq 6}
    \P\Big(\forall\ n\geq 1:\Big|\frac{1}{n}\sum_{i=1}^nZ_i-\frac{1}{n}\sum_{i=1}^n\E[Z_i|Z_{i-1},...,Z_{1}]\Big|\leq \frac{u_{\delta}(V_n)}{n}\Big)\geq 1-2\delta,
\end{align}
where $(V_i)_{i\geq 1}$ is a variance process adapted to the filtration $\{\sigma(Z_1,...,Z_i)\}_{i\geq 1}$. The function $u_\delta$ is essentially a sub-exponential uniform boundary and is user-defined (see Sections 2 and 3 of \citeA{howard2021time} for more details). It is obvious that, when $\E[Z_i|Z_{i-1},...,Z_1]$ is invariant with respect to $i$ (e.g. IID), \eqref{eq 6} becomes an anytime EB inequality of marginal mean. Because \eqref{eq 6} is designed to be stopping-time valid, its fixed-time instantiation typically incurs an iterated-logarithm (or other uniform-boundary) overhead relative to a sharp fixed-$n$ bound. Consequently, obtaining sharp fixed-$n$ empirical Bernstein bounds for the marginal mean is nontrivial. A sharp resolution in the bounded case was later provided by \citeA{waudby2024estimating}. More specifically, Theorem 2 and Remark 1 of \citeauthor{waudby2024estimating} give a nice answer to this question (their results include more general weighted empirical means, but we here only focus on the simple empirical mean for the ease of comparison). By assuming that the $Z_i$'s are bounded random variables taking values in $[0,1]$, they proved the following EB inequality of marginal mean for all $\delta>0$ when $\E[Z_i|Z_{i-1},...,Z_1]=\mu$,
\begin{align}
    \label{eq 7}
    \P\left(\frac{1}{n}\sum_{i=1}^nZ_i-\mu\leq \sqrt{\frac{2\log (\frac{1}{\delta})V_{n,\delta}}{n}}\right)\geq 1-\delta.
\end{align}
Here $V_{n,\delta}$ is a quantity depending on the sample and the level $\delta$. More specifically, when the conditional variance is invariant, i.e., $\text{Var}(Z_i|Z_{i-1},...,Z_1)=\sigma^2$, $V_{n,\delta}$ converges to $\sigma^2$ almost surely. Thus, \eqref{eq 7} is sharp in the sense of \eqref{eq 5}. To our knowledge, \citeA{waudby2024estimating} were the first to prove that their EB inequality is asymptotically sharp.

In summary, the results of Ramdas and co-authors are indeed elegant. But, as we discussed for the fixed-$n$ EB inequality of marginal mean, they primarily target the predictable mean (the arithmetic average of conditional expectations). Unfortunately, in general mixing settings, there is no clear high-probability control of the gap between this predictable mean and the marginal mean, unless one assumes that all conditional expectations are constant, which essentially reduces to the martingale-difference case. Their recent important works, e.g. \citeA{chugg2025variational} and \citeA{wang2024sharp}, still rely on the assumption of constant conditional expectation. More specifically, it seems that their sharpness relies on a carefully tuned variance proxy satisfying the following three points simultaneously: (i) it is predictable, (ii) it ensures the super-martingale property, and (iii) it converges to the true variance. Thus, the design of this proxy is non-trivial even for martingale differences or i.i.d.\ data. For mixing processes, whose relevant variance of partial sums is the long-run variance rather than the marginal variance, the difficulty is likely to become more significant.

\par As for EB inequalities of dependent data, \citeA{mirzaei2025empirical} derived EB inequalities for $\beta$-mixing Hilbert-valued processes based on block techniques. Their bounds depend explicitly on the knowledge of mixing coefficients and the block length, and the effective sample size is reduced to the number of blocks, which in general prevents the parametric rate (see Theorems 1--3 there). In their bounds, the $\beta$-mixing coefficient enters through a reduced confidence level, and the leading term involves an uncentred second moment within blocks, so that it is not sharp in the sense of \eqref{eq 5}. Moreover, $\beta$-mixing implies strong mixing, so that processes that are not strongly mixing, such as the autoregression of Example \ref{example 1}, are not covered.

\subsection{Our Contributions}
\label{sec 1.3}
The main contributions of this paper are:

 \begin{itemize}
     \item [C1] (\textbf{Sharp EB inequalities for MDS with exact coverage}) For bounded martingale difference sequences with a deterministic predictable quadratic variation, Theorem \ref{th 1} gives an empirical Bernstein inequality whose radius only involves the observable $[\M]_n$, whose coverage probability is no less than $1-\alpha$ for every $n\geq 2$, and whose leading term is sharp in the sense of \eqref{eq 5}, without identically distributed observations or a convergent average variance. Corollary \ref{corollary 1} handles an unknown mean with the sample variance, with the same coverage and sharpness. Proposition \ref{prop 1} shows that the deterministic predictable quadratic variation cannot be dropped, and Theorems \ref{th 2} and \ref{th 3} relax it to a variance band and to moment conditions, respectively, which serve as templates for dependent data.

     \item [C2] (\textbf{Self-normalized inequalities for the marginal mean beyond strong mixing}) We focus on self-normalized concentration inequalities for the marginal mean in the presence of temporal dependence, rather than on the predictable (drifting) mean. In general mixing models, the gap between the two kinds of means is hard to control. Theorem \ref{th 4} covers $\tilde{\phi}$-mixing processes and functions of bounded variation, Theorem \ref{th 5} covers $\tau$-mixing processes and unbounded Lipschitz functions with moments of order $q>4$, and Theorem \ref{th 6} covers $\tau$-mixing processes and bounded Lipschitz functions. These classes include processes that are not strongly mixing, such as the autoregression of Example \ref{example 1} (Proposition \ref{prop B.1}), and the functions may change with the sample size. Part (i) of Theorems \ref{th 4}, \ref{th 5} and \ref{th 6} requires neither stationarity nor any assumption on the decay rate of the dependence coefficients. Under summability, their leading terms only involve a block sample variance and are sharp in the sense of \eqref{eq 5}, with the long-run variance in place of the marginal variance (part (ii) of Theorems \ref{th 4}, \ref{th 5} and \ref{th 6}), and the radii of Theorems \ref{th 5} and \ref{th 6} contain no dependence coefficient.

     \item [C3] (\textbf{Direct evaluation of time-series forecasting rules}) The losses of a forecasting rule along a time series are autocorrelated, so that EB inequalities for independent observations or martingale differences do not apply to them. Theorem \ref{th 7} and Corollary \ref{corollary 2} apply Theorem \ref{th 6} to these losses and give a generalization bound for the rule: its risk is bounded from above by the empirical risk plus a radius that is computed from the observed losses with a single block length, contains no dependence coefficient, and is sharp in the sense of \eqref{eq 5}, with the long-run variance in place of the marginal variance. In the simulation of Section \ref{sec application 2}, with sample sizes between $3000$ and $10000$, both the bound and its leading term alone hold in more than $95\%$ of the samples in all four designs, and the leading term approaches its oracle counterpart with known long-run variance (Figure \ref{fig sim}).
 \end{itemize}

\section{Self-Normalized Inequalities With Sample Variance Only}
\label{sec 3}
This section presents our sharp self-normalized inequalities under different dependence and moment conditions. Section \ref{sec 3.1} focuses on martingale difference sequences and Section \ref{sec 3.2} focuses on $\tilde{\phi}$- and $\tau$-mixing processes.
\subsection{Martingale Difference Sequence (MDS)}
\label{sec 3.1}
Given a probability space $(\Omega,\mathcal{F},\P)$, let $(\mathcal{A}_i)_{i\geq 0}$ be a filtration satisfying $\mathcal{A}_i\subset\mathcal{F}$. For a martingale difference sequence $(Z_i,\mathcal{A}_i)_{i\geq 1}$, we denote
\begin{align}
\label{eq notation MDS}
    \M_n=\sum_{i=1}^nZ_i,\quad \<\M\>_n=\sum_{i=1}^n\E[Z_i^2|\mathcal{A}_{i-1}],\quad [\M]_n=\sum_{i=1}^nZ_i^2.
\end{align}
Since an MDS has mean zero, these quantities coincide with those introduced in Section \ref{sec 1.1}. Note that $[\M]_n$ is completely observable, whereas $\<\M\>_n$ is not. The following theorem is the basic building block of the paper.

\begin{theorem}
\label{th 1}
Given a probability space $(\Omega,\mathcal{F},\P)$, let $(\mathcal{A}_i)_{i\geq 0}$ be a filtration satisfying $\mathcal{A}_i\subset\mathcal{F}$. Suppose $(Z_i,\mathcal{A}_i)_{i\geq 1}$ is a martingale difference sequence such that $|Z_i|\leq b$ holds for every $i\geq 1$, and let $\M_n$, $\<\M\>_n$ and $[\M]_n$ be defined as in \eqref{eq notation MDS}.
\begin{enumerate}[(i)]
    \item Suppose that, for the given $n\geq 1$, $\<\M\>_n$ is almost surely equal to a deterministic constant. Then, for any $t_1,t_2>0$,
    \begin{align}
    \label{eq th 1 1}
        \P\Big(|\M_n|\leq \sqrt{2t_1[\M]_n+b^2t_1t_2}+b\sqrt{t_1t_2}+\frac{bt_1}{3}\Big)\geq 1-2e^{-t_1}-e^{-t_2},
    \end{align}
    and the one-sided event $\big\{\M_n\leq \sqrt{2t_1[\M]_n+b^2t_1t_2}+b\sqrt{t_1t_2}+\frac{bt_1}{3}\big\}$ has probability no less than $1-e^{-t_1}-e^{-t_2}$. In particular, for any $\alpha\in(0,1)$ and $n\geq 2$, the choice
    \begin{align}
    \label{eq th 1 2}
        t_1=\log\Big(\frac{2\sqrt{n}}{(\sqrt{n}-1)\alpha}\Big),\qquad t_2=\log\Big(\frac{\sqrt{n}}{\alpha}\Big)
    \end{align}
    makes the right-hand side of \eqref{eq th 1 1} equal to $1-\alpha$.
    \item Suppose that the assumption of (i) holds for every $n\geq 2$. Denote $\sigma_n^2=\frac{1}{n}\sum_{i=1}^n\text{Var}(Z_i)$ and let
    \begin{align*}
        L_n=\frac{1}{n}\Big(\sqrt{2t_1[\M]_n+b^2t_1t_2}+b\sqrt{t_1t_2}+\frac{bt_1}{3}\Big)
    \end{align*}
    with $t_1$ and $t_2$ chosen as in \eqref{eq th 1 2}. If $\liminf_{n}\sigma_n^2>0$, then
    \begin{align}
    \label{eq th 1 3}
        \frac{\sqrt{n}L_n}{\sigma_n\sqrt{2\log(2/\alpha)}}\xrightarrow[]{\P}1.
    \end{align}
\end{enumerate}
\end{theorem}

Apart from bounded increments, the only structural requirement of Theorem \ref{th 1} is that the predictable quadratic variation $\<\M\>_n$ is deterministic. This holds whenever $Z_1,\dots,Z_n$ are independent, not necessarily identically distributed. More generally, it holds whenever every conditional variance $\E[Z_i^2|\mathcal{A}_{i-1}]$ is deterministic, while the remaining features of the conditional law of $Z_i$, such as its conditional skewness or its support, may depend on the past. A simple bounded example of the latter kind is the following. Let $(a_i)_{i\geq 1}$ be predictable with values in $[a_{\min},a_{\max}]\subset(0,\infty)$ and, conditionally on $\mathcal{A}_{i-1}$, let $Z_i$ take the value $a_i$ with probability $1/(1+a_i^2)$ and the value $-1/a_i$ otherwise. Then $\E[Z_i|\mathcal{A}_{i-1}]=0$, $\E[Z_i^2|\mathcal{A}_{i-1}]=1$ and $|Z_i|\leq \max\{a_{\max},1/a_{\min}\}$, whereas $(Z_i)_{i\geq 1}$ is in general not a sequence of independent random variables.

 The deterministic character of $\<\M\>_n$ is used twice in the proof, both for $\M_n$ and for the lower tail of $[\M]_n$: in both places, it allows the Chernoff parameter to be fixed before the data are observed. This requirement cannot be dropped: Proposition \ref{prop 1} below shows that, once $\<\M\>_n$ is allowed to be random, the bound in \eqref{eq th 1 1} fails with probability tending to one, both for fixed $t_1,t_2$ and for the choice \eqref{eq th 1 2}. Part (ii) shows that the leading term of \eqref{eq th 1 1} is sharp in the two-sided sense of \eqref{eq 5}: the $\sqrt{n}$-scaled radius behaves like $\sigma_n\sqrt{2\log(2/\alpha)}$, while the coverage probability is no less than $1-\alpha$ for every $n\geq 2$. In particular, no inflation of the confidence parameter survives in the limit. For heteroscedastic sequences, the relevant scale is the average variance $\sigma_n^2$, and no convergence of the conditional variances needs to be assumed, because $[\M]_n/n-\sigma_n^2=O_{\P}(n^{-1/2})$ holds automatically under the assumptions of Theorem \ref{th 1} (see Step 5 of its proof). The price of exact coverage is the lower-order term $b\sqrt{t_1t_2}=O(b\sqrt{\log n})$.

 Three elementary ingredients deliver the constants in \eqref{eq th 1 1}. First, the error probability is split unequally: the event controlling the variance receives only $\alpha/\sqrt{n}$, so that the price of estimating the variance moves into the lower-order term. Splitting the error probability equally, as in \citeA{maurer2009empirical}, inflates the leading constant permanently. Unequal splitting is in the spirit of the plug-in construction of \citeA{wang2024sharp} for independent observations. Any allocation of $o(1)$ probability to the variance event preserves sharpness, and the choice $\alpha/\sqrt{n}$ in \eqref{eq th 1 2} balances the first-order loss in $t_1$ against the growth of $t_2$. Second, the lower tail of $[\M]_n$ is controlled by the elementary inequality $e^{-x}\leq 1-x+\frac{x^2}{2}$, $x\geq 0$, which yields a sub-Gaussian bound without a linear term (Lemma \ref{lemma A.2}). A Bernstein-type bound would add a term of order $b^2t_2$ inside the square root. Third, the quadratic inequality relating $\sqrt{\<\M\>_n}$ to $[\M]_n$ is solved exactly instead of splitting the square root, which asymptotically halves the coefficient of the cross term $b\sqrt{t_1t_2}$.

\begin{proposition}
\label{prop 1}
Let $b>0$ and let $(t_{1,n})_{n\geq 1}$ and $(t_{2,n})_{n\geq 1}$ be positive sequences such that $\sup_{n}t_{1,n}<\infty$ and $t_{2,n}/n\to 0$. Then there exist a filtration $(\mathcal{A}_i)_{i\geq 0}$ and, for every $n\geq 1$, a martingale difference sequence $(Z_{n,i},\mathcal{A}_i)_{i\geq 1}$ such that $|Z_{n,i}|\leq b$ for every $i\geq 1$ and, with $\M_n$, $[\M]_n$ and $\<\M\>_n$ formed from $Z_{n,1},\dots,Z_{n,n}$, $\<\M\>_n=[\M]_n$ almost surely and
\begin{align}
\label{eq prop 1 1}
    \lim_{n\to\infty}\P\Big(\M_n>\sqrt{2t_{1,n}[\M]_n+b^2t_{1,n}t_{2,n}}+b\sqrt{t_{1,n}t_{2,n}}+\frac{bt_{1,n}}{3}\Big)=1.
\end{align}
If $t_{1,n}$ and $t_{2,n}$ do not depend on $n$, the martingale difference sequence can be chosen not to depend on $n$ either.
\end{proposition}

Proposition \ref{prop 1} covers both fixed $t_1,t_2>0$ and the choice \eqref{eq th 1 2}, for which $\sup_{n\geq 2}t_{1,n}\leq\log\big(\frac{2\sqrt{2}}{(\sqrt{2}-1)\alpha}\big)$ and $t_{2,n}/n=\log(\sqrt{n}/\alpha)/n\to 0$. In both cases, the probability in \eqref{eq th 1 1} tends to zero for the MDS of Proposition \ref{prop 1}, although $|Z_i|\leq b$. Hence, the assumption that $\<\M\>_n$ is deterministic cannot be removed from Theorem \ref{th 1}. Since $\<\M\>_n=[\M]_n$ in the construction, and since the radius in \eqref{eq prop 1 1} is no smaller than $\sqrt{2t_{1,n}\<\M\>_n}+bt_{1,n}/3$, the same conclusion also holds for the oracle bound in which $[\M]_n$ is replaced by $\<\M\>_n$: for bounded MDS with random $\<\M\>_n$, a fixed-$n$ bound of the form $\P\big(\M_n\leq\sqrt{2\<\M\>_nt}+bt/3\big)\geq 1-e^{-t}$ cannot hold in general.

The construction is a simple random walk stopped at its first crossing of the radius in \eqref{eq prop 1 1}, and the proof only uses that such a walk exceeds $K\sqrt{k}$ somewhere on a long time window with probability tending to one (Lemma \ref{lemma A.4}). This is the mechanism behind the law of the iterated logarithm, which also underlies the iterated-logarithm overhead of bounds that remain valid for random $\<\M\>_n$, such as the time-uniform bounds of \citeA{howard2021time}.

One focus of this paper is to derive sharp self-normalized concentration inequalities for commonly used classes of processes that are not strongly mixing, such as $\tilde{\phi}$-mixing and $\tau$-mixing processes (see Section \ref{sec 3.2}). For such processes, the blocking construction of Section \ref{sec 3.2} produces martingale differences whose predictable quadratic variation is random, but nearly deterministic. Theorems \ref{th 2} and \ref{th 3} below relax the deterministic $\<\M\>_n$ of Theorem \ref{th 1} accordingly, and can be regarded as the template theorems for $\tilde{\phi}$-mixing and $\tau$-mixing, respectively. 

\begin{theorem}
\label{th 2}
Let $(Z_i,\mathcal{A}_i)_{i\geq 1}$ be a martingale difference sequence such that $|Z_i|\leq b$ holds for every $i\geq 1$, and let $\M_n$, $\<\M\>_n$ and $[\M]_n$ be defined as in \eqref{eq notation MDS}. For the given $n\geq 1$, let $v>0$ and $\Gamma\geq 0$ be deterministic. Then, for any $t_1,t_2>0$,
\begin{align}
    &\P\Big(|\M_n|\leq \sqrt{2t_1\big([\M]_n+\Gamma\big)+b^2t_1t_2}+b\sqrt{t_1t_2}+\frac{bt_1}{3}\Big)\notag\\
    &\qquad\geq 1-2e^{-t_1}-e^{-t_2}-\P\big(\<\M\>_n\notin[v-\Gamma,v]\big),\label{eq th 2 1}
\end{align}
and the corresponding one-sided event has probability no less than $1-e^{-t_1}-e^{-t_2}-\P\big(\<\M\>_n\notin[v-\Gamma,v]\big)$.
\end{theorem}

Theorem \ref{th 2} is the template for the extension to $\tilde{\phi}$-mixing in Section \ref{sec 3.2}. There, the $\tilde{\phi}$-coefficients control conditional expectations in $L^\infty$, so that the blocking construction produces martingale differences that are bounded almost surely, and Theorem \ref{th 2} reduces the dependent-data inequality to a single input: a variance band $[v-\Gamma,v]$ for $\<\M\>_n$ with $\Gamma=o(v)$, derived from the $\tilde{\phi}$-coefficients of the process. Since this control is also in $L^\infty$, the band holds almost surely and the term $\P(\<\M\>_n\notin[v-\Gamma,v])$ in \eqref{eq th 2 1} vanishes. Once the band is available, the radius, its leading term $\sqrt{2t_1[\M]_n}$ and the coverage probability are inherited from \eqref{eq th 2 1}, without any further exponential argument. The level $v$ itself may be unknown, since it does not enter the radius.

 Only an upper bound of $\<\M\>_n$ can be fed into a Chernoff bound, and only a lower bound of $\<\M\>_n$ is needed to replace it by $[\M]_n$. $v$ is the former and $v-\Gamma$ is the latter. Indeed, both exponential super-martingales in the proof control $\M_n$ and $\<\M\>_n-[\M]_n$ only after $\<\M\>_n$ has been replaced by a deterministic number no smaller than $\<\M\>_n$, and the Chernoff parameters are tuned at that number. The lower end of the band is used once, to turn the unknown $v$ into the observable $[\M]_n+\Gamma$. Centring the band brings nothing: a symmetric band $|\<\M\>_n-v'|\leq \Gamma'$ is the case $v=v'+\Gamma'$, $\Gamma=2\Gamma'$, and unequal deviations $\Gamma_-$ below and $\Gamma_+$ above $v'$ correspond to $v=v'+\Gamma_+$, $\Gamma=\Gamma_-+\Gamma_+$.

\begin{theorem}
\label{th 3}
Let $(Z_i,\mathcal{A}_i)_{i\geq 1}$ be a martingale difference sequence such that $\E[Z_i^4]<\infty$ holds for every $i\geq 1$, and let $\M_n$, $\<\M\>_n$ and $[\M]_n$ be defined as in \eqref{eq notation MDS}. Denote $\sigma_n^2=\frac{1}{n}\sum_{i=1}^n\text{Var}(Z_i)$, so that $n\sigma_n^2=\E[\<\M\>_n]=\E[[\M]_n]$. Then, for any $t>0$ and any deterministic $x>0$,
\begin{align}
\label{eq th 3 1}
    \P\Big(|\M_n|\leq\sqrt{2t\big([\M]_n+2x\big)}\Big)\geq 1-2e^{-t}-\frac{2}{x}\E\big|\<\M\>_n-n\sigma_n^2\big|-\frac{4}{x^2}\sum_{i=1}^n\E[Z_i^4],
\end{align}
and the corresponding one-sided event has probability no less than the same lower bound with $2e^{-t}$ replaced by $e^{-t}$.
\end{theorem}

 Theorem \ref{th 3} is the template for the extension to $\tau$-mixing in Section \ref{sec 3.2}. There, the martingale differences obtained from the blocking construction are bounded only in $L^p$, so that the Bernstein-type engine of Theorems \ref{th 1} and \ref{th 2} is not available. Theorem \ref{th 3} replaces it by a self-normalized super-martingale that requires square integrability only (Lemma \ref{lemma A.5}, the condition ``SN I'' in Table 3 of \citeA{howard2020time}), so that neither bounded increments nor truncation are needed. The dependent-data inequality is then reduced to two moment-level inputs: the coupling error $\E|\<\M\>_n-n\sigma_n^2|$ between the conditional and the marginal variance, which is controlled by the dependence coefficients of the process, and the fourth moments of the martingale differences. Once these two inputs are available, the radius, its leading term and the coverage probability are inherited from \eqref{eq th 3 1}.

The radius in \eqref{eq th 3 1} contains no linear term, and its leading term is $\sqrt{2t[\M]_n}$ as soon as $x=o(n\sigma_n^2)$. The coverage probability tends to $1-2e^{-t}$ if, in addition, $\E|\<\M\>_n-n\sigma_n^2|=o(x)$ and $\sum_{i=1}^n\E[Z_i^4]=o(x^2)$. For instance, if $\sum_{i=1}^n\E[Z_i^4]=O(n)$ and $\E|\<\M\>_n-n\sigma_n^2|=O(\sqrt{n})$, any $x$ with $\sqrt{n}=o(x)$ and $x=o(n)$ works. Markov's inequality with higher moments of $\<\M\>_n-n\sigma_n^2$ can be used in place of the first moment. The price of dropping bounded increments is that the cost of controlling $[\M]_n$ is moment-level rather than exponential. This is why Theorems \ref{th 1} and \ref{th 2} keep the Bernstein-type engine when the increments are bounded.

A setting closely related to Theorem \ref{th 1} is that $\E[Z_i|\mathcal{A}_{i-1}]=\mu$ holds for every $i$, which makes $(Z_i-\mu,\mathcal{A}_i)_{i\geq 1}$ an MDS. When $\mu$ is unknown, $\sum_{i=1}^n(Z_i-\mu)^2$ is no longer observable. Replacing $\mu$ by the sample mean yields the following empirical Bernstein inequality.

\begin{corollary}
\label{corollary 1}
Given a probability space $(\Omega,\mathcal{F},\P)$, let $(\mathcal{A}_i)_{i\geq 0}$ be a filtration satisfying $\mathcal{A}_i\subset\mathcal{F}$. Suppose $(Z_i)_{i\geq 1}$ is adapted to $(\mathcal{A}_i)_{i\geq 0}$ and satisfies $\E[Z_i|\mathcal{A}_{i-1}]=\mu$ and $|Z_i-\mu|\leq b$ for every $i\geq 1$. Denote $\overline{\mathbf{M}}_n=\frac{1}{n}\sum_{i=1}^nZ_i$ and $[\mathbb{V}]_n=\sum_{i=1}^n(Z_i-\overline{\mathbf{M}}_n)^2$.
\begin{enumerate}[(i)]
    \item Suppose that, for the given $n\geq 1$, $\sum_{i=1}^n\E[(Z_i-\mu)^2|\mathcal{A}_{i-1}]$ is almost surely equal to a deterministic constant. Then, for any $t_1,t_2>0$ such that $2t_1<n$,
    \begin{align}
    \label{eq corollary 1 1}
    \begin{split}
        &\P\Big(\big|\overline{\mathbf{M}}_n-\mu\big|\leq L_n\Big)\geq 1-2e^{-t_1}-e^{-t_2},\\
        &L_n=\frac{1}{n-2t_1}\Big(\sqrt{2t_1[\mathbb{V}]_n+b^2t_1t_2}+\Big(1+\sqrt{\frac{2t_1}{n}}\Big)\Big(b\sqrt{t_1t_2}+\frac{bt_1}{3}\Big)\Big),
    \end{split}
    \end{align}
    and each of the one-sided events $\{\overline{\mathbf{M}}_n-\mu\leq L_n\}$ and $\{\mu-\overline{\mathbf{M}}_n\leq L_n\}$ has probability no less than $1-e^{-t_1}-e^{-t_2}$. In particular, for any $\alpha\in(0,1)$ and $n\geq 2$, the choice \eqref{eq th 1 2} makes the right-hand side of \eqref{eq corollary 1 1} equal to $1-\alpha$ whenever $2t_1<n$.
    \item Suppose that the assumption of (i) holds for every $n\geq 2$. Denote $\sigma_n^2=\frac{1}{n}\sum_{i=1}^n\text{Var}(Z_i)$ and let $L_n$ be the radius in \eqref{eq corollary 1 1} with $t_1$ and $t_2$ chosen as in \eqref{eq th 1 2}. If $\liminf_n\sigma_n^2>0$, then
    \begin{align*}
        \frac{\sqrt{n}L_n}{\sigma_n\sqrt{2\log(2/\alpha)}}\xrightarrow[]{\P}1.
    \end{align*}
\end{enumerate}
\end{corollary}

Corollary \ref{corollary 1} inherits from Theorem \ref{th 1} the exact coverage and the sharp leading term. Replacing $\mu$ by $\overline{\mathbf{M}}_n$ creates the term $n(\overline{\mathbf{M}}_n-\mu)^2$ in the identity $\sum_{i=1}^n(Z_i-\mu)^2=[\mathbb{V}]_n+n(\overline{\mathbf{M}}_n-\mu)^2$. Splitting the square root to isolate this term would inflate the whole radius by the factor $(1-\sqrt{2t_1/n})^{-1}=1+O(n^{-1/2})$. Instead, the resulting quadratic inequality in $|\overline{\mathbf{M}}_n-\mu|$ is solved exactly, which yields the factor $n/(n-2t_1)=1+O(n^{-1})$ in front of the leading term. The boundedness condition concerns the centred variables: if $Z_i$ takes values in a known interval of length $b$, then $|Z_i-\mu|\leq b$. If only $|Z_i|\leq b'$ is known, one may take $b=2b'$. The deterministic condition holds, for example, for independent bounded observations with a common mean and possibly different variances.


\subsection{$\tilde{\phi}$- and $\tau$-Mixing}
\label{sec 3.2}

\subsubsection{Review of Dependence Conditions}
\label{sec 3.2.1}
Throughout Section \ref{sec 3.2}, $\mathbf{Z}=(Z_i)_{i\in\mathbb{Z}}$ is a real-valued process and, for $-\infty\leq a\leq b\leq\infty$, $\mathcal{F}_a^b$ denotes the $\sigma$-algebra generated by $Z_i$, $a\leq i\leq b$. For a $\sigma$-algebra $\mathcal{M}\subset\mathcal{F}$ and a random vector $X=(X_1,\dots,X_l)$ with values in $\mathbb{R}^l$, $\P_{X|\mathcal{M}}$ denotes a regular conditional distribution of $X$ given $\mathcal{M}$, and $\{X\leq x\}=\bigcap_{m=1}^l\{X_m\leq x_m\}$ for $x=(x_1,\dots,x_l)\in\mathbb{R}^l$. We use the following coefficients. Definition \ref{def tau k} recalls the coefficient $\tau$ of \citeA{dedecker2005new}, its bivariate version $\tau_2$ and its version $\tau_\infty$ for finite-dimensional random vectors. Definition \ref{def phi tilde k} recalls the coefficient $\tilde{\phi}$ of \citeA{dedecker2005new} and introduces a bivariate version $\tilde{\phi}_2$, which is discussed after the definitions.

\begin{definition}
\label{def phi tilde k}
For a $\sigma$-algebra $\mathcal{M}$ and a random vector $X$ with values in $\mathbb{R}^l$, $l\in\{1,2\}$, let
\begin{align*}
    \tilde{\phi}(\mathcal{M},X)=\sup_{x\in\mathbb{R}^l}\big\|\P(X\leq x|\mathcal{M})-\P(X\leq x)\big\|_{\infty}.
\end{align*}
For $r\geq 1$, define
\begin{align}
\label{eq def phi tilde k}
    \tilde{\phi}_1(r)&=\sup\Big\{\tilde{\phi}\big(\mathcal{F}_{-\infty}^p,Z_j\big):\ p\in\mathbb{Z},\ j\geq p+r\Big\},\notag\\
    \tilde{\phi}_2(r)&=\max\Big\{\tilde{\phi}_1(r),\ \sup\Big\{\tilde{\phi}\big(\mathcal{F}_{-\infty}^p,(Z_i,Z_j)\big):\ p\in\mathbb{Z},\ p+r\leq i<j\Big\}\Big\}.
\end{align}
We say that $\mathbf{Z}$ is $\tilde{\phi}$-mixing if $\tilde{\phi}_1(r)\to 0$ as $r\to\infty$.
\end{definition}

\begin{definition}
\label{def tau k}
Equip $\mathbb{R}^l$, $l\geq 1$, with the distance $\|x-y\|_1=\sum_{m=1}^l|x_m-y_m|$ and let $\Lambda_1(\mathbb{R}^l)$ be the set of $1$-Lipschitz functions $f:\mathbb{R}^l\to\mathbb{R}$. For a $\sigma$-algebra $\mathcal{M}$ and a random vector $X$ with values in $\mathbb{R}^l$ such that $\E\|X\|_1<\infty$, let
\begin{align*}
    \tau(\mathcal{M},X)=\Big\|\sup_{f\in\Lambda_1(\mathbb{R}^l)}\Big|\int f\,d\P_{X|\mathcal{M}}-\int f\,d\P_X\Big|\Big\|_1.
\end{align*}
For $r\geq 1$, define
\begin{align}
\label{eq def tau k}
    \tau_1(r)&=\sup\Big\{\tau\big(\mathcal{F}_{-\infty}^p,Z_j\big):\ p\in\mathbb{Z},\ j\geq p+r\Big\},\notag\\
    \tau_2(r)&=\max\Big\{\tau_1(r),\ \frac{1}{2}\sup\Big\{\tau\big(\mathcal{F}_{-\infty}^p,(Z_i,Z_j)\big):\ p\in\mathbb{Z},\ p+r\leq i<j\Big\}\Big\},\notag\\
    \tau_\infty(r)&=\sup_{l\geq 1}\frac{1}{l}\sup\Big\{\tau\big(\mathcal{F}_{-\infty}^p,(Z_{j_1},\dots,Z_{j_l})\big):\ p\in\mathbb{Z},\ p+r\leq j_1<\dots<j_l\Big\}.
\end{align}
Clearly, $\tau_1\leq\tau_2\leq\tau_\infty$, and the three coefficients are non-increasing in $r$. We say that $\mathbf{Z}$ is $\tau$-mixing if $\E|Z_i|<\infty$ for every $i$ and $\tau_1(r)\to 0$ as $r\to\infty$.
\end{definition}

 For real-valued $X$, $\tau(\mathcal{M},X)=\int_{\mathbb{R}}\|\P(X\leq t|\mathcal{M})-\P(X\leq t)\|_1dt$, and $\tau_2$ and $\tau_\infty$ are the case $k=2$ of the coefficient $\tau_k$ and the coefficient $\tau_\infty=\sup_k\tau_k$ introduced in Section 7 of \citeA{dedecker2005new}, who call the corresponding sequences $\tau$-dependent. The coefficient $\tilde{\phi}(\mathcal{M},X)$ for real-valued $X$ is also taken from \citeA{dedecker2005new}. The bivariate coefficient $\tilde{\phi}_2$, however, is not taken from the literature. \citeA{dedecker2005new} point out in their Section 7 that coefficients based on conditional distribution functions are difficult to extend to higher dimensions, because the notion of bounded variation is already delicate in $\mathbb{R}^2$, and they extend $\tau$ instead. A different multivariate extension, based on products of centred indicator functions, is used by \citeA{dedecker2007empirical}. Our definition of $\tilde{\phi}_2$ is dictated by the proof technique of Section \ref{sec 3.2.2}: after blocking and the exact coboundary representation, the conditional second moments of block sums involve products $h_n(Z_i)h_n(Z_{i'})$, and the variance band required by Theorem \ref{th 2} needs a control of their conditional expectations that is uniform in $\omega$. The coefficient $\tilde{\phi}_2$ provides this control for functions of bounded variation (Lemma \ref{lemma A.6}), and it only involves lower quadrants of $\mathbb{R}^2$. Note that $\tilde{\phi}$ controls conditional expectations uniformly in $\omega$, whereas $\tau$ controls them in $L^1$ only.

 The bivariate extension preserves the generality of $\tilde{\phi}$. In particular, it still contains processes that are not strongly mixing.

\begin{example}
\label{example 1}
Let $(\varepsilon_t)_{t\in\mathbb{Z}}$ be i.i.d.\ with $\P(\varepsilon_t=0)=\P(\varepsilon_t=1)=\frac{1}{2}$ and let $Z_t=\frac{1}{2}(Z_{t-1}+\varepsilon_t)$, whose stationary solution $Z_t=\sum_{j\geq 0}2^{-(j+1)}\varepsilon_{t-j}$ is uniformly distributed on $[0,1]$. Since the binary digits of $Z_t$ are $\varepsilon_t,\varepsilon_{t-1},\dots$, we have $Z_p=2^rZ_{p+r}-\lfloor 2^rZ_{p+r}\rfloor$ almost surely. Hence the event $\{Z_p\leq\frac{1}{2}\}$ belongs to both $\mathcal{F}_{-\infty}^p$ and $\mathcal{F}_{p+r}^{\infty}$, and the strong mixing coefficient of $\mathbf{Z}$ is at least $\frac{1}{4}$ for every $r$, that is, $\mathbf{Z}$ is not strongly mixing (see \citeA{andrews1984non}, and Example 1 in Section 4.1 of \citeA{dedecker2005new}). Nevertheless, $\tilde{\phi}_1(r)\leq 2^{-r}$, as shown in the same example of \citeA{dedecker2005new}, and $\tilde{\phi}_2(r)\leq 2^{1-r}$ for every $r\geq 1$ (Proposition \ref{prop B.1} (iv) in Appendix \ref{app B}).
\end{example}

Examples of $\tau$-mixing processes are given in Sections 4 and 7.2 of \citeA{dedecker2005new}, including causal linear processes, iterated random functions with contracting maps, and nonlinear autoregressions $Z_t=g(Z_{t-1})+\varepsilon_t$ with i.i.d.\ integrable innovations and $g$ Lipschitz with constant smaller than one, for which $\tau_\infty$, and hence $\tau_2$, decays geometrically.

The two coefficients do not control each other (Proposition \ref{prop B.1} in Appendix \ref{app B}). On the one hand, the Gaussian autoregression $Z_t=aZ_{t-1}+\varepsilon_t$ with $0<|a|<1$ is $\tau$-mixing by the results just cited, but it is not $\tilde{\phi}$-mixing, because $\tilde{\phi}$ requires a control that is uniform in $\omega$, which fails on an unbounded state space. On the other hand, $\tilde{\phi}_1$ and $\tilde{\phi}_2$ do not increase under strictly increasing transformations of the process and require no moment condition, whereas $\tau$ is only defined for integrable variables. For instance, a strictly increasing transformation of Example \ref{example 1} has Cauchy marginals and the same bounds on $\tilde{\phi}_1$ and $\tilde{\phi}_2$. Moreover, even for bounded processes, $\tilde{\phi}_2$ controls $\tau_2$ only through a square root, and no linear bound is possible. The two coefficients also serve different purposes. The coefficient $\tilde{\phi}_2$ controls conditional expectations of products of functions of bounded variation, including indicator functions, uniformly in $\omega$ and without moment conditions. This leads to Theorem \ref{th 4}, whose coverage probability is of exponential type and which applies, for instance, to empirical distribution functions and to $0$-$1$ losses. The coefficient $\tau_2$ controls Lipschitz functions of possibly unbounded processes, but only in $L^1$, which leads to Theorem \ref{th 5}, whose coverage probability is of moment type. For bounded Lipschitz functions, the coefficient $\tau_\infty$ allows a coupling argument instead, which leads to Theorem \ref{th 6}. We therefore treat the two classes separately.

\subsubsection{Inequalities}
\label{sec 3.2.2}
Throughout this subsection, for the given $n\geq 1$, $h_n:\mathbb{R}\to\mathbb{R}$ is a measurable function such that $\E|h_n(Z_i)|<\infty$ and $\E[h_n(Z_i)]=\mu_n$ hold for every $1\leq i\leq n$, and $\P_n(h_n)=\frac{1}{n}\sum_{i=1}^nh_n(Z_i)$. For a deterministic block length $l\geq 1$, we assume that $n=ml$ for an integer $m\geq 1$ and partition $\{1,\dots,n\}$ into the blocks $B_j=\{(j-1)l+1,\dots,jl\}$, $1\leq j\leq m$. If $l$ does not divide $n$, the last $n-l\lfloor n/l\rfloor$ observations are discarded. The block empirical long-run variance is
\begin{align}
\label{eq V hat}
    \hat{V}_n=\frac{1}{n}\sum_{j=1}^m\Big(\sum_{i\in B_j}\big(h_n(Z_i)-\P_n(h_n)\big)\Big)^2.
\end{align}
\paragraph{$\tilde{\phi}$-mixing.} For $\tilde{\phi}$-mixing processes, $h_n$ takes values in an interval $[a_n,b_n]$ and is of bounded variation with total variation $\|h_n\|_{TV}$. The dependence of $\mathbf{Z}$ then enters only through
\begin{align}
\label{eq Phi Psi}
    \tilde{\Phi}_n=\sum_{r=1}^n\tilde{\phi}_1(r),\qquad \Psi_n=\sum_{r=1}^nr\,\tilde{\phi}_2(r),
\end{align}
and we denote
\begin{align}
\label{eq b bar D}
    \bar{b}_n&=l(b_n-a_n)+2\|h_n\|_{TV}\tilde{\Phi}_n,\notag\\
    \Gamma_n&=2m\Big[\big(4\|h_n\|_{TV}^2+10(b_n-a_n)\|h_n\|_{TV}\big)\Psi_n+2\|h_n\|_{TV}^2\tilde{\Phi}_n^2\Big].
\end{align}

\begin{theorem}
\label{th 4}
Let $l$ and $m$ be as above, let $h_n:\mathbb{R}\to[a_n,b_n]$ be of bounded variation, and let $\hat{V}_n$, $\tilde{\Phi}_n$, $\Psi_n$, $\bar{b}_n$ and $\Gamma_n$ be given in \eqref{eq V hat}--\eqref{eq b bar D}.
\begin{enumerate}[(i)]
    \item For any $t_1,t_2>0$ such that $2t_1<m$,
    \begin{align}
    \label{eq th 4 1}
        \P\Big(\big|\P_n(h_n)-\mu_n\big|\leq L_n\Big)\geq 1-2e^{-t_1}-e^{-t_2},
    \end{align}
    where
    \begin{align}
        L_n=\frac{1}{n-2t_1l}\bigg(&\sqrt{2t_1n\hat{V}_n+2t_1\Gamma_n+\bar{b}_n^2t_1t_2}\notag\\
        &+\Big(1+\sqrt{\frac{2t_1}{m}}\Big)\Big(\bar{b}_n\sqrt{t_1t_2}+\frac{\bar{b}_nt_1}{3}+2\sqrt{2t_1m}\,\|h_n\|_{TV}\tilde{\Phi}_n\Big)\bigg),\label{eq th 4 2}
    \end{align}
    and each of the one-sided events $\{\P_n(h_n)-\mu_n\leq L_n\}$ and $\{\mu_n-\P_n(h_n)\leq L_n\}$ has probability no less than $1-e^{-t_1}-e^{-t_2}$. In particular, for any $\alpha\in(0,1)$, the choice \eqref{eq th 1 2} with $n$ replaced by $m$ makes the right-hand side of \eqref{eq th 4 1} equal to $1-\alpha$ whenever $2t_1<m$.
    \item Let $l=l_n$, with $m=n/l_n$ for the values of $n$ considered, and suppose that $l_n\to\infty$ and $l_n^2\log n/n\to 0$ as $n\to\infty$, that $\sup_n\max\{b_n-a_n,\|h_n\|_{TV}\}<\infty$ and that $\sum_{r=1}^\infty r\,\tilde{\phi}_2(r)<\infty$. Denote $\sigma_n^2=\frac{1}{n}\text{Var}\big(\sum_{i=1}^nh_n(Z_i)\big)$ and let $L_n$ be the radius in \eqref{eq th 4 2} with $t_1$ and $t_2$ chosen as in \eqref{eq th 1 2} with $n$ replaced by $m$. If $\liminf_n\sigma_n^2>0$, then
    \begin{align*}
        \frac{\sqrt{n}L_n}{\sigma_n\sqrt{2\log(2/\alpha)}}\xrightarrow[]{\P}1.
    \end{align*}
\end{enumerate}
\end{theorem}

Theorem \ref{th 4} follows from Theorem \ref{th 2} in three steps. First, an auxiliary triangular array, whose in-sample part has the distribution of $(Z_1,\dots,Z_n)$ and whose past and future are independent of it, turns the coboundary decomposition of $\sum_{i=1}^n(h_n(Z_i)-\mu_n)$ along the block ends into an exact representation as a sum of $m$ martingale differences: the boundary terms at both ends of the sample vanish, so that no remainder term is left (Lemma \ref{lemma A.7}). Each martingale difference is the centred block sum plus two boundary terms of size at most $\|h_n\|_{TV}\tilde{\Phi}_n$, so that it is bounded by $\bar{b}_n$. Second, $\tilde{\phi}_2$ controls the conditional second moments of block sums uniformly in $\omega$, so that the predictable quadratic variation lies almost surely in a band of width $\Gamma_n$ around $\text{Var}(\sum_{i=1}^nh_n(Z_i))$ (Lemma \ref{lemma A.8}), and the probability term in \eqref{eq th 2 1} vanishes. Third, the boundary terms are removed from the realized quadratic variation by the triangle inequality in $\mathbb{R}^m$, which produces the term $2\sqrt{2t_1m}\,\|h_n\|_{TV}\tilde{\Phi}_n$, and the unknown mean in the block sums is handled by solving the resulting quadratic inequality exactly, as in Corollary \ref{corollary 1}. This yields the factor $n/(n-2t_1l)=m/(m-2t_1)$.

Part (i) uses neither stationarity nor summability of the coefficients: besides bounded variation, it only requires the mean of $h_n(Z_i)$ to be constant in $1\leq i\leq n$. The coefficients enter the radius only through $\tilde{\Phi}_n$ and $\Psi_n$, and $L_n$ is non-decreasing in both, so that they may be replaced by any known upper bounds. In part (ii), $\sigma_n^2$ is the long-run variance at sample size $n$, which converges to the usual long-run variance under stationarity and summable autocovariances. Hence, the leading term $\sqrt{2t_1\hat{V}_n/n}$ is sharp in the two-sided sense of \eqref{eq 5}, with the long-run variance in place of the marginal variance. The condition $l_n\to\infty$ is needed because $\Gamma_n/n$ and $\sqrt{m/n}=l_n^{-1/2}$ must vanish, while $l_n^2\log n=o(n)$ makes the terms involving $\bar{b}_n=O(l_n)$ negligible.

\paragraph{$\tau$-mixing.} For $\tau$-mixing processes, $h_n$ is Lipschitz, with smallest Lipschitz constant denoted by $\|h_n\|_{\text{Lip}}$, and it need not be bounded. Since $\tau$ controls conditional expectations in $L^1$ only, the block martingale differences are no longer bounded, and their predictable quadratic variation is close to its mean in $L^1$ only. We therefore use Theorem \ref{th 3} in place of Theorem \ref{th 2}, and we bound its moment inputs through moments of order $q>4$. Let
\begin{align*}
    \kappa_n=\max_{1\leq i\leq n}\big\|h_n(Z_i)-\mu_n\big\|_q<\infty,\qquad \beta=\frac{q-2}{q-1},\qquad \gamma=\frac{q-4}{4(q-1)},
\end{align*}
and denote
\begin{align}
\label{eq tau quantities}
    \Psi_n^\tau=\sum_{r=1}^nr\,\tau_2(r)^\beta,\qquad \Phi_n^\tau=\sum_{r=1}^n\tau_1(r)^\gamma,\qquad c_n=\kappa_n^{q/(q-1)}\|h_n\|_{\text{Lip}}^\beta,
\end{align}
and, for $x>0$,
\begin{align}
\label{eq eta}
    \eta_n(x)=\frac{112\,m\,c_n\Psi_n^\tau}{x}+\frac{32\big(nl^3\kappa_n^4+16\,m\,\kappa_n^{4-4\gamma}\|h_n\|_{\text{Lip}}^{4\gamma}(\Phi_n^\tau)^4\big)}{x^2}.
\end{align}

\begin{theorem}
\label{th 5}
Let $\mathbf{Z}$ be $\tau$-mixing, let $l$ and $m$ be as above, let $q>4$, and let $h_n$ be Lipschitz with $\kappa_n<\infty$. Let $\hat{V}_n$ and $\eta_n$ be given in \eqref{eq V hat} and \eqref{eq eta}.
\begin{enumerate}[(i)]
    \item For any $t>0$ such that $2t<m$ and any deterministic $x>0$,
    \begin{align}
    \label{eq th 5 1}
    \begin{split}
        &\P\Big(\big|\P_n(h_n)-\mu_n\big|\leq L_n\Big)\geq 1-2e^{-t}-\eta_n(x),\\
        &L_n=\frac{\sqrt{2t\big(n\hat{V}_n+2x\big)}+\Big(1+\sqrt{\frac{2t}{m}}\Big)\sqrt{2tx}}{n-2tl},
    \end{split}
    \end{align}
    and each of the one-sided events $\{\P_n(h_n)-\mu_n\leq L_n\}$ and $\{\mu_n-\P_n(h_n)\leq L_n\}$ has probability no less than $1-e^{-t}-\eta_n(x)$. In particular, for any $\alpha\in(0,1)$, the choice $t=\log(2/\alpha)$ makes the right-hand side of \eqref{eq th 5 1} equal to $1-\alpha-\eta_n(x)$.
    \item Let $l=l_n$ and $x=x_n$ be deterministic, with $m=n/l_n$ for the values of $n$ considered, and suppose that, as $n\to\infty$, $l_n\to\infty$, $l_n^3/n\to 0$, $x_n/n\to 0$, $l_nx_n/n\to\infty$ and $x_n^2/(nl_n^3)\to\infty$. Suppose further that $\sup_n\|h_n\|_{\text{Lip}}<\infty$, $\sup_n\kappa_n<\infty$, $\sum_{r=1}^\infty r\,\tau_2(r)^\beta<\infty$ and $\sum_{r=1}^\infty\tau_1(r)^\gamma<\infty$. Then $\eta_n(x_n)\to 0$. Moreover, denote $\sigma_n^2=\frac{1}{n}\text{Var}\big(\sum_{i=1}^nh_n(Z_i)\big)$ and let $L_n$ be the radius in \eqref{eq th 5 1} with $t=\log(2/\alpha)$ and $x=x_n$. If $\liminf_n\sigma_n^2>0$, then
    \begin{align*}
        \frac{\sqrt{n}L_n}{\sigma_n\sqrt{2\log(2/\alpha)}}\xrightarrow[]{\P}1.
    \end{align*}
\end{enumerate}
\end{theorem}

 The blocking construction and the exact coboundary representation are the same as for Theorem \ref{th 4}. Theorem \ref{th 3} then reduces the problem to its two moment inputs, the coupling error between the predictable quadratic variation and $\text{Var}(\sum_{i=1}^nh_n(Z_i))$, and the fourth moments of the block martingale differences, which are bounded through $\tau_1$, $\tau_2$ and $\kappa_n$ (Lemma \ref{lemma A.9}). Since $h_n$ may be unbounded, products $h_n(Z_i)h_n(Z_{i'})$ are controlled by truncating $h_n$ at a level that balances $\tau_2$ against the moments of order $q$, which explains the exponent $\beta<1$. The boundary terms of the coboundary representation are controlled in $L^4$ by interpolation between $L^1$, where $\tau_1$ applies, and $L^q$. This explains the exponent $\gamma$ and the requirement $q>4$. Finally, the boundary terms are removed from the realized quadratic variation on an event whose probability is bounded by Markov's inequality, with the same slack $x$, and the unknown mean is handled by solving the resulting quadratic inequality exactly, as in Theorem \ref{th 4}.

The radius in Theorem \ref{th 5} contains neither dependence coefficients nor bounds on $h_n$: apart from the deterministic slack $x$, it is computed from the data. The dependence of $\mathbf{Z}$ enters only through $\eta_n(x)$, that is, through the coverage probability. Part (i) requires neither stationarity nor summability of the coefficients. In part (ii), $x_n=o(n)$ keeps the leading term $\sqrt{2t\hat{V}_n/n}$ sharp in the two-sided sense of \eqref{eq 5}, with the long-run variance in place of the marginal variance, while $l_nx_n/n\to\infty$ and $x_n^2/(nl_n^3)\to\infty$ make $\eta_n(x_n)$ vanish. Such $x_n$ exist whenever $l_n\to\infty$ and $l_n^3=o(n)$, for instance $l_n=n^{1/4}$ and $x_n=n/\log n$. Since $\eta_n(x)$ is of moment type, a small coverage loss requires a slack $x$ that is large compared with $m$ and with $\sqrt{nl^3}$. This is the price of dropping bounded increments, see the discussion after Theorem \ref{th 3}.

\paragraph{Bounded Lipschitz functions.} If $h_n$ is Lipschitz and bounded, the coverage loss of moment type in Theorem \ref{th 5} can be replaced by a coupling term that involves $\tau_\infty$ at a single lag. To this end, we fix an integer $d$ with $1\leq d\leq l/2$ and split each block $B_j$ into a big block $I_j$, formed by its first $l-d$ indices, and a small block $J_j$, formed by its last $d$ indices. Denote
\begin{align}
\label{eq U W}
    U_j=\sum_{i\in I_j}h_n(Z_i),\qquad W_j=\sum_{i\in J_j}h_n(Z_i),\qquad 1\leq j\leq m,
\end{align}
and let $[\mathbb{V}]_m^U=\sum_{j=1}^m(U_j-\frac{1}{m}\sum_{k=1}^mU_k)^2$ and $[\mathbb{V}]_m^W=\sum_{j=1}^m(W_j-\frac{1}{m}\sum_{k=1}^mW_k)^2$ be the empirical variances of the big-block sums and of the small-block sums, in the notation of Corollary \ref{corollary 1}. The center $\P_n(h_n)=\frac{1}{n}\sum_{j=1}^m(U_j+W_j)$ and the two empirical variances use all observations.

\begin{theorem}
\label{th 6}
Let $\mathbf{Z}$ be $\tau$-mixing, let $l$ and $m$ be as above, let $d$ be an integer with $1\leq d\leq l/2$, let $h_n:\mathbb{R}\to[a_n,b_n]$ be Lipschitz, and let $[\mathbb{V}]_m^U$ and $[\mathbb{V}]_m^W$ be as above.
\begin{enumerate}[(i)]
    \item For any $t_1,t_2,t_3,t_4>0$ such that $2t_1<m$ and $2t_3<m$, and any deterministic $x>0$,
    \begin{align}
    \label{eq th 6 1}
        \P\Big(\big|\P_n(h_n)-\mu_n\big|\leq L_n\Big)\geq 1-2e^{-t_1}-e^{-t_2}-2e^{-t_3}-e^{-t_4}-\frac{\|h_n\|_{\text{Lip}}\,n\,\tau_\infty(d)}{x},
    \end{align}
    where, with $b_U=(l-d)(b_n-a_n)$ and $b_W=d\,(b_n-a_n)$,
    \begin{align}
        &L_n=\frac{1}{n-2t_1l}\bigg(\sqrt{2t_1[\mathbb{V}]_m^U+b_U^2t_1t_2}+\sqrt{2t_1}\,x+\Big(1+\sqrt{\frac{2t_1}{m}}\Big)\Big(b_U\sqrt{t_1t_2}+\frac{b_Ut_1}{3}\Big)\bigg)\notag\\
        &\quad+\frac{1}{n-2t_3l}\bigg(\sqrt{2t_3[\mathbb{V}]_m^W+b_W^2t_3t_4}+\sqrt{2t_3}\,x+\Big(1+\sqrt{\frac{2t_3}{m}}\Big)\Big(b_W\sqrt{t_3t_4}+\frac{b_Wt_3}{3}\Big)\bigg)+\frac{x}{n},\label{eq th 6 2}
    \end{align}
    and each of the one-sided events $\{\P_n(h_n)-\mu_n\leq L_n\}$ and $\{\mu_n-\P_n(h_n)\leq L_n\}$ has probability no less than $1-e^{-t_1}-e^{-t_2}-e^{-t_3}-e^{-t_4}-\|h_n\|_{\text{Lip}}\,n\,\tau_\infty(d)/x$. Note that the second line of \eqref{eq th 6 2} stems from the sums over the small blocks and from the slack $x$. Under the assumptions of (ii), it is of smaller order than the leading term $\sqrt{2t_1[\mathbb{V}]_m^U}/n$ and vanishes after multiplication by $\sqrt{n}$.
    \item Let $l=l_n$, $d=d_n$ and $x=x_n$ be deterministic, with $m=n/l_n$ for the values of $n$ considered and $1\leq d_n\leq l_n/2$, and suppose that, as $n\to\infty$, $d_n\log n/l_n\to 0$, $l_n^2\log n/n\to 0$, $x_n\sqrt{\log n}/\sqrt{n}\to 0$ and $n\tau_\infty(d_n)/x_n\to 0$. Suppose further that $\sup_n\max\{b_n-a_n,\|h_n\|_{\text{Lip}}\}<\infty$ and $\sum_{r=1}^\infty\tau_1(r)<\infty$. Denote $\sigma_n^2=\frac{1}{n}\text{Var}\big(\sum_{i=1}^nh_n(Z_i)\big)$ and for $\alpha\in(0,1)$, let $L_n$ be the radius in \eqref{eq th 6 2} with $t_1=\log(2/\alpha)$, $t_2=t_3=t_4=\log n$ and $x=x_n$. Then the right-hand side of \eqref{eq th 6 1} tends to $1-\alpha$. If, moreover, $\liminf_n\sigma_n^2>0$, then
    \begin{align*}
        \frac{\sqrt{n}L_n}{\sigma_n\sqrt{2\log(2/\alpha)}}\xrightarrow[]{\P}1.
    \end{align*}
\end{enumerate}
\end{theorem}

Theorem \ref{th 6} does not use the martingale representation behind Theorems \ref{th 4} and \ref{th 5}. Instead, Lemma 5 of \citeA{dedecker2004coupling} provides, for a real integrable random variable $X$ and a $\sigma$-algebra $\mathcal{M}$, a copy of $X$ that is independent of $\mathcal{M}$ and lies at $L^1$-distance $\tau(\mathcal{M},X)$ from $X$, as soon as a uniform random variable independent of $X$ and $\mathcal{M}$ is available. Applied successively to the big-block sums, it yields independent random variables with the same distributions as these sums. Since a big-block sum is a Lipschitz function of $l-d$ observations, this step requires the coefficient $\tau_\infty$ rather than $\tau_1$ or $\tau_2$. The small-block sums are coupled in the same way, and the total $L^1$-distance of both couplings is at most $\|h_n\|_{\text{Lip}}\,n\,\tau_\infty(d)$. The same construction underlies the Bennett-type inequality in Theorem 1 of \citeA{dedecker2004coupling}. Here, both families of independent copies satisfy the assumptions of Corollary \ref{corollary 1}, which supplies the empirical variances, the exact solution in the unknown mean and the sharp leading term. The coupling error is moved into the radius on an event whose probability is bounded by Markov's inequality, which produces the slack $x$ and the term $\|h_n\|_{\text{Lip}}\,n\,\tau_\infty(d)/x$. Since the small blocks are controlled by their own empirical variance, their contribution to the radius is of order $\sqrt{d\log n/l}$ relative to the leading term, rather than of the order $d\sqrt{\log n/l}$ that a bound based on the range alone, such as the inequality of \citeA{hoeffding1963probability}, would give.

 Part (i) uses neither stationarity nor summability of the coefficients. The radius contains no dependence coefficient: apart from the slack $x$ and the length $d$ of the small blocks, it only involves the range $b_n-a_n$ and the data. The dependence of $\mathbf{Z}$ enters only through the coverage loss $\|h_n\|_{\text{Lip}}\,n\,\tau_\infty(d)/x$, in which $\tau_\infty(d)$ may be replaced by any known upper bound. In part (ii), the leading term $\sqrt{2t_1[\mathbb{V}]_m^U}/n$ is sharp in the two-sided sense of \eqref{eq 5}, with the long-run variance in place of the marginal variance, as in Theorems \ref{th 4} and \ref{th 5}. The summability of $\tau_1$ only serves to identify the scale of the big-block sums with $\sigma_n^2$. If $\tau_\infty(r)\leq C\rho^r$ for some $C>0$ and $\rho\in(0,1)$, as for the autoregressions in Section \ref{sec 3.2.1}, the choice $l_n=\lfloor n^{1/3}\rfloor$, $d_n=\lceil 2\log n/\log(1/\rho)\rceil$ and $x_n=\sqrt{n}/\log n$ satisfies the assumptions of (ii) for $n$ large enough, and the coupling term in the coverage loss is at most $\|h_n\|_{\text{Lip}}\,C\,n^{-3/2}\log n$. By contrast, the coverage loss $\eta_n(x)$ in Theorem \ref{th 5} is of moment type and requires a slack $x$ that is large compared with $m$ and with $\sqrt{nl^3}$. The two theorems are complementary: Theorem \ref{th 5} applies to unbounded Lipschitz functions with moments of order $q>4$ under $\tau_1$ and $\tau_2$, whereas Theorem \ref{th 6} requires bounded Lipschitz functions and the stronger coefficient $\tau_\infty$, in exchange for a coverage loss that is polynomially small already for small blocks of logarithmic length. Functions of bounded variation that are not Lipschitz, such as indicator functions, are covered by Theorem \ref{th 4} under $\tilde{\phi}$-mixing.

\subsection{Comparison With Other Results}
\label{sec comparison}
\paragraph{Sharpness.} The notion of sharpness in \eqref{eq 5} is that of \citeA{wang2024sharp}: the first-order term of an empirical Bernstein inequality, including constants, asymptotically matches the oracle Bernstein inequality with known variance. The inequality of \citeA{maurer2009empirical} is not sharp in this sense, and \citeA{waudby2024estimating} were the first to prove sharpness. Theorem \ref{th 1} and Corollary \ref{corollary 1} are sharp whenever $\<\M\>_n$ is deterministic, without identically distributed observations or a convergent average variance, and Theorems \ref{th 4}, \ref{th 5} and \ref{th 6} are sharp with the long-run variance.

\paragraph{Weak dependence.} Concentration inequalities for weakly dependent processes involve population quantities and constants that depend on the dependence structure, so that they cannot be evaluated from the sample alone. The Bernstein-type inequality of \citeA{merlevede2011bernstein} for strongly mixing and $\tau$-mixing sequences with coefficients of order $\exp(-cr^{\gamma_1})$ involves a variance factor defined through covariances of truncated variables and constants that are not explicit. The Bernstein inequalities of \citeA{paulin2015concentration} for stationary Markov chains with a spectral gap involve the asymptotic variance, the marginal variance and the spectral gap. The Fuk--Nagaev inequality of \citeA{babii2024high} for $\tau$-mixing processes involves an unspecified constant, the moments and the decay rate of the coefficients, and the Bennett-type inequality of \citeA{dedecker2004coupling} involves the variance of block sums. Theorems \ref{th 4}, \ref{th 5} and \ref{th 6} cover $\tilde{\phi}$- and $\tau$-mixing processes, including processes that are not strongly mixing. Their leading terms only involve a block sample variance and are sharp with the long-run variance, and the radii of Theorems \ref{th 5} and \ref{th 6} contain no dependence coefficient. Additionally, the empirical Bernstein inequalities of \citeA{waudby2024estimating} and \citeA{wang2024sharp} require a constant conditional mean, $\E[X_i| \mathcal{F}_{i-1}] = \mu$, which excludes autocorrelated processes, and their radii are calibrated by the marginal rather than the long-run variance.

\paragraph{General functions.} Theorems \ref{th 4}, \ref{th 5} and \ref{th 6} are stated for triangular arrays $h_n$ from the outset, with explicit dependence on the range or the moments and on the total variation or the Lipschitz constant of $h_n$. Functions that change with $n$, such as kernel smoothers $z\mapsto K((x-z)/b_n)/b_n$, therefore enter directly, with a closed-form radius. Betting-based intervals are designed for the mean of observations in a known bounded interval. For such triangular arrays, the choice of the betting fractions, the truncation and the initial values, for which there is no universal rule, may need to be adapted to the bandwidth.

\section{Application: Expected Loss of a Forecasting Rule}
\label{sec application}
A forecasting rule is judged by its risk, the expected loss it incurs on new observations of the process. The risk is unknown, and what we observe are the losses of the rule along a sample path. These losses inherit the dependence of the time series, so that inequalities for independent observations or martingale differences do not apply to them. Theorem \ref{th 6} does apply, and this section turns it into an upper bound on the risk that is computed from the observed losses alone.

We focus on time series whose dependence coefficients decay exponentially, $\tau_\infty(r)\leq C\rho^r$ for some $C>0$ and $\rho\in(0,1)$. Many time series models satisfy this condition. For instance, if $Z_t=g(Z_{t-1})+\xi_t$ with i.i.d.\ integrable innovations $\xi_t$ and a $\kappa$-Lipschitz function $g$ with $\kappa<1$, then $\tau_\infty(r)\leq 2\E|Z_0|\kappa^r$ by Section 7.2 of \citeA{dedecker2005new}. This covers linear autoregressions of order one, threshold autoregressions with $g(z)=a_1z\mathbf{1}\{z\geq 0\}+a_2z\mathbf{1}\{z<0\}$ and $|a_1|,|a_2|<1$, and the bounded autoregression of Example \ref{example 1}, which is not strongly mixing. Theorem \ref{th 7} itself only requires summable coefficients. Exponential decay makes its coverage loss vanish quickly.

\subsection{Setting and Guarantee}
\label{sec application 1}
Let $\mathbf{Z}$ be strictly stationary, let $p\geq 0$ be an integer and let $f:\mathbb{R}^p\to\mathbb{R}$ be a forecasting rule that predicts $Z_t$ from $(Z_{t-p},\dots,Z_{t-1})$. The accuracy of $f$ is measured by a loss $\ell:\mathbb{R}^{p+1}\to[a,b]$ that is $L$-Lipschitz with respect to $\|\cdot\|_1$, where $a<b$ and $L>0$ are known. For instance, if $f$ is $L_f$-Lipschitz, the clipped absolute error $\ell(z_1,\dots,z_{p+1})=\min\{|z_{p+1}-f(z_1,\dots,z_p)|,c\}$ takes values in $[0,c]$ and is $\max\{1,L_f\}$-Lipschitz. For $k\in\mathbb{Z}$, denote $\ell_k=\ell(Z_{k-p},\dots,Z_k)$. The risk of $f$ is its expected loss $\mathcal{R}=\E[\ell_k]$, which does not depend on $k$ by stationarity, and the long-run variance of the losses is $\sigma^2=\sum_{k\in\mathbb{Z}}\mathrm{Cov}(\ell_0,\ell_k)$. A sample $Z_{1-p},\dots,Z_n$ yields the $n$ losses $\ell_1,\dots,\ell_n$ and the empirical risk $\hat{\mathcal{R}}_n=\frac{1}{n}\sum_{i=1}^n\ell_i$.

The radius is computed from blocks of losses. The block length is $l=\lfloor n^{1/3}\rfloor$ and the number of blocks is $m=\lfloor n/l\rfloor$. If $l$ does not divide $n$, the last $n-ml$ losses are discarded, as in Section \ref{sec 3.2.2}, and we assume $n=ml$ to simplify notation. Each block $B_j=\{(j-1)l+1,\dots,jl\}$ is split into the big block $I_j$, formed by its first $l-d$ indices, and the small block $J_j$, formed by its last $d=\lceil\log^2n\rceil$ indices. With the big-block sums $U_j=\sum_{i\in I_j}\ell_i$ and their empirical variance $[\mathbb{V}]_m^U=\sum_{j=1}^m(U_j-\frac{1}{m}\sum_{k=1}^mU_k)^2$, the leading term of the radius is
\begin{align}
\label{eq lead}
    \bar{L}_n(\alpha)=\frac{\sqrt{2\log(1/\alpha)[\mathbb{V}]_m^U}}{n},\qquad\alpha\in(0,1).
\end{align}
The small blocks separate consecutive big blocks, so that the big-block sums are nearly independent and their empirical variance estimates the long-run variance. The losses in the small blocks enter the radius of Theorem \ref{th 6} only through lower-order terms. The choices of $l$ and $d$ make these terms negligible, and the slack $x$ of Theorem \ref{th 6} is the only free parameter.

\begin{theorem}
\label{th 7}
Let $\mathbf{Z}$ be strictly stationary and $\tau$-mixing with $\sum_{r=1}^\infty\tau_\infty(r)<\infty$, let $\sigma^2>0$ and let $\alpha\in(0,1)$. There exists $N_0$, depending only on $\alpha$, $b-a$, $L$, $p$, $\sigma^2$ and the sequence $(\tau_\infty(r))_{r\geq 1}$, such that, for every $n\geq N_0$ and every deterministic $x$ with $0<x\leq\sqrt{n}/\log^{3/2}n$,
\begin{align}
\label{eq th 7 1}
    \P\Big(\mathcal{R}\leq\hat{\mathcal{R}}_n+\Big(1+\frac{1}{\sqrt{\log n}}\Big)\bar{L}_n(\alpha)\Big)\geq 1-\alpha-\beta_n(x),
\end{align}
where
\begin{align*}
    \beta_n(x)=\frac{5}{n}+\frac{(p+1)L\,(n+2pm)\,\tau_\infty(d-p)}{x}.
\end{align*}
Moreover, $\sqrt{n}\,\bar{L}_n(\alpha)\xrightarrow[]{\P}\sigma\sqrt{2\log(1/\alpha)}$ as $n\to\infty$, that is, the leading term is sharp in the sense of \eqref{eq 5}, with the long-run variance in place of the marginal variance.
\end{theorem}

 Theorem \ref{th 7} bounds the risk from above by the empirical risk plus a radius. The radius is the leading term $\bar{L}_n(\alpha)$, inflated by the factor $1+1/\sqrt{\log n}$. It is computed from the observed losses and contains no dependence coefficient. Since the inflation factor tends to one, the last statement of Theorem \ref{th 7} shows that the radius itself is sharp in the sense of \eqref{eq 5}, with the long-run variance in place of the marginal variance. The coverage loss $\beta_n(x)$ consists of the failure probability $5/n$ of auxiliary events and of a coupling term, through which the dependence of $\mathbf{Z}$ enters at the single lag $d-p$. The slack $x$ governs the trade-off between the two sides: a larger $x$ reduces the coupling term, while the restriction $x\leq\sqrt{n}/\log^{3/2}n$ keeps the inflation of the radius below $1/\sqrt{\log n}$, so that the upper end of this range gives the smallest coverage loss. Under the exponential decay considered in this section, since $d\geq\log^2n$, the coupling term with $x=\sqrt{n}/\log^{3/2}n$ decays faster than any power of $n$.

The following corollary replaces the big-block sums by the full block sums, so that the leading term only involves the block length $l$ and is defined whenever $m\geq 2$.

\begin{corollary}
\label{corollary 2}
Let the assumptions of Theorem \ref{th 7} hold, let $Y_j=\sum_{i\in B_j}\ell_i$, $1\leq j\leq m$, let $[\mathbb{V}]_m^Y=\sum_{j=1}^m(Y_j-\frac{1}{m}\sum_{k=1}^mY_k)^2$ and let
\begin{align}
\label{eq lead full}
    \tilde{L}_n(\alpha)=\frac{\sqrt{2\log(1/\alpha)[\mathbb{V}]_m^Y}}{n}.
\end{align}
There exists $N_0$, depending only on $\alpha$, $b-a$, $L$, $p$, $\sigma^2$ and the sequence $(\tau_\infty(r))_{r\geq 1}$, such that, for every $n\geq N_0$ and every deterministic $x$ with $0<x\leq\sqrt{n}/\log^{3/2}n$,
\begin{align}
\label{eq cor 2}
    \P\Big(\mathcal{R}\leq\hat{\mathcal{R}}_n+\Big(1+\frac{1}{\sqrt{\log n}}\Big)\tilde{L}_n(\alpha)\Big)\geq 1-\alpha-\beta_n(x),
\end{align}
where $\beta_n(x)$ is given in Theorem \ref{th 7}. Moreover, $\sqrt{n}\,\tilde{L}_n(\alpha)\xrightarrow[]{\P}\sigma\sqrt{2\log(1/\alpha)}$ as $n\to\infty$.
\end{corollary}

The upper bound $\hat{\mathcal{R}}_n+(1+1/\sqrt{\log n})\tilde{L}_n(\alpha)$ is computed from the observed losses with the single block length $l=\lfloor n^{1/3}\rfloor$. It involves neither the length $d$ of the small blocks nor the slack $x$ nor any dependence coefficient, and it can be computed for every sample size with at least two blocks. The small blocks only serve to decouple the blocks in the proof: since they form a fraction $d/l$ of each block, replacing the big-block sums by the full block sums changes the leading term only by a factor of smaller order than the inflation $1/\sqrt{\log n}$. The guarantee is therefore the same as in Theorem \ref{th 7}, with the same coverage loss $\beta_n(x)$. Section \ref{sec application 2} evaluates this upper bound for sample sizes between $3000$ and $10000$.

\subsection{Simulation}
\label{sec application 2}
Theorem \ref{th 7} and Corollary \ref{corollary 2} guarantee the upper bound for $n\geq N_0$ and show that its leading term is sharp as $n\to\infty$. The simulation examines how the construction performs at sample sizes of practical interest. We evaluate the upper bound of Corollary \ref{corollary 2} for four strictly stationary processes:
\begin{enumerate}[(a)]
    \item the bounded autoregression $Z_t=0.5Z_{t-1}+0.5\varepsilon_t$ of Example \ref{example 1}, with i.i.d.\ innovations satisfying $\P(\varepsilon_t=0)=\P(\varepsilon_t=1)=\frac{1}{2}$, which is not strongly mixing,
    \item the bounded autoregression $Z_t=0.8Z_{t-1}+0.2\varepsilon_t$ with the same innovations,
    \item the threshold autoregression $Z_t=0.8Z_{t-1}\mathbf{1}\{Z_{t-1}\geq 0\}+0.3Z_{t-1}\mathbf{1}\{Z_{t-1}<0\}+\varepsilon_t$ with i.i.d.\ standard normal innovations,
    \item the bounded autoregression $Z_t=0.9Z_{t-1}+0.1\varepsilon_t$ with the same innovations as in (a) and (b).
\end{enumerate}
Each process is of the form $Z_t=g(Z_{t-1})+\xi_t$ with i.i.d.\ $\xi_t$ and a $\kappa$-Lipschitz function $g$, where $\kappa=0.5$, $0.8$, $0.8$ and $0.9$, so that $\tau_\infty(r)\leq 2\E|Z_0|\kappa^r$ by Section 7.2 of \citeA{dedecker2005new}. The forecasting rule uses the last observation, $p=1$. In (a), (b) and (d), $f(z)=\frac{1}{2}+\frac{\rho}{2}(z-\frac{1}{2})$, where $\rho$ is the autoregressive coefficient, so that $f$ shrinks the conditional mean $\frac{1}{2}+\rho(z-\frac{1}{2})$ halfway towards the unconditional mean. In (c), the linear rule $f(z)=z/2$ ignores the threshold. The loss is the clipped squared error $\ell(z_1,z_2)=\min\{(z_2-f(z_1))^2,c\}$, with $c=1$ in (a), (b) and (d), where the clipping is not active, and $c=4$ in (c). It takes values in $[0,c]$ and is $2\sqrt{c}$-Lipschitz. Since the rules are misspecified, the losses are positively autocorrelated, and their long-run variance exceeds their marginal variance by a factor of $1.34$, $3.04$, $1.66$ and $6.64$ in (a), (b), (c) and (d). In (a), (b) and (d), the risk $\mathcal{R}$ and the long-run variance $\sigma^2$ are available in closed form; in (c), they are computed from the autocovariances of the losses up to lag $100$, estimated from $2000$ independent paths of length $5\times 10^4$. For each process and each $n\in\{3000,5000,8000,10000\}$, we draw $5000$ independent samples after a burn-in of $500$ observations and compute the upper bound of Corollary \ref{corollary 2} with $\alpha=0.05$ and $l=\lfloor n^{1/3}\rfloor$, that is, with block lengths between $14$ and $21$ and between $214$ and $476$ blocks. The simulation addresses three questions.
\begin{itemize}
    \item \textbf{Validity.} Does the upper bound $\hat{\mathcal{R}}_n+(1+1/\sqrt{\log n})\tilde{L}_n(\alpha)$ of Corollary \ref{corollary 2} hold with probability at least $1-\alpha$ for these sample sizes?
    \item \textbf{Leading term.} Is the inflation factor needed in practice, that is, does $\hat{\mathcal{R}}_n+\tilde{L}_n(\alpha)$ hold with probability at least $1-\alpha$?
    \item \textbf{Sharpness.} How close is the leading term $\tilde{L}_n(\alpha)$ to the oracle term $\sqrt{2\log(1/\alpha)}\,\sigma/\sqrt{n}$, which uses the long-run variance?
\end{itemize}

\begin{figure}[ht]
\centering
\includegraphics[width=0.49\textwidth]{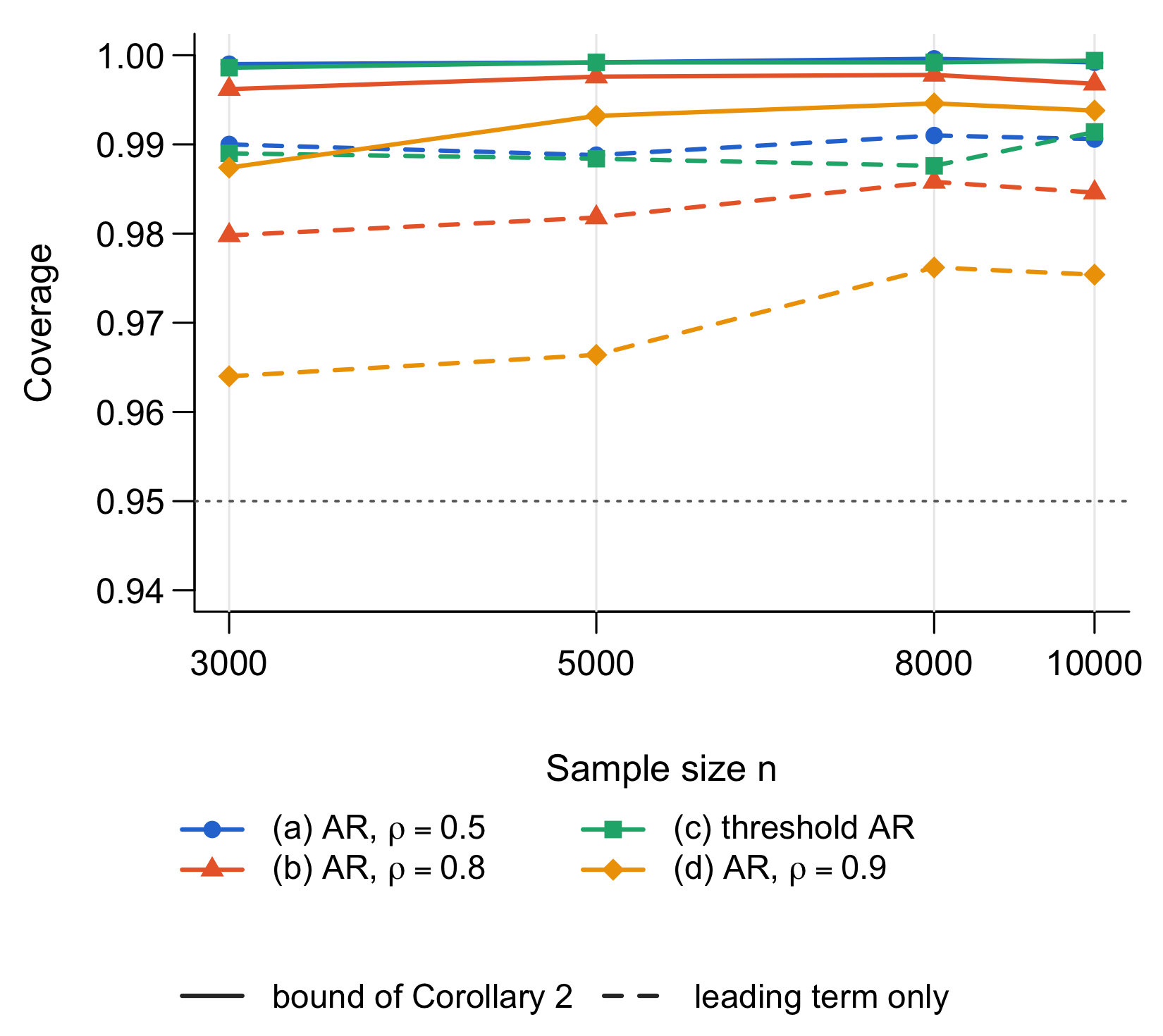}\hfill
\includegraphics[width=0.49\textwidth]{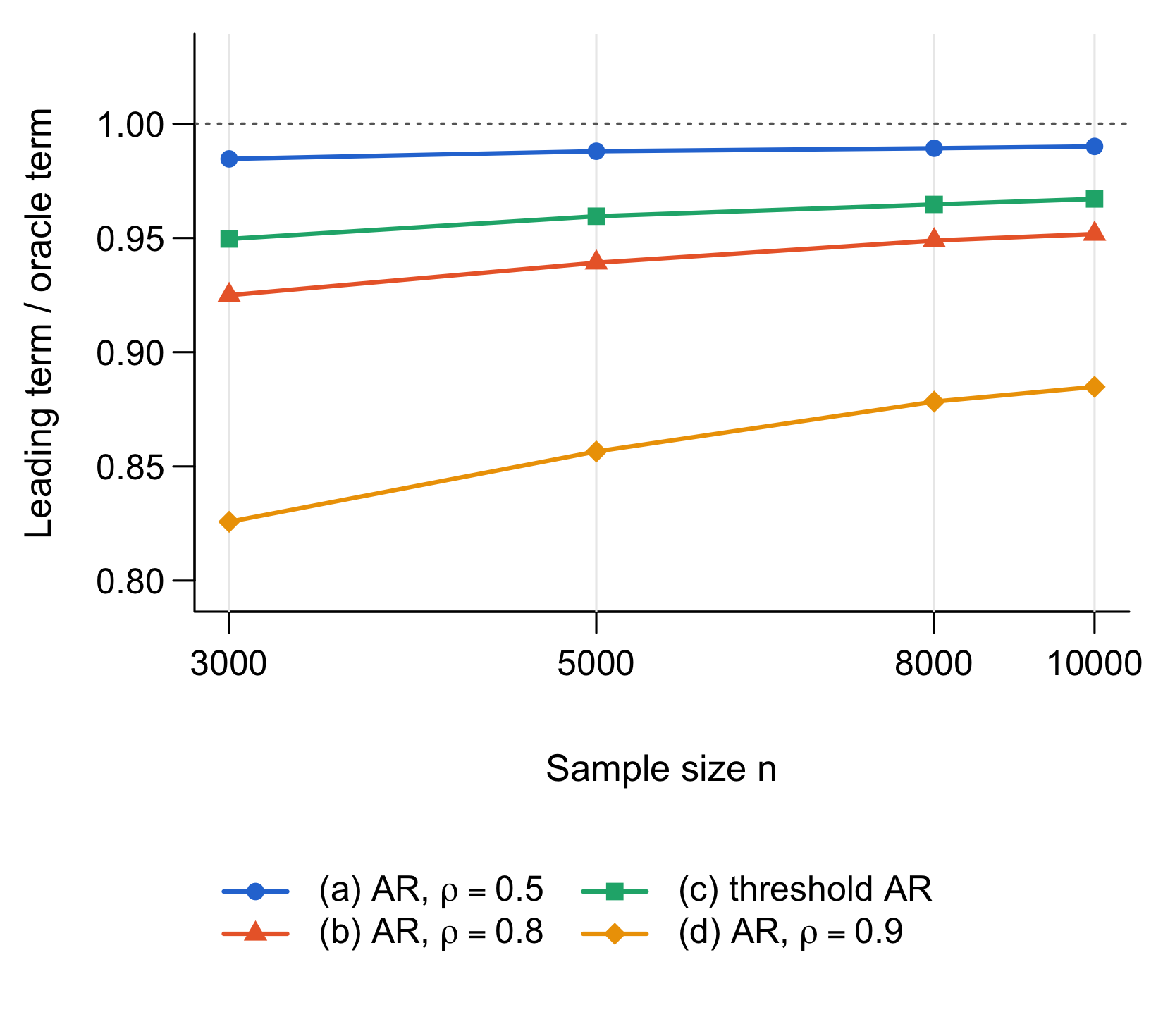}
\caption{Left: frequency, over $5000$ samples, with which the upper bound of Corollary \ref{corollary 2} (solid) and the empirical risk plus the leading term alone (dashed) lie above the risk $\mathcal{R}$, for the designs (a)--(d). The dotted line is $1-\alpha=0.95$. Right: average ratio $\tilde{L}_n(\alpha)/(\sqrt{2\log(1/\alpha)}\,\sigma/\sqrt{n})$ of the leading term to the oracle term. The dotted line is $1$.}
\label{fig sim}
\end{figure}

Figure \ref{fig sim} answers the three questions. The upper bound of Corollary \ref{corollary 2} holds in at least $99.6\%$ of the samples for (a), (b) and (c), and in at least $98.7\%$ for (d), for all sample sizes. Without the inflation factor, the upper bound holds in at least $97.9\%$ of the samples for (a), (b) and (c), and in at least $96.4\%$ for (d), where the lowest frequency occurs at $n=3000$. All Monte Carlo standard errors are below $0.3$ percentage points. The ratio of the leading term to the oracle term increases with $n$ for all designs, from $0.985$ to $0.990$ for (a), from $0.950$ to $0.967$ for (c), from $0.925$ to $0.952$ for (b) and from $0.826$ to $0.885$ for (d). It is below $1$ because $\text{Var}(Y_1)/l=\sum_{|k|<l}(1-|k|/l)\,\mathrm{Cov}(\ell_0,\ell_k)$ down-weights the autocovariances within a block and omits those at lags beyond it; in these designs, the down-weighting within the block accounts for most of the gap. With at least $214$ blocks, the loss of one degree of freedom in $[\mathbb{V}]_m^Y$ reduces the ratio by less than $0.3\%$. The underestimation grows with the strength of the dependence and is strongest for (d), whose long-run variance is about $6.6$ times its marginal variance. The frequencies nevertheless exceed $1-\alpha$, since the leading term carries the constant $\sqrt{2\log(1/\alpha)}$ of the oracle Bernstein inequality, which bounds the tail probability for every sample size rather than asymptotically. Hence, in all four designs, the empirical frequencies exceed $1-\alpha$ from $n=3000$ on, both for the upper bound of Corollary \ref{corollary 2} and for the leading term alone, so that the inflation factor is not needed in these designs.

\section{Conclusions}
This paper develops self-normalized and empirical Bernstein inequalities for the marginal mean of dependent data that use only a sample variance. For martingale differences, the inequality is available in closed form, holds with probability at least $1-\alpha$ for every sample size and is sharp under heteroscedasticity. For $\tilde{\phi}$- and $\tau$-mixing processes, including processes that are not strongly mixing, the leading terms only involve a block sample variance and are sharp with respect to the long-run variance, and for $\tau$-mixing processes the radii contain no dependence coefficient. Applied to the losses of a forecasting rule, these results yield a sharp upper bound on its risk that is computed from the observed losses alone, and in simulations with sample sizes between $3000$ and $10000$, the bound and its leading term alone hold in more than $95\%$ of the samples, also under strong dependence.

Several extensions are natural. First, Proposition \ref{prop 1} shows that the self-normalized radius can fail when the predictable quadratic variation is random, so that conditionally heteroscedastic sequences, such as GARCH-type processes, call for a different normalization. Second, the risk bound of Section \ref{sec application} concerns a single forecasting rule, and extending it to classes of rules, as required for model selection, is a natural next step. Third, the functions in Theorems \ref{th 4}, \ref{th 5} and \ref{th 6} may depend on the sample size, so that the inequalities apply directly to kernel smoothers and other nonparametric estimators under dependence. Finally, the block length in Section \ref{sec application} is fixed at $\lfloor n^{1/3}\rfloor$, and a data-driven choice that preserves the guarantee remains open.

\section{Acknowledgement}
The author is grateful to Prof. Dr. Aaditya Ramdas for pointing out an important mistake in the earlier version of this work. The author also thanks him for the guide of literature review.

\bibliographystyle{\bibstylename}
\bibliography{reference}

\appendix
\setcounter{equation}{0}
\renewcommand\theequation{A.\arabic{equation}}
\renewcommand\thetheorem{A.\arabic{theorem}}
\renewcommand\thelemma{A.\arabic{lemma}}
\renewcommand\theremark{A.\arabic{remark}}

\section{Proof of Main Results and Technical Lemmas}

\subsection{Technical Lemmas}
Throughout this subsection, $\mathcal{A}\subset\mathcal{F}$ is a sub $\sigma$-algebra and $b>0$. As in \citeA{freedman1975tail}, we denote
\begin{align}
\label{eq psi}
    \psi(\lambda)=\frac{\lambda^2}{2(1-\frac{\lambda b}{3})},\qquad \lambda\in[0,3/b).
\end{align}

\begin{lemma}
    \label{lemma A.1}
    Suppose $Z$ is a random variable such that $\E[Z|\mathcal{A}]=0$ and $Z\leq b$ almost surely. Then, for every $\lambda\in[0,3/b)$,
    \begin{align*}
        \E[e^{\lambda Z}|\mathcal{A}]\leq \exp\big(\psi(\lambda)\E[Z^2|\mathcal{A}]\big)\ \ \text{a.s.}
    \end{align*}
\end{lemma}

\begin{lemma}
    \label{lemma A.2}
    Suppose $Z$ is a random variable such that $|Z|\leq b$ almost surely. Then, for every $\lambda\geq 0$,
    \begin{align*}
        \E[e^{-\lambda Z^2}|\mathcal{A}]\leq \exp\Big(-\lambda\E[Z^2|\mathcal{A}]+\frac{\lambda^2b^2}{2}\E[Z^2|\mathcal{A}]\Big)\ \ \text{a.s.}
    \end{align*}
\end{lemma}

\begin{lemma}
    \label{lemma A.3}
    For any $\nu>0$ and $x>0$, the number
    \begin{align*}
        \lambda(\nu,x):=\frac{\sqrt{2x/\nu}}{1+\frac{b}{3}\sqrt{2x/\nu}}\in(0,3/b)
    \end{align*}
    satisfies
    \begin{align}
    \label{eq lemma A.3 1}
        \lambda(\nu,x)\Big(\sqrt{2\nu x}+\frac{b}{3}x\Big)-\psi(\lambda(\nu,x))\nu=x.
    \end{align}
\end{lemma}

\begin{lemma}
    \label{lemma A.4}
    Let $(\varepsilon_i)_{i\geq 1}$ be independent random variables with $\P(\varepsilon_i=1)=\P(\varepsilon_i=-1)=\frac{1}{2}$ and denote $S_k=\sum_{i=1}^k\varepsilon_i$. Then, for every $K>0$ and every sequence $(a_n)_{n\geq 1}$ such that $1\leq a_n\leq n$ and $n/a_n\to\infty$,
    \begin{align*}
        \lim_{n\to\infty}\P\Big(\exists\, k\in[a_n,n]:\ S_k>K\sqrt{k}\Big)=1.
    \end{align*}
\end{lemma}

\begin{lemma}
    \label{lemma A.5}
    Suppose $Z$ is a random variable such that $\E[Z^2]<\infty$ and $\E[Z|\mathcal{A}]=0$. Then, for every $\lambda\in\mathbb{R}$,
    \begin{align*}
        \E\Big[\exp\Big(\lambda Z-\frac{\lambda^2}{6}\big(Z^2+2\E[Z^2|\mathcal{A}]\big)\Big)\Big|\mathcal{A}\Big]\leq 1\ \ \text{a.s.}
    \end{align*}
\end{lemma}

\begin{lemma}
    \label{lemma A.6}
    Let $\mathcal{M}\subset\mathcal{F}$ be a $\sigma$-algebra, let $X$ and $Y$ be real-valued random variables, and let $f:\mathbb{R}\to\mathbb{R}$ be a function of bounded variation. Denote $\|f\|_\infty=\sup_{x\in\mathbb{R}}|f(x)|$ and let $\tilde{\phi}(\cdot,\cdot)$ be given in Definition \ref{def phi tilde k}. Then, almost surely,
    \begin{enumerate}[(i)]
        \item $\big|\E[f(X)|\mathcal{M}]-\E[f(X)]\big|\leq\|f\|_{TV}\,\tilde{\phi}(\mathcal{M},X)$,
        \item with $\bar{\phi}=\max\{\tilde{\phi}(\mathcal{M},X),\tilde{\phi}(\mathcal{M},Y),\tilde{\phi}(\mathcal{M},(X,Y))\}$,
        \begin{align*}
            \big|\E[f(X)f(Y)|\mathcal{M}]-\E[f(X)f(Y)]\big|\leq\big(\|f\|_{TV}^2+2\|f\|_\infty\|f\|_{TV}\big)\bar{\phi}.
        \end{align*}
    \end{enumerate}
    Both bounds remain true if $\mathcal{M}$ is replaced on the left-hand side by any $\sigma$-algebra $\mathcal{M}'\subset\mathcal{M}$.
\end{lemma}

\noindent Lemmas \ref{lemma A.7}--\ref{lemma A.9} are stated in the setting of Section \ref{sec 3.2.2} and concern the following auxiliary triangular array. For the given $n\geq 1$, let $\mathbf{Z}_n=(Z_{nk})_{k\in\mathbb{Z}}$ be a triangular array such that
\begin{itemize}
    \item [\ding{172}] $(Z_{n1},\dots,Z_{nn})$ has the same distribution as $(Z_1,\dots,Z_n)$,
    \item [\ding{173}] $\sigma(Z_{nk}:k\leq 0)$, $\sigma(Z_{n1},\dots,Z_{nn})$ and $\sigma(Z_{nk}:k\geq n+1)$ are mutually independent,
    \item [\ding{174}] $(Z_{nk})_{k\geq n+1}$ are copies of $Z_{nn}$.
\end{itemize}
Such an array exists on a product probability space. Denote $\mathcal{F}_{ni}=\sigma(Z_{nk}:k\leq i)$ for $i\in\mathbb{Z}$ and, for $1\leq i\leq n$ and $1\leq j\leq m$,
\begin{align}
\label{eq notation block}
    Y_{ni}=h_n(Z_{ni})-\mu_n,\qquad S_n=\sum_{i=1}^nY_{ni},\qquad T_{nj}=\sum_{i\in B_j}Y_{ni},\qquad M_{nj}=\sum_{i\in B_j}D_{ni},
\end{align}
where
\begin{align*}
    D_{ni}=\sum_{k=i}^n\big(\E[h_n(Z_{nk})|\mathcal{F}_{ni}]-\E[h_n(Z_{nk})|\mathcal{F}_{n,i-1}]\big),\qquad \varepsilon_{ni}=\sum_{k>i}\E\big[h_n(Z_{nk})-\mu_n\big|\mathcal{F}_{ni}\big],\quad 0\leq i\leq n.
\end{align*}
We also write $\P_n(h_n)$ and $\hat{V}_n$ for the quantities of Section \ref{sec 3.2.2} computed from $(Z_{n1},\dots,Z_{nn})$.

\begin{lemma}
    \label{lemma A.7}
    In the setting above, the following statements hold.
    \begin{enumerate}[(i)]
        \item Let $0\leq i\leq n$ and let $g:\mathbb{R}^n\to\mathbb{R}$ be measurable with $\E|g(Z_1,\dots,Z_n)|<\infty$. Denote
        \begin{align*}
            g_n=g(Z_{n1},\dots,Z_{nn}),\qquad g_0=g(Z_1,\dots,Z_n).
        \end{align*}
        Then, almost surely, $\E[g_n|\mathcal{F}_{ni}]=\E[g_n|\sigma(Z_{n1},\dots,Z_{ni})]$, where $\sigma(Z_{n1},\dots,Z_{ni})$ is trivial for $i=0$. Moreover, for every $s\in[1,\infty]$,
        \begin{align*}
            \big\|\E[g_n|\mathcal{F}_{ni}]-\E[g_n]\big\|_s\leq\big\|\E[g_0|\mathcal{F}_{-\infty}^i]-\E[g_0]\big\|_s.
        \end{align*}
        \item $\varepsilon_{ni}=\sum_{k=i+1}^n\E[Y_{nk}|\mathcal{F}_{ni}]$ for $0\leq i\leq n$, $\varepsilon_{n0}=\varepsilon_{nn}=0$, and $M_{nj}=T_{nj}+\varepsilon_{n,jl}-\varepsilon_{n,(j-1)l}$ for $1\leq j\leq m$. Moreover, $(M_{nj},\mathcal{F}_{n,jl})_{1\leq j\leq m}$ is a martingale difference sequence and $\sum_{j=1}^mM_{nj}=S_n=n(\P_n(h_n)-\mu_n)$.
        \item If $h_n:\mathbb{R}\to[a_n,b_n]$ is of bounded variation, then $|\varepsilon_{ni}|\leq\|h_n\|_{TV}\tilde{\Phi}_n$ for $0\leq i\leq n$, and $|M_{nj}|\leq\bar{b}_n$ for $1\leq j\leq m$.
        \item $\sum_{j=1}^mT_{nj}^2=n\hat{V}_n+nl\big(\P_n(h_n)-\mu_n\big)^2$.
    \end{enumerate}
\end{lemma}

\begin{lemma}
    \label{lemma A.8}
    In the setting above, let $h_n:\mathbb{R}\to[a_n,b_n]$ be of bounded variation. Then, for every $1\leq j\leq m$,
    \begin{align*}
        \Big|\E\big[M_{nj}^2\big|\mathcal{F}_{n,(j-1)l}\big]-\E\big[M_{nj}^2\big]\Big|\leq\frac{\Gamma_n}{2m}\ \ \text{a.s.}
    \end{align*}
    Consequently, $\big|\sum_{j=1}^m\E[M_{nj}^2|\mathcal{F}_{n,(j-1)l}]-\text{Var}(S_n)\big|\leq \Gamma_n/2$ almost surely.
\end{lemma}

\begin{lemma}
    \label{lemma A.9}
    In the setting above, let $\mathbf{Z}$ be $\tau$-mixing, let $q>4$, let $h_n$ be Lipschitz with $\kappa_n<\infty$, and let $\beta$, $\gamma$, $\Psi_n^\tau$, $\Phi_n^\tau$ and $c_n$ be given in Section \ref{sec 3.2.2}. Then the following statements hold.
    \begin{enumerate}[(i)]
        \item $\E[\varepsilon_{ni}^2]\leq 2c_n\Psi_n^\tau$ and $\|\varepsilon_{ni}\|_4\leq\kappa_n^{1-\gamma}\|h_n\|_{\text{Lip}}^\gamma\Phi_n^\tau$ for $0\leq i\leq n$.
        \item $\E\big|\sum_{j=1}^m\E[M_{nj}^2|\mathcal{F}_{n,(j-1)l}]-\text{Var}(S_n)\big|\leq 52\,m\,c_n\Psi_n^\tau$.
        \item $\sum_{j=1}^m\E\big[(\varepsilon_{n,jl}-\varepsilon_{n,(j-1)l})^2\big]\leq 8\,m\,c_n\Psi_n^\tau$.
        \item $\sum_{j=1}^m\E[M_{nj}^4]\leq 8nl^3\kappa_n^4+128\,m\,\kappa_n^{4-4\gamma}\|h_n\|_{\text{Lip}}^{4\gamma}(\Phi_n^\tau)^4$.
    \end{enumerate}
\end{lemma}

\noindent \textbf{Proof of Lemma \ref{lemma A.1}} For $\lambda\geq 0$, define $g(x)=(e^{\lambda x}-1-\lambda x)/x^2$ for $x\neq 0$ and $g(0)=\lambda^2/2$. Taylor's formula with integral remainder gives $g(x)=\lambda^2\int_0^1(1-s)e^{s\lambda x}ds$, so $g$ is non-decreasing on $\mathbb{R}$. Since $Z\leq b$, we obtain $e^{\lambda Z}\leq 1+\lambda Z+Z^2g(b)$. Taking conditional expectations and using $\E[Z|\mathcal{A}]=0$ and $1+y\leq e^y$,
\begin{align*}
    \E[e^{\lambda Z}|\mathcal{A}]\leq 1+g(b)\E[Z^2|\mathcal{A}]\leq \exp\big(g(b)\E[Z^2|\mathcal{A}]\big).
\end{align*}
Finally, since $k!\geq 2\cdot 3^{k-2}$ for every $k\geq 2$ and $\lambda b<3$,
\begin{align*}
    g(b)=\sum_{k=2}^\infty\frac{\lambda^kb^{k-2}}{k!}\leq \frac{\lambda^2}{2}\sum_{k=2}^{\infty}\Big(\frac{\lambda b}{3}\Big)^{k-2}=\psi(\lambda).
\end{align*}

\noindent \textbf{Proof of Lemma \ref{lemma A.2}} The function $q(x)=1-x+\frac{x^2}{2}-e^{-x}$ satisfies $q(0)=0$ and $q'(x)=x-1+e^{-x}\geq 0$ for $x\geq 0$, so $e^{-x}\leq 1-x+\frac{x^2}{2}$ holds for all $x\geq 0$. Applying it to $x=\lambda Z^2$ and using $Z^4\leq b^2Z^2$ yields
\begin{align*}
    \E[e^{-\lambda Z^2}|\mathcal{A}]\leq 1-\lambda\E[Z^2|\mathcal{A}]+\frac{\lambda^2}{2}\E[Z^4|\mathcal{A}]\leq 1-\lambda\E[Z^2|\mathcal{A}]+\frac{\lambda^2b^2}{2}\E[Z^2|\mathcal{A}].
\end{align*}
The claim follows from $1+y\leq e^y$.\\

\noindent \textbf{Proof of Lemma \ref{lemma A.3}} Write $u=\sqrt{2x/\nu}$ and $c=b/3$, so that $\lambda:=\lambda(\nu,x)=\frac{u}{1+cu}$, $1-c\lambda=\frac{1}{1+cu}$ and $\psi(\lambda)=\frac{u^2}{2(1+cu)}$. Since $\sqrt{2\nu x}=\nu u$ and $\nu u^2=2x$,
\begin{align*}
    \lambda\big(\sqrt{2\nu x}+cx\big)-\psi(\lambda)\nu=\frac{\nu u^2+cux}{1+cu}-\frac{\nu u^2}{2(1+cu)}=\frac{x+cux}{1+cu}=x.
\end{align*}

\noindent \textbf{Proof of Lemma \ref{lemma A.4}} Fix $K>0$ and an integer $\theta\geq\max\{2,K^{-2}\}$. Put
\begin{align*}
    \kappa=\frac{1}{K^2\theta}\in(0,1],\qquad q=\frac{1}{2}\P\big(\mathcal{N}>2\sqrt{2}K\big)>0,
\end{align*}
where $\mathcal{N}$ denotes a standard normal random variable. Let $a_n'=\sqrt{na_n}$, so that $a_n\leq a_n'\leq n$, $a_n'\geq\sqrt{n}$ and $n/a_n'=\sqrt{n/a_n}\to\infty$. Define the set of indices
\begin{align*}
    J_n=\big\{j\geq 1:\ a_n'\leq\theta^{j-1},\ \theta^j\leq n\big\},
\end{align*}
whose cardinality satisfies $|J_n|\geq\log_\theta(n/a_n')-2\to\infty$. For $j\in J_n$, denote $\Delta_j=S_{\theta^j}-S_{\theta^{j-1}}$ and $H_j=\{\Delta_j>2K\theta^{j/2}\}$.\\

\noindent \textbf{Step 1} (The events $H_j$) The events $H_j$, $j\in J_n$, are independent, and $H_j$ is independent of $\sigma(\varepsilon_1,\dots,\varepsilon_{\theta^{j-1}})$. Moreover, $\Delta_j$ is a sum of $\theta^j-\theta^{j-1}\geq\theta^{j-1}\geq\sqrt{n}$ independent Rademacher variables, and
\begin{align*}
    \frac{2K\theta^{j/2}}{\sqrt{\theta^j-\theta^{j-1}}}=\frac{2K}{\sqrt{1-1/\theta}}\leq 2\sqrt{2}K.
\end{align*}
By the central limit theorem applied with the fixed threshold $2K/\sqrt{1-1/\theta}$, there exists $n_0=n_0(K,\theta)$ such that $\P(H_j)\geq q$ for all $j\in J_n$ and $n\geq n_0$.\\

\noindent \textbf{Step 2} (Decomposition) Let $j_0=\min J_n$ and $j^*=\min\{j\in J_n: H_j\ \text{occurs}\}$, with $\min\emptyset=\infty$. If $j^*=j<\infty$ and $S_{\theta^{j-1}}\geq -K\theta^{j/2}$, then
\begin{align*}
    S_{\theta^j}=\Delta_j+S_{\theta^{j-1}}>2K\theta^{j/2}-K\theta^{j/2}=K\sqrt{\theta^j},\qquad \theta^j\in[a_n',n]\subset[a_n,n].
\end{align*}
Hence
\begin{align}
\label{eq proof of lemma A.4 1}
    \Big\{\forall\,k\in[a_n,n]:\ S_k\leq K\sqrt{k}\Big\}\subset\{j^*=\infty\}\cup\bigcup_{j\in J_n}\Big(\{j^*=j\}\cap\big\{S_{\theta^{j-1}}<-K\theta^{j/2}\big\}\Big).
\end{align}

\noindent \textbf{Step 3} (Bounds) Let $n\geq n_0$. By Step 1, $\P(j^*=\infty)=\prod_{j\in J_n}(1-\P(H_j))\leq(1-q)^{|J_n|}$. For $j\in J_n$, we have $\{j^*=j\}=H_j\cap C_j$ with $C_j=\bigcap_{i\in J_n,\,i<j}H_i^c\in\sigma(\varepsilon_1,\dots,\varepsilon_{\theta^{j-1}})$. Thus, Step 1 and Chebyshev's inequality yield
\begin{align*}
    \P\Big(j^*=j,\ S_{\theta^{j-1}}<-K\theta^{j/2}\Big)&=\P(H_j)\,\P\Big(C_j\cap\big\{S_{\theta^{j-1}}<-K\theta^{j/2}\big\}\Big)\\
    &\leq\min\Big\{(1-q)^{j-j_0},\ \frac{\theta^{j-1}}{K^2\theta^j}\Big\}=\min\big\{(1-q)^{j-j_0},\kappa\big\}.
\end{align*}
Let $i_0=\lceil\log(1/\kappa)/q\rceil$. Since $(1-q)^i\leq e^{-qi}\leq\kappa$ for $i\geq i_0$,
\begin{align*}
    \sum_{j\in J_n}\min\big\{(1-q)^{j-j_0},\kappa\big\}\leq\sum_{i=0}^\infty\min\big\{(1-q)^i,\kappa\big\}\leq i_0\kappa+\sum_{i\geq i_0}(1-q)^i\leq\kappa\Big(1+\frac{1+\log(1/\kappa)}{q}\Big).
\end{align*}

\noindent \textbf{Step 4} (Conclusion) By \eqref{eq proof of lemma A.4 1} and Step 3, for every $n\geq n_0$,
\begin{align*}
    \P\Big(\forall\,k\in[a_n,n]:\ S_k\leq K\sqrt{k}\Big)\leq(1-q)^{|J_n|}+\frac{1}{K^2\theta}\Big(1+\frac{1+\log(K^2\theta)}{q}\Big).
\end{align*}
Since $|J_n|\to\infty$, the $\limsup_n$ of the left-hand side is no larger than the second term on the right-hand side. As $q$ does not depend on $\theta$, letting $\theta\to\infty$ finishes the proof.\\

\noindent \textbf{Proof of Lemma \ref{lemma A.5}} We first show that
\begin{align}
\label{eq proof of lemma A.5 1}
    e^{x-\frac{x^2}{6}}\leq 1+x+\frac{x^2}{3},\qquad x\in\mathbb{R}.
\end{align}
Since $1+x+\frac{x^2}{3}>0$ for every $x\in\mathbb{R}$, the function $g(x)=\log(1+x+\frac{x^2}{3})-x+\frac{x^2}{6}$ is well defined, and a direct computation yields
\begin{align*}
    g'(x)=\frac{1+\frac{2x}{3}}{1+x+\frac{x^2}{3}}-1+\frac{x}{3}=\frac{x^3}{9\big(1+x+\frac{x^2}{3}\big)}.
\end{align*}
Hence $g$ is non-increasing on $(-\infty,0]$ and non-decreasing on $[0,\infty)$, so that $g(x)\geq g(0)=0$, which is \eqref{eq proof of lemma A.5 1}. Applying \eqref{eq proof of lemma A.5 1} with $x=\lambda Z$ and using $\E[Z|\mathcal{A}]=0$ and $1+y\leq e^y$,
\begin{align*}
    \E\Big[e^{\lambda Z-\frac{\lambda^2}{6}Z^2}\Big|\mathcal{A}\Big]\leq 1+\frac{\lambda^2}{3}\E[Z^2|\mathcal{A}]\leq\exp\Big(\frac{\lambda^2}{3}\E[Z^2|\mathcal{A}]\Big),
\end{align*}
where the left-hand side is finite because $x-\frac{x^2}{6}\leq\frac{3}{2}$. Since $\E[Z^2|\mathcal{A}]$ is measurable with respect to $\mathcal{A}$, the claim follows.\\

\noindent \textbf{Proof of Lemma \ref{lemma A.6}} Since $f$ is of bounded variation, it has at most countably many discontinuities, and the limits $f(x+)$, $f(x-)$ and $f(+\infty)$ exist. Let $\mu_2$ be the purely atomic signed measure with mass $f(y+)-f(y)$ at every point $y$ at which $f$ is not right-continuous, and let $\mu_1$ be the signed measure such that $\mu_1+\mu_2$ is the Lebesgue--Stieltjes measure of the right-continuous function $x\mapsto f(x+)$. Then $\mu_1$ consists of the continuous part of $df$ and of the left jumps $f(y)-f(y-)$, while $\mu_2$ consists of the right jumps $f(y+)-f(y)$, so that $|\mu_1|(\mathbb{R})+|\mu_2|(\mathbb{R})=\|f\|_{TV}$. Since $f(+\infty)-f(x+)=(\mu_1+\mu_2)((x,\infty))$ and $f(x)=f(x+)-\mu_2(\{x\})$, we obtain
\begin{align}
\label{eq proof of lemma A.6 1}
    f(x)=f(+\infty)-A(x),\qquad A(x)=\int1\{x<t\}\,\mu_1(dt)+\int1\{x\leq t\}\,\mu_2(dt),\qquad x\in\mathbb{R}.
\end{align}

\noindent \textbf{Step 1} (Proof of (i)) Let $\P_{X|\mathcal{M}}$ be a regular conditional distribution of $X$ given $\mathcal{M}$ and denote, for $t\in\mathbb{R}$,
\begin{align*}
    \Delta_{\leq}(t)=\P_{X|\mathcal{M}}\big((-\infty,t]\big)-\P(X\leq t),\qquad \Delta_{<}(t)=\P_{X|\mathcal{M}}\big((-\infty,t)\big)-\P(X<t).
\end{align*}
By \eqref{eq proof of lemma A.6 1} and Fubini's theorem,
\begin{align*}
    \E[f(X)|\mathcal{M}]-\E[f(X)]=-\int\Delta_{<}(t)\,\mu_1(dt)-\int\Delta_{\leq}(t)\,\mu_2(dt)\ \ \text{a.s.}
\end{align*}
By the definition of $\tilde{\phi}(\mathcal{M},X)$, there exists a null set $N$ such that $|\Delta_{\leq}(q)|\leq\tilde{\phi}(\mathcal{M},X)$ holds on $N^c$ for every rational $q$. On $N^c$, both $t\mapsto\P_{X|\mathcal{M}}((-\infty,t])$ and $t\mapsto\P(X\leq t)$ are distribution functions. Hence, by right-continuity, $|\Delta_{\leq}(t)|\leq\tilde{\phi}(\mathcal{M},X)$ holds on $N^c$ for every $t\in\mathbb{R}$, and, by taking left limits, so does $|\Delta_{<}(t)|\leq\tilde{\phi}(\mathcal{M},X)$. Therefore, on $N^c$,
\begin{align*}
    \big|\E[f(X)|\mathcal{M}]-\E[f(X)]\big|\leq\big(|\mu_1|(\mathbb{R})+|\mu_2|(\mathbb{R})\big)\tilde{\phi}(\mathcal{M},X)=\|f\|_{TV}\,\tilde{\phi}(\mathcal{M},X).
\end{align*}
The same argument shows that $|\E[A(X)|\mathcal{M}]-\E[A(X)]|\leq\|f\|_{TV}\,\tilde{\phi}(\mathcal{M},X)$ almost surely.\\

\noindent \textbf{Step 2} (Proof of (ii)) By \eqref{eq proof of lemma A.6 1},
\begin{align*}
    f(X)f(Y)=f(+\infty)^2-f(+\infty)\big(A(X)+A(Y)\big)+A(X)A(Y).
\end{align*}
Since $|f(+\infty)|\leq\|f\|_\infty$, Step 1 shows that the conditional expectation of the middle term deviates from its expectation by at most $\|f\|_\infty\|f\|_{TV}(\tilde{\phi}(\mathcal{M},X)+\tilde{\phi}(\mathcal{M},Y))$ almost surely. For the last term, let $\P_{(X,Y)|\mathcal{M}}$ be a regular conditional distribution of $(X,Y)$ given $\mathcal{M}$. By \eqref{eq proof of lemma A.6 1} and Fubini's theorem, $\E[A(X)A(Y)|\mathcal{M}]-\E[A(X)A(Y)]$ is a sum of four integrals with respect to the product measures $\mu_a\otimes\mu_{a'}$, $a,a'\in\{1,2\}$, whose integrands are of the form $\P_{(X,Y)|\mathcal{M}}(I\times J)-\P((X,Y)\in I\times J)$ with $I\in\{(-\infty,s),(-\infty,s]\}$ and $J\in\{(-\infty,t),(-\infty,t]\}$. As in Step 1, outside a null set, these integrands are bounded in absolute value by $\tilde{\phi}(\mathcal{M},(X,Y))$ for $I=(-\infty,s]$ and $J=(-\infty,t]$ with rational $s,t$, hence for all real $s,t$ by continuity from above, and for all four types of $I\times J$ by continuity from below. Therefore,
\begin{align*}
    \big|\E[A(X)A(Y)|\mathcal{M}]-\E[A(X)A(Y)]\big|&\leq\big(|\mu_1|(\mathbb{R})+|\mu_2|(\mathbb{R})\big)^2\tilde{\phi}(\mathcal{M},(X,Y))\\
    &=\|f\|_{TV}^2\,\tilde{\phi}(\mathcal{M},(X,Y))\ \ \text{a.s.},
\end{align*}
and (ii) follows.\\

\noindent \textbf{Step 3} (Smaller $\sigma$-algebras) Let $\mathcal{M}'\subset\mathcal{M}$ and let $W$ be a bounded random variable. Then $\E[W|\mathcal{M}']-\E[W]=\E\big[\E[W|\mathcal{M}]-\E[W]\big|\mathcal{M}'\big]$, and $|\E[V|\mathcal{M}']|\leq\|V\|_\infty$ almost surely for every bounded $V$. Applying this to $W=f(X)$ and $W=f(X)f(Y)$ proves the last claim.\\

\noindent \textbf{Proof of Lemma \ref{lemma A.7}} (i) Denote $\mathcal{G}=\sigma(Z_{nk}:k\leq 0)$ and $\mathcal{H}_i=\sigma(Z_{n1},\dots,Z_{ni})$, so that $\mathcal{F}_{ni}=\mathcal{G}\vee\mathcal{H}_i$ for $0\leq i\leq n$. By \ding{173}, $\mathcal{G}$ is independent of $\sigma(Z_{n1},\dots,Z_{nn})$, which contains $\sigma(g_n)\vee\mathcal{H}_i$. Hence, $\E[g_n|\mathcal{F}_{ni}]=\E[g_n|\mathcal{H}_i]$ almost surely. Next, let $w$ be a measurable function such that $\E[g_0|\mathcal{F}_1^i]=w(Z_1,\dots,Z_i)$ almost surely, where $\mathcal{F}_1^0=\{\emptyset,\Omega\}$ and $w$ is constant for $i=0$. For every bounded measurable function $u$, \ding{172} yields
\begin{align*}
    \E\big[g_nu(Z_{n1},\dots,Z_{ni})\big]&=\E\big[g_0u(Z_1,\dots,Z_i)\big]=\E\big[w(Z_1,\dots,Z_i)u(Z_1,\dots,Z_i)\big]\\
    &=\E\big[w(Z_{n1},\dots,Z_{ni})u(Z_{n1},\dots,Z_{ni})\big],
\end{align*}
so that $\E[g_n|\mathcal{H}_i]=w(Z_{n1},\dots,Z_{ni})$ almost surely. Since $(Z_{n1},\dots,Z_{ni})$ and $(Z_1,\dots,Z_i)$ have the same distribution and $\E[g_n]=\E[g_0]$, we obtain, for $s\in[1,\infty]$,
\begin{align*}
    \big\|\E[g_n|\mathcal{F}_{ni}]-\E[g_n]\big\|_s=\big\|\E[g_0|\mathcal{F}_1^i]-\E[g_0]\big\|_s\leq\big\|\E[g_0|\mathcal{F}_{-\infty}^i]-\E[g_0]\big\|_s,
\end{align*}
where the inequality follows from $\E[g_0|\mathcal{F}_1^i]-\E[g_0]=\E\big[\E[g_0|\mathcal{F}_{-\infty}^i]-\E[g_0]\big|\mathcal{F}_1^i\big]$ and the contraction of conditional expectations in $L^s$, because $\mathcal{F}_1^i\subset\mathcal{F}_{-\infty}^i$.\\
(ii) Let $0\leq i\leq n$. For $k\geq n+1$, \ding{173} implies that $Z_{nk}$ is independent of $\mathcal{F}_{ni}\subset\mathcal{G}\vee\sigma(Z_{n1},\dots,Z_{nn})$, so that \ding{174} and \ding{172} yield $\E[h_n(Z_{nk})-\mu_n|\mathcal{F}_{ni}]=\E[h_n(Z_{nn})]-\mu_n=\E[h_n(Z_n)]-\mu_n=0$. Hence $\varepsilon_{ni}=\sum_{k=i+1}^n\E[Y_{nk}|\mathcal{F}_{ni}]$ and, in particular, $\varepsilon_{nn}=0$. By (i) with $i=0$ and by \ding{172}, $\E[Y_{nk}|\mathcal{F}_{n0}]=\E[Y_{nk}]=\E[h_n(Z_k)]-\mu_n=0$ for $1\leq k\leq n$, so that $\varepsilon_{n0}=0$. Since $Y_{nk}$ is measurable with respect to $\mathcal{F}_{ni}$ for $k\leq i$, we have $\E[S_n|\mathcal{F}_{ni}]=\sum_{k=1}^iY_{nk}+\varepsilon_{ni}$. Moreover, for $1\leq k<i\leq n$, $h_n(Z_{nk})$ is measurable with respect to $\mathcal{F}_{n,i-1}$, so that $D_{ni}=\E[S_n|\mathcal{F}_{ni}]-\E[S_n|\mathcal{F}_{n,i-1}]$. Therefore,
\begin{align*}
    M_{nj}=\E\big[S_n\big|\mathcal{F}_{n,jl}\big]-\E\big[S_n\big|\mathcal{F}_{n,(j-1)l}\big]=T_{nj}+\varepsilon_{n,jl}-\varepsilon_{n,(j-1)l},\qquad 1\leq j\leq m,
\end{align*}
which is measurable with respect to $\mathcal{F}_{n,jl}$ and has zero conditional mean given $\mathcal{F}_{n,(j-1)l}$ by the tower property. Finally, $\sum_{j=1}^mM_{nj}=\E[S_n|\mathcal{F}_{nn}]-\E[S_n|\mathcal{F}_{n0}]=(S_n+\varepsilon_{nn})-\varepsilon_{n0}=S_n$.\\
(iii) For $0\leq i<k\leq n$, apply (i) to $g(z_1,\dots,z_n)=h_n(z_k)-\mu_n$, and then Lemma \ref{lemma A.6} (i) to $f=h_n-\mu_n$ with $\mathcal{M}=\mathcal{F}_{-\infty}^i$. Since $\E[f(Z_k)]=0$, $\|f\|_{TV}=\|h_n\|_{TV}$ and $\tilde{\phi}(\mathcal{F}_{-\infty}^i,Z_k)\leq\tilde{\phi}_1(k-i)$ by \eqref{eq def phi tilde k}, we obtain $|\E[Y_{nk}|\mathcal{F}_{ni}]|\leq\|h_n\|_{TV}\tilde{\phi}_1(k-i)$ almost surely. By (ii), summing over $k$ yields $|\varepsilon_{ni}|\leq\|h_n\|_{TV}\sum_{r=1}^{n-i}\tilde{\phi}_1(r)\leq\|h_n\|_{TV}\tilde{\Phi}_n$. Moreover, $|Y_{ni}|\leq b_n-a_n$, because $h_n(Z_{ni})$ and $\mu_n$ both lie in $[a_n,b_n]$. Hence $|T_{nj}|\leq l(b_n-a_n)$ and, by (ii), $|M_{nj}|\leq l(b_n-a_n)+2\|h_n\|_{TV}\tilde{\Phi}_n=\bar{b}_n$.\\
(iv) Write $T_{nj}=\sum_{i\in B_j}(h_n(Z_{ni})-\P_n(h_n))+l(\P_n(h_n)-\mu_n)$. Since $n=ml$, the first sums add up to zero over $1\leq j\leq m$. Hence, by \eqref{eq V hat},
\begin{align*}
    \sum_{j=1}^mT_{nj}^2=\sum_{j=1}^m\Big(\sum_{i\in B_j}\big(h_n(Z_{ni})-\P_n(h_n)\big)\Big)^2+ml^2\big(\P_n(h_n)-\mu_n\big)^2=n\hat{V}_n+nl\big(\P_n(h_n)-\mu_n\big)^2.
\end{align*}

\noindent \textbf{Proof of Lemma \ref{lemma A.8}} By the first statement of Lemma \ref{lemma A.7} (i), $\varepsilon_{nl}$ has a version that is measurable with respect to $\sigma(Z_{n1},\dots,Z_{nl})$. Hence, by Lemma \ref{lemma A.7} (ii) and (iii), $M_{n1}$ is almost surely equal to a bounded measurable function of $(Z_{n1},\dots,Z_{nn})$, and Lemma \ref{lemma A.7} (i) with $i=0$ yields $\E[M_{n1}^2|\mathcal{F}_{n0}]=\E[M_{n1}^2]$. Thus there is nothing to prove for $j=1$. Fix $2\leq j\leq m$ and denote $p=(j-1)l$ and $\mathcal{M}'=\mathcal{F}_{np}$.\\

\noindent \textbf{Step 1} (Decomposition) By Lemma \ref{lemma A.7} (ii), $\E[M_{nj}|\mathcal{M}']=0$. Since $\varepsilon_{np}$ is measurable with respect to $\mathcal{M}'$, this means $\varepsilon_{np}=\E[T_{nj}+\varepsilon_{n,jl}|\mathcal{M}']$, so that
\begin{align*}
    \E\big[M_{nj}^2\big|\mathcal{M}'\big]=\E\big[(T_{nj}+\varepsilon_{n,jl})^2\big|\mathcal{M}'\big]-\varepsilon_{np}^2,\qquad \E\big[M_{nj}^2\big]=\E\big[(T_{nj}+\varepsilon_{n,jl})^2\big]-\E\big[\varepsilon_{np}^2\big].
\end{align*}
Hence $\E[M_{nj}^2|\mathcal{M}']-\E[M_{nj}^2]=\Delta_1+2\Delta_2+\Delta_3-\Delta_4$, where
\begin{align*}
    \Delta_1&=\E\big[T_{nj}^2\big|\mathcal{M}'\big]-\E\big[T_{nj}^2\big],\qquad \Delta_2=\E\big[T_{nj}\varepsilon_{n,jl}\big|\mathcal{M}'\big]-\E\big[T_{nj}\varepsilon_{n,jl}\big],\\
    \Delta_3&=\E\big[\varepsilon_{n,jl}^2\big|\mathcal{M}'\big]-\E\big[\varepsilon_{n,jl}^2\big],\qquad \Delta_4=\varepsilon_{np}^2-\E\big[\varepsilon_{np}^2\big].
\end{align*}
By Lemma \ref{lemma A.7} (iii), $\varepsilon_{n,jl}^2$ and $\varepsilon_{np}^2$ take values in $[0,\|h_n\|_{TV}^2\tilde{\Phi}_n^2]$, so that $|\Delta_3|+|\Delta_4|\leq 2\|h_n\|_{TV}^2\tilde{\Phi}_n^2$ almost surely. Moreover, since $Y_{ni}$, $i\in B_j$, is measurable with respect to $\mathcal{F}_{n,jl}\supset\mathcal{M}'$, Lemma \ref{lemma A.7} (ii) and the tower property yield $\E[T_{nj}\varepsilon_{n,jl}|\mathcal{M}']=\sum_{i\in B_j}\sum_{k=jl+1}^n\E[Y_{ni}Y_{nk}|\mathcal{M}']$. Denoting $\Delta(i,k)=\E[Y_{ni}Y_{nk}|\mathcal{M}']-\E[Y_{ni}Y_{nk}]$, we thus obtain
\begin{align*}
    \Delta_1+2\Delta_2=\sum_{i\in B_j}\Big(\Delta(i,i)+2\sum_{k=i+1}^n\Delta(i,k)\Big).
\end{align*}

\noindent \textbf{Step 2} (Bounds on $\Delta(i,k)$) Let $f=h_n-\mu_n$, so that $\|f\|_{TV}=\|h_n\|_{TV}$ and $\|f\|_\infty\leq b_n-a_n$, and let $i=p+d\in B_j$ with $1\leq d\leq l$. We apply Lemma \ref{lemma A.7} (i) with $i$ replaced by $p$, and then Lemma \ref{lemma A.6} with $\mathcal{M}=\mathcal{F}_{-\infty}^p$, together with \eqref{eq def phi tilde k}.
\begin{enumerate}[(a)]
    \item Since $|f(x)^2-f(y)^2|\leq 2\|f\|_\infty|f(x)-f(y)|$, the function $f^2$ is of bounded variation with $\|f^2\|_{TV}\leq 2(b_n-a_n)\|h_n\|_{TV}$. Lemma \ref{lemma A.6} (i) applied to $f^2$ yields $|\Delta(i,i)|\leq 2(b_n-a_n)\|h_n\|_{TV}\tilde{\phi}_1(d)$.
    \item For $i<k\leq n$, the indices of $Z_i$, $Z_k$ and $(Z_i,Z_k)$ are no smaller than $p+d$. Lemma \ref{lemma A.6} (ii) yields $|\Delta(i,k)|\leq A_n\tilde{\phi}_2(d)$ with $A_n=\|h_n\|_{TV}^2+2(b_n-a_n)\|h_n\|_{TV}$.
    \item For $i<k\leq n$, we have $\mathcal{M}'\subset\mathcal{F}_{ni}$, and the tower property yields
    \begin{align*}
        \E[Y_{ni}Y_{nk}|\mathcal{M}']=\E\big[Y_{ni}\E[Y_{nk}|\mathcal{F}_{ni}]\big|\mathcal{M}'\big].
    \end{align*}
    As in the proof of Lemma \ref{lemma A.7} (iii), the inner conditional expectation is bounded by $\|h_n\|_{TV}\tilde{\phi}_1(k-i)$ in absolute value. Hence the conditional and the unconditional expectations of $Y_{ni}Y_{nk}$ are both bounded by $(b_n-a_n)\|h_n\|_{TV}\tilde{\phi}_1(k-i)$ in absolute value, so that $|\Delta(i,k)|\leq 2(b_n-a_n)\|h_n\|_{TV}\tilde{\phi}_1(k-i)$.
\end{enumerate}

\noindent \textbf{Step 3} (Summation) By \eqref{eq def phi tilde k}, $\tilde{\phi}_1\leq\tilde{\phi}_2$, and $\tilde{\phi}_2$ is non-increasing. Since $2(b_n-a_n)\|h_n\|_{TV}\leq A_n$, (b) and (c) yield $|\Delta(i,k)|\leq A_n\tilde{\phi}_2(\max\{d,k-i\})$ for $i<k\leq n$. For every $r\geq 1$, at most $2r-1$ pairs $(d,e)\in\{1,\dots,l\}\times\{1,\dots,n-1\}$ satisfy $\max\{d,e\}=r$, and $\max\{d,e\}\leq n$ for all of them. Hence
\begin{align*}
    \sum_{i\in B_j}\sum_{k=i+1}^n|\Delta(i,k)|\leq A_n\sum_{r=1}^n(2r-1)\tilde{\phi}_2(r)\leq 2A_n\Psi_n.
\end{align*}
Together with (a) and $\sum_{d=1}^l\tilde{\phi}_1(d)\leq\tilde{\Phi}_n\leq\Psi_n$, we obtain
\begin{align*}
    |\Delta_1+2\Delta_2|\leq 2(b_n-a_n)\|h_n\|_{TV}\Psi_n+4A_n\Psi_n=\big(4\|h_n\|_{TV}^2+10(b_n-a_n)\|h_n\|_{TV}\big)\Psi_n.
\end{align*}
Combined with Step 1, this yields $|\E[M_{nj}^2|\mathcal{M}']-\E[M_{nj}^2]|\leq\Gamma_n/(2m)$ almost surely.\\

\noindent \textbf{Step 4} (Consequence) Since $\E[S_n]=0$, Lemma \ref{lemma A.7} (ii) and the orthogonality of martingale differences yield $\sum_{j=1}^m\E[M_{nj}^2]=\E[S_n^2]=\text{Var}(S_n)$. Summing the bound of Step 3 over $1\leq j\leq m$ yields the last claim.\\

\noindent \textbf{Proof of Lemma \ref{lemma A.9}} Denote $L=\|h_n\|_{\text{Lip}}$, $\kappa=\kappa_n$, $q'=q/(q-1)$ and $f=h_n-\mu_n$, so that $Y_{ni}=f(Z_{ni})$, $\E[f(Z_i)]=0$ and $\|f(Z_i)\|_q\leq\kappa$ for $1\leq i\leq n$. By \ding{172}, $Z_{ni}$ and $Z_i$ have the same distribution, so that every moment bound below that concerns $f(Z_i)$ also holds for $f(Z_{ni})$. If $\kappa=0$ or $L=0$, then $Y_{ni}=0$ almost surely for $1\leq i\leq n$ (for $L=0$, $h_n$ is constant and equal to $\mu_n$), and (i)--(iv) are trivial. We therefore assume $\kappa>0$ and $L>0$. Since $\tau$-mixing only requires $\E|Z_i|<\infty$ for every $i$, the coefficients $\tau_1(r)$ and $\tau_2(r)$ may be infinite. If $\tau_2(r)=\infty$ for some $r\leq n$, then $\Psi_n^\tau=\infty$ and the bounds involving $\Psi_n^\tau$, namely the first bound in (i), (ii) and (iii), are trivial. If $\tau_1(r)=\infty$ for some $r\leq n$, then also $\Phi_n^\tau=\infty$, and all bounds are trivial. In the steps below, we may therefore assume that every $\tau$-coefficient entering a bound is finite. For $s\in[1,q]$, let $\theta_s=(q-s)/(s(q-1))$, so that $1/s=\theta_s+(1-\theta_s)/q$. Then $\theta_{q'}=\beta$ and $\theta_4=\gamma$. By H\"older's inequality, $\|W\|_s\leq\|W\|_1^{\theta_s}\|W\|_q^{1-\theta_s}$ for every random variable $W$. Moreover, $\kappa^{2-\beta}L^\beta=c_n$, because $2-\beta=q/(q-1)$.\\

\noindent \textbf{Step 1} (Single conditional expectations) Let $0\leq i<k\leq n$ and $g_{i,k}=\E[Y_{nk}|\mathcal{F}_{ni}]$. Lemma \ref{lemma A.7} (i) with $s=1$ yields $\|g_{i,k}\|_1\leq\|\E[f(Z_k)|\mathcal{F}_{-\infty}^i]\|_1$. Since $f(Z_k)$ is integrable, $\E[f(Z_k)|\mathcal{F}_{-\infty}^i]=\int f\,d\P_{Z_k|\mathcal{F}_{-\infty}^i}$ almost surely, and $\E[f(Z_k)]=\int f\,d\P_{Z_k}=0$. Since $f/L\in\Lambda_1(\mathbb{R})$, Definition \ref{def tau k} yields $\|\E[f(Z_k)|\mathcal{F}_{-\infty}^i]\|_1\leq L\tau(\mathcal{F}_{-\infty}^i,Z_k)\leq L\tau_1(k-i)$. Moreover, $\|g_{i,k}\|_q\leq\|Y_{nk}\|_q\leq\kappa$ by Jensen's inequality for conditional expectations. Hence
\begin{align}
\label{eq proof of lemma A.9 1}
    \|g_{i,k}\|_{q'}\leq\kappa^{1-\beta}\big(L\tau_1(k-i)\big)^\beta,\qquad \|g_{i,k}\|_4\leq\kappa^{1-\gamma}\big(L\tau_1(k-i)\big)^\gamma.
\end{align}

\noindent \textbf{Step 2} (Proof of (i) and (iii)) By Lemma \ref{lemma A.7} (ii), $\varepsilon_{ni}=\sum_{k=i+1}^ng_{i,k}$. By Minkowski's inequality and \eqref{eq proof of lemma A.9 1}, $\|\varepsilon_{ni}\|_4\leq\kappa^{1-\gamma}L^\gamma\sum_{r=1}^n\tau_1(r)^\gamma=\kappa^{1-\gamma}L^\gamma\Phi_n^\tau$. For $i<k\leq k'\leq n$, H\"older's inequality and \eqref{eq proof of lemma A.9 1} yield $|\E[g_{i,k}g_{i,k'}]|\leq\|g_{i,k}\|_q\|g_{i,k'}\|_{q'}\leq c_n\tau_1(k'-i)^\beta$. For every $r\geq 1$, at most $2r-1$ pairs $(k,k')$ satisfy $\max\{k,k'\}-i=r$. Since $\tau_1\leq\tau_2$, we obtain
\begin{align*}
    \E[\varepsilon_{ni}^2]=\sum_{k,k'=i+1}^n\E[g_{i,k}g_{i,k'}]\leq c_n\sum_{r=1}^n(2r-1)\tau_1(r)^\beta\leq 2c_n\Psi_n^\tau.
\end{align*}
This proves (i), and (iii) follows from $(a-b)^2\leq 2a^2+2b^2$.\\

\noindent \textbf{Step 3} (Products) Fix $2\leq j\leq m$, denote $p=(j-1)l$ and $\mathcal{M}'=\mathcal{F}_{np}$, and, for $i\in B_j$ and $i\leq k\leq n$, let $\Delta(i,k)=\E[Y_{ni}Y_{nk}|\mathcal{M}']-\E[Y_{ni}Y_{nk}]$ and $d=i-p$. For $K>0$, let $f_K=\max\{-K,\min\{K,f\}\}$, which is $L$-Lipschitz and bounded by $K$, and note that $|f-f_K|\leq|f|1\{|f|>K\}$ and $|f_K|\leq|f|$. Since $|f|^s1\{|f|>K\}\leq|f|^qK^{s-q}$ for $s\leq q$, we have, for $1\leq i\leq n$ and, by \ding{172}, also with $Z_i$ replaced by $Z_{ni}$,
\begin{align}
\label{eq proof of lemma A.9 2}
    \E\big[f(Z_i)^21\{|f(Z_i)|>K\}\big]\leq\kappa^qK^{2-q},\qquad \big\|f(Z_i)1\{|f(Z_i)|>K\}\big\|_{q'}\leq\kappa^{q-1}K^{2-q}.
\end{align}
Below, we use Lemma \ref{lemma A.7} (i) with $s=1$ and $i$ replaced by $p$, together with Definition \ref{def tau k}. The bounds in (a) and (b) below hold for every $K>0$ and have the form $C_1KL\tau_0+C_2\kappa^qK^{2-q}$ for a finite $\tau$-coefficient $\tau_0\geq 0$ and constants $C_1,C_2>0$. If $\tau_0>0$, we choose $K=\kappa(\kappa/(L\tau_0))^{1/(q-1)}$, for which $KL\tau_0=\kappa^qK^{2-q}=c_n\tau_0^\beta$, so that the bound equals $(C_1+C_2)c_n\tau_0^\beta$. If $\tau_0=0$, the infimum over $K>0$ is zero, and the same conclusion holds.
\begin{enumerate}[(a)]
    \item The function $f_K^2$ is $2KL$-Lipschitz, so that the conditional expectation of $f_K(Z_{ni})^2$ deviates from its mean by at most $2KL\tau(\mathcal{F}_{-\infty}^p,Z_i)\leq 2KL\tau_1(d)$ in $L^1$. Since $0\leq f^2-f_K^2\leq f^21\{|f|>K\}$, \eqref{eq proof of lemma A.9 2} yields $\|\Delta(i,i)\|_1\leq 2KL\tau_1(d)+2\kappa^qK^{2-q}$, and the choice $\tau_0=\tau_1(d)$ yields $\|\Delta(i,i)\|_1\leq 4c_n\tau_1(d)^\beta$.
    \item Let $i<k$. The function $(x,y)\mapsto f_K(x)f_K(y)$ is $KL$-Lipschitz with respect to $\|\cdot\|_1$, so that the corresponding part of $\Delta(i,k)$ is bounded by $KL\tau(\mathcal{F}_{-\infty}^p,(Z_i,Z_k))\leq 2KL\tau_2(d)$ in $L^1$. Since $|f(x)f(y)-f_K(x)f_K(y)|\leq|f(x)|1\{|f(x)|>K\}|f(y)|+|f(x)||f(y)|1\{|f(y)|>K\}$, H\"older's inequality and \eqref{eq proof of lemma A.9 2} bound the remaining part by $4\kappa^qK^{2-q}$ in $L^1$. The choice $\tau_0=2\tau_2(d)$ yields $\|\Delta(i,k)\|_1\leq 5c_n(2\tau_2(d))^\beta\leq 10c_n\tau_2(d)^\beta$.
    \item Let $i<k$. Since $\mathcal{M}'\subset\mathcal{F}_{ni}$, $\E[Y_{ni}Y_{nk}|\mathcal{M}']=\E[Y_{ni}g_{i,k}|\mathcal{M}']$ with $g_{i,k}$ as in Step 1. H\"older's inequality and \eqref{eq proof of lemma A.9 1} yield $\E|Y_{ni}g_{i,k}|\leq\kappa\|g_{i,k}\|_{q'}\leq c_n\tau_1(k-i)^\beta$, so that $\|\Delta(i,k)\|_1\leq 2c_n\tau_1(k-i)^\beta$.
\end{enumerate}
Since $\tau_1\leq\tau_2$ and $\tau_2$ is non-increasing, (b) and (c) yield $\|\Delta(i,k)\|_1\leq 10c_n\tau_2(\max\{d,k-i\})^\beta$ for $i<k\leq n$.\\

\noindent \textbf{Step 4} (Proof of (ii)) For $j=1$, the first statement of Lemma \ref{lemma A.7} (i) shows that $\varepsilon_{nl}$ has a version that is measurable with respect to $\sigma(Z_{n1},\dots,Z_{nl})$, so that, by Lemma \ref{lemma A.7} (ii), $M_{n1}$ is almost surely equal to a measurable function of $(Z_{n1},\dots,Z_{nn})$. Moreover, $M_{n1}\in L^2$ by (i) and $\|T_{n1}\|_2\leq l\kappa$. Hence Lemma \ref{lemma A.7} (i) with $i=0$ yields $\E[M_{n1}^2|\mathcal{F}_{n0}]=\E[M_{n1}^2]$. For $2\leq j\leq m$, we use the decomposition derived in Step 1 of the proof of Lemma \ref{lemma A.8}, but not the bound on $\Delta_3$ and $\Delta_4$ given there, which requires $h_n$ to be bounded. The derivation of the decomposition only uses $\E[M_{nj}|\mathcal{M}']=0$, the $\mathcal{M}'$-measurability of $\varepsilon_{np}$ and the tower property, which apply because $Y_{ni}\in L^q$ and $\varepsilon_{ni}\in L^2$ by (i). It yields $\E[M_{nj}^2|\mathcal{M}']-\E[M_{nj}^2]=\Delta_1+2\Delta_2+\Delta_3-\Delta_4$ with
\begin{align*}
    \Delta_1+2\Delta_2=\sum_{i\in B_j}\Big(\Delta(i,i)+2\sum_{k=i+1}^n\Delta(i,k)\Big),\qquad \|\Delta_3\|_1+\|\Delta_4\|_1\leq 2\E[\varepsilon_{n,jl}^2]+2\E[\varepsilon_{np}^2]\leq 8c_n\Psi_n^\tau,
\end{align*}
where the last inequality follows from (i). By Step 3 and the counting argument of Step 3 of the proof of Lemma \ref{lemma A.8},
\begin{align*}
    \|\Delta_1+2\Delta_2\|_1\leq 4c_n\sum_{d=1}^l\tau_1(d)^\beta+20c_n\sum_{r=1}^n(2r-1)\tau_2(r)^\beta\leq 44c_n\Psi_n^\tau.
\end{align*}
Hence $\|\E[M_{nj}^2|\mathcal{F}_{n,(j-1)l}]-\E[M_{nj}^2]\|_1\leq 52c_n\Psi_n^\tau$ for every $1\leq j\leq m$. Since $\sum_{j=1}^m\E[M_{nj}^2]=\text{Var}(S_n)$ as in Step 4 of the proof of Lemma \ref{lemma A.8}, (ii) follows.\\

\noindent \textbf{Step 5} (Proof of (iv)) By Minkowski's inequality, $\|T_{nj}\|_4\leq l\max_i\|Y_{ni}\|_4\leq l\kappa$, because $q>4$. Since $(a+b)^4\leq 8a^4+8b^4$, Lemma \ref{lemma A.7} (ii) and statement (i) of the present lemma yield
\begin{align*}
    \E[M_{nj}^4]\leq 8\E[T_{nj}^4]+64\big(\E[\varepsilon_{n,jl}^4]+\E[\varepsilon_{n,(j-1)l}^4]\big)\leq 8l^4\kappa^4+128\,\kappa^{4-4\gamma}L^{4\gamma}(\Phi_n^\tau)^4.
\end{align*}
Summing over $1\leq j\leq m$ and using $ml^4=nl^3$ proves (iv).\\

\subsection{Proof of Results in Section 2}

\noindent \textbf{Proof of Theorem \ref{th 1}} \\
\noindent \textbf{Step 0} (Preparation) Denote by $v$ the deterministic value of $\<\M\>_n$. Since $\E[Z_i|\mathcal{A}_{i-1}]=0$, we have $\E[Z_i]=0$ and hence
\begin{align}
\label{eq proof of th 1 0}
    v=\E[\<\M\>_n]=\sum_{i=1}^n\E[Z_i^2]=\sum_{i=1}^n\text{Var}(Z_i)=n\sigma_n^2.
\end{align}
If $v=0$, then $Z_i=0$ almost surely for every $i\leq n$ and \eqref{eq th 1 1} holds trivially. We thus assume $v>0$. Denote $V_i:=\E[Z_i^2|\mathcal{A}_{i-1}]-Z_i^2$. Then $V_i$ is measurable with respect to $\mathcal{A}_i$, $\E[V_i|\mathcal{A}_{i-1}]=0$ and $\<\M\>_n-[\M]_n=\sum_{i=1}^nV_i$.\\

\noindent \textbf{Step 1} (Deviation of $\M_n$) For $\lambda\in[0,3/b)$ and $0\leq k\leq n$, define
\begin{align*}
    L_k(\lambda)=\exp\Big(\lambda\M_k-\psi(\lambda)\<\M\>_k\Big),\qquad \M_0=\<\M\>_0=0,
\end{align*}
where $\psi$ is given in \eqref{eq psi}. Since $\E[Z_k^2|\mathcal{A}_{k-1}]$ is measurable with respect to $\mathcal{A}_{k-1}$, Lemma \ref{lemma A.1} yields
\begin{align*}
    \E[L_k(\lambda)|\mathcal{A}_{k-1}]=L_{k-1}(\lambda)\exp\big(-\psi(\lambda)\E[Z_k^2|\mathcal{A}_{k-1}]\big)\E[e^{\lambda Z_k}|\mathcal{A}_{k-1}]\leq L_{k-1}(\lambda)\ \ \text{a.s.}
\end{align*}
Thus, $(L_k(\lambda),\mathcal{A}_k)_{0\leq k\leq n}$ is a non-negative super-martingale and $\E[L_n(\lambda)]\leq L_0(\lambda)=1$. Since $\<\M\>_n=v$ almost surely, this means $\E[e^{\lambda\M_n}]\leq e^{\psi(\lambda)v}$. Letting $\lambda=\lambda(v,t_1)$ from Lemma \ref{lemma A.3}, Markov's inequality and \eqref{eq lemma A.3 1} yield
\begin{align*}
    \P\Big(\M_n\geq \sqrt{2vt_1}+\frac{b}{3}t_1\Big)\leq \exp\Big(-\lambda\Big(\sqrt{2vt_1}+\frac{b}{3}t_1\Big)+\psi(\lambda)v\Big)=e^{-t_1}.
\end{align*}
Since $(-Z_i,\mathcal{A}_i)_{i\geq 1}$ is an MDS bounded by $b$ with the same predictable quadratic variation, the same argument applies to $-\M_n$. Therefore,
\begin{align}
\label{event 1}
    \P(\mathbf{E}_1^{\pm})\geq 1-e^{-t_1},\quad \mathbf{E}_1^{\pm}:=\Big\{\pm\M_n<\sqrt{2vt_1}+\frac{b}{3}t_1\Big\}.
\end{align}

\noindent \textbf{Step 2} (Lower tail of $[\M]_n$) For $\lambda\geq 0$ and $0\leq k\leq n$, define
\begin{align*}
    K_k(\lambda)=\exp\Big(\lambda\big(\<\M\>_k-[\M]_k\big)-\frac{\lambda^2b^2}{2}\<\M\>_k\Big)=\exp\Big(\sum_{i=1}^k\Big(\lambda V_i-\frac{\lambda^2b^2}{2}\E[Z_i^2|\mathcal{A}_{i-1}]\Big)\Big).
\end{align*}
By Lemma \ref{lemma A.2},
\begin{align*}
    \E[K_k(\lambda)|\mathcal{A}_{k-1}]=K_{k-1}(\lambda)\exp\Big(\lambda\E[Z_k^2|\mathcal{A}_{k-1}]-\frac{\lambda^2b^2}{2}\E[Z_k^2|\mathcal{A}_{k-1}]\Big)\E[e^{-\lambda Z_k^2}|\mathcal{A}_{k-1}]\leq K_{k-1}(\lambda)\ \ \text{a.s.}
\end{align*}
Thus, $(K_k(\lambda),\mathcal{A}_k)_{0\leq k\leq n}$ is a non-negative super-martingale and, since $\<\M\>_n=v$ almost surely, $\E[e^{\lambda(v-[\M]_n)}]\leq e^{\lambda^2b^2v/2}$. Letting $\lambda=\sqrt{2t_2/(b^2v)}$, Markov's inequality yields
\begin{align*}
    \P\Big(v-[\M]_n\geq b\sqrt{2vt_2}\Big)\leq \exp\Big(-\lambda b\sqrt{2vt_2}+\frac{\lambda^2b^2v}{2}\Big)=e^{-2t_2+t_2}=e^{-t_2},
\end{align*}
which indicates
\begin{align}
\label{event 2}
    \P(\mathbf{E}_2)\geq 1-e^{-t_2},\quad \mathbf{E}_2:=\Big\{v<[\M]_n+b\sqrt{2vt_2}\Big\}.
\end{align}

\noindent \textbf{Step 3} (Exact root) Regarding $\sqrt{v}$ as the non-negative number $\mathbf{Y}$, we have
\begin{align*}
    \mathbf{E}_2=\Big\{\mathbf{Y}^2-b\sqrt{2t_2}\,\mathbf{Y}-[\M]_n<0\Big\}=:\{f(\mathbf{Y})<0\}.
\end{align*}
Since $[\M]_n\geq 0$, the quadratic function $f$ has a non-negative largest root
\begin{align*}
    y_{\text{root}}=\frac{b\sqrt{2t_2}+\sqrt{2b^2t_2+4[\M]_n}}{2}=b\sqrt{\frac{t_2}{2}}+\sqrt{[\M]_n+\frac{b^2t_2}{2}},
\end{align*}
and $f(y)<0$ with $y\geq 0$ implies $y<y_{\text{root}}$. Hence
\begin{align}
\label{event 3}
    \mathbf{E}_2\subset\mathbf{E}_3:=\Big\{\sqrt{v}<b\sqrt{\frac{t_2}{2}}+\sqrt{[\M]_n+\frac{b^2t_2}{2}}\Big\}.
\end{align}

\noindent \textbf{Step 4} (Proof of (i)) On $\mathbf{E}_1^{+}\cap\mathbf{E}_3$, we have
\begin{align*}
    \M_n<\sqrt{2t_1}\sqrt{v}+\frac{b}{3}t_1<\sqrt{2t_1}\Big(b\sqrt{\frac{t_2}{2}}+\sqrt{[\M]_n+\frac{b^2t_2}{2}}\Big)+\frac{b}{3}t_1=\sqrt{2t_1[\M]_n+b^2t_1t_2}+b\sqrt{t_1t_2}+\frac{bt_1}{3}.
\end{align*}
The same holds for $-\M_n$ on $\mathbf{E}_1^{-}\cap\mathbf{E}_3$. Together with \eqref{event 1}--\eqref{event 3}, the one-sided event has probability no less than $\P(\mathbf{E}_1^{+}\cap\mathbf{E}_2)\geq 1-e^{-t_1}-e^{-t_2}$, and the two-sided event in \eqref{eq th 1 1} has probability no less than $\P(\mathbf{E}_1^{+}\cap\mathbf{E}_1^{-}\cap\mathbf{E}_2)\geq 1-2e^{-t_1}-e^{-t_2}$. Finally, under \eqref{eq th 1 2},
\begin{align*}
    2e^{-t_1}+e^{-t_2}=\frac{(\sqrt{n}-1)\alpha}{\sqrt{n}}+\frac{\alpha}{\sqrt{n}}=\alpha.
\end{align*}

\noindent \textbf{Step 5} (Proof of (ii)) Denote by $t_{1,n}$ and $t_{2,n}$ the choice in \eqref{eq th 1 2}, so that $t_{1,n}\to\log(2/\alpha)$ and $t_{2,n}/n\to 0$. By \eqref{eq proof of th 1 0}, $v=n\sigma_n^2$ for every $n$. Since $(V_i,\mathcal{A}_i)_{i\geq 1}$ is a bounded MDS and $\E[V_i^2]\leq\E[Z_i^4]\leq b^2\E[Z_i^2]$, the orthogonality of martingale differences yields
\begin{align*}
    \E\Big[\big([\M]_n-n\sigma_n^2\big)^2\Big]=\E\Big[\Big(\sum_{i=1}^nV_i\Big)^2\Big]=\sum_{i=1}^n\E[V_i^2]\leq b^2n\sigma_n^2,
\end{align*}
that is, $[\M]_n/n-\sigma_n^2=O_{\P}(n^{-1/2})$, because $\sigma_n^2\leq b^2$. Now,
\begin{align*}
    \sqrt{n}L_n=\sqrt{\frac{2t_{1,n}[\M]_n}{n}+\frac{b^2t_{1,n}t_{2,n}}{n}}+b\sqrt{\frac{t_{1,n}t_{2,n}}{n}}+\frac{bt_{1,n}}{3\sqrt{n}}.
\end{align*}
Using $|\sqrt{x}-\sqrt{y}|\leq\sqrt{|x-y|}$ for $x,y\geq 0$, we obtain
\begin{align*}
    \Big|\sqrt{n}L_n-\sigma_n\sqrt{2\log(2/\alpha)}\Big|\leq& \sqrt{2t_{1,n}\Big|\frac{[\M]_n}{n}-\sigma_n^2\Big|+2\big|t_{1,n}-\log(2/\alpha)\big|\sigma_n^2+\frac{b^2t_{1,n}t_{2,n}}{n}}\\
    &+b\sqrt{\frac{t_{1,n}t_{2,n}}{n}}+\frac{bt_{1,n}}{3\sqrt{n}}\xrightarrow[]{\P}0.
\end{align*}
Since $\liminf_n\sigma_n^2>0$, dividing by $\sigma_n\sqrt{2\log(2/\alpha)}$ yields \eqref{eq th 1 3}. Thus, we finish the proof of Theorem \ref{th 1}.\\

\noindent \textbf{Proof of Proposition \ref{prop 1}} Denote $T=\sup_nt_{1,n}<\infty$, $K=\sqrt{2T}+1$ and
\begin{align*}
    d_n=2\sqrt{Tt_{2,n}}+\frac{T}{3},\qquad a_n=\min\big\{\max\{d_n^2,1\},n\big\}.
\end{align*}
Since $t_{2,n}/n\to 0$, we have $d_n^2/n\to 0$. Hence $1\leq a_n\leq n$ and $n/a_n\to\infty$. Let $(\varepsilon_i)_{i\geq 1}$ and $S_k$ be as in Lemma \ref{lemma A.4}, let $\mathcal{A}_0$ be the trivial $\sigma$-algebra and $\mathcal{A}_i=\sigma(\varepsilon_1,\dots,\varepsilon_i)$ for $i\geq 1$. For every $n\geq 1$, define
\begin{align*}
    \rho_n(k)=\sqrt{2t_{1,n}k+t_{1,n}t_{2,n}}+\sqrt{t_{1,n}t_{2,n}}+\frac{t_{1,n}}{3},\qquad k\geq 1,
\end{align*}
together with $\tau_n=\min\{k\in\{1,\dots,n\}:S_k>\rho_n(k)\}$, where $\min\emptyset=\infty$, and
\begin{align*}
    Z_i=b\,\varepsilon_i1\{\tau_n\geq i\},\qquad i\geq 1.
\end{align*}
Here $Z_i$ stands for $Z_{n,i}$, and we suppress the index $n$ to simplify notation.

\noindent \textbf{Step 1} (MDS structure) Since $\tau_n$ is a stopping time with respect to $(\mathcal{A}_i)_{i\geq 0}$, we have $\{\tau_n\geq i\}=\{\tau_n\leq i-1\}^c\in\mathcal{A}_{i-1}$. Since $\varepsilon_i$ is independent of $\mathcal{A}_{i-1}$ and has mean zero, $Z_i$ is measurable with respect to $\mathcal{A}_i$,
\begin{align*}
    \E[Z_i|\mathcal{A}_{i-1}]=b1\{\tau_n\geq i\}\E[\varepsilon_i]=0,\qquad |Z_i|\leq b,\qquad \E[Z_i^2|\mathcal{A}_{i-1}]=b^21\{\tau_n\geq i\}=Z_i^2.
\end{align*}
Hence $(Z_i,\mathcal{A}_i)_{i\geq 1}$ is an MDS bounded by $b$, $\<\M\>_n=[\M]_n=b^2(\tau_n\wedge n)$ and $\M_n=bS_{\tau_n\wedge n}$.\\

\noindent \textbf{Step 2} (Violation on $\{\tau_n\leq n\}$) On $\{\tau_n\leq n\}$, we have $[\M]_n=b^2\tau_n$ and
\begin{align*}
    \M_n=bS_{\tau_n}>b\rho_n(\tau_n)=\sqrt{2t_{1,n}[\M]_n+b^2t_{1,n}t_{2,n}}+b\sqrt{t_{1,n}t_{2,n}}+\frac{bt_{1,n}}{3}.
\end{align*}
Thus, the probability in \eqref{eq prop 1 1} is no less than $\P(\tau_n\leq n)=\P(\exists\,k\leq n: S_k>\rho_n(k))$.\\

\noindent \textbf{Step 3} (Reduction to Lemma \ref{lemma A.4}) Let $n$ be large enough such that $d_n^2\leq n$, so that $a_n\geq d_n^2$. For every $k\in[a_n,n]$, using $\sqrt{x+y}\leq\sqrt{x}+\sqrt{y}$ and $t_{1,n}\leq T$, we obtain
\begin{align*}
    \rho_n(k)\leq\sqrt{2Tk}+2\sqrt{Tt_{2,n}}+\frac{T}{3}=\sqrt{2Tk}+d_n\leq\sqrt{2Tk}+\sqrt{k}=K\sqrt{k}.
\end{align*}
Hence $\{\exists\,k\in[a_n,n]:S_k>K\sqrt{k}\}\subset\{\tau_n\leq n\}$, and Lemma \ref{lemma A.4} yields \eqref{eq prop 1 1}.\\

\noindent \textbf{Step 4} (MDS free of $n$) If $t_{1,n}=t_1$ and $t_{2,n}=t_2$ do not depend on $n$, neither does $\rho_n=:\rho$. Let $\tau=\inf\{k\geq 1:S_k>\rho(k)\}$ and $Z_i=b\,\varepsilon_i1\{\tau\geq i\}$. For every $n$, $1\{\tau\geq i\}=1\{\tau_n\geq i\}$ holds for all $i\leq n$, because $\tau_n=\tau$ on $\{\tau\leq n\}$, while on $\{\tau>n\}$ both indicators are equal to one for every $i\leq n$. Thus, $(Z_i)_{i\leq n}$ coincides with the construction above for every $n$, and the MDS does not depend on $n$. Thus, we finish the proof of Proposition \ref{prop 1}.\\

\noindent \textbf{Proof of Theorem \ref{th 2}} Denote $\mathbf{B}=\{v-\Gamma\leq\<\M\>_n\leq v\}$. The processes $(L_k(\lambda))_{0\leq k\leq n}$, $\lambda\in[0,3/b)$, and $(K_k(\lambda))_{0\leq k\leq n}$, $\lambda\geq 0$, defined in Steps 1 and 2 of the proof of Theorem \ref{th 1} are non-negative super-martingales. This part of that proof does not use that $\<\M\>_n$ is deterministic. Hence $\E[L_n(\lambda)]\leq 1$ and $\E[K_n(\lambda)]\leq 1$.\\

\noindent \textbf{Step 1} (Deviation of $\M_n$) Let $\lambda=\lambda(v,t_1)$ from Lemma \ref{lemma A.3} and
\begin{align*}
    \mathbf{F}_1^{\pm}=\Big\{\pm\M_n\geq\sqrt{2vt_1}+\frac{b}{3}t_1,\ \<\M\>_n\leq v\Big\}.
\end{align*}
On $\mathbf{F}_1^{+}$, \eqref{eq lemma A.3 1} yields
\begin{align*}
    \lambda\M_n-\psi(\lambda)\<\M\>_n\geq\lambda\Big(\sqrt{2vt_1}+\frac{b}{3}t_1\Big)-\psi(\lambda)v=t_1,
\end{align*}
so that $\mathbf{F}_1^{+}\subset\{L_n(\lambda)\geq e^{t_1}\}$ and Markov's inequality yields $\P(\mathbf{F}_1^{+})\leq e^{-t_1}$. Applying the same argument to the MDS $(-Z_i,\mathcal{A}_i)_{i\geq 1}$ yields $\P(\mathbf{F}_1^{-})\leq e^{-t_1}$.\\

\noindent \textbf{Step 2} (Lower tail of $[\M]_n$) Let $\lambda=\sqrt{2t_2/(b^2v)}$ and
\begin{align*}
    \mathbf{F}_2=\Big\{\<\M\>_n-[\M]_n\geq b\sqrt{2vt_2},\ \<\M\>_n\leq v\Big\}.
\end{align*}
On $\mathbf{F}_2$,
\begin{align*}
    \lambda\big(\<\M\>_n-[\M]_n\big)-\frac{\lambda^2b^2}{2}\<\M\>_n\geq\lambda b\sqrt{2vt_2}-\frac{\lambda^2b^2v}{2}=2t_2-t_2=t_2,
\end{align*}
so that $\mathbf{F}_2\subset\{K_n(\lambda)\geq e^{t_2}\}$ and $\P(\mathbf{F}_2)\leq e^{-t_2}$.\\

\noindent \textbf{Step 3} (Conclusion) On $\mathbf{B}\setminus\mathbf{F}_2$, we have $v-\Gamma\leq\<\M\>_n<[\M]_n+b\sqrt{2vt_2}$, that is, $\mathbf{Y}^2-b\sqrt{2t_2}\,\mathbf{Y}-([\M]_n+\Gamma)<0$ with $\mathbf{Y}=\sqrt{v}$. As in Step 3 of the proof of Theorem \ref{th 1}, this implies
\begin{align*}
    \sqrt{v}<b\sqrt{\frac{t_2}{2}}+\sqrt{[\M]_n+\Gamma+\frac{b^2t_2}{2}}.
\end{align*}
Since $\<\M\>_n\leq v$ on $\mathbf{B}$, we obtain on $\mathbf{B}\setminus(\mathbf{F}_1^{+}\cup\mathbf{F}_2)$ that
\begin{align*}
    \M_n<\sqrt{2t_1}\sqrt{v}+\frac{b}{3}t_1<\sqrt{2t_1\big([\M]_n+\Gamma\big)+b^2t_1t_2}+b\sqrt{t_1t_2}+\frac{bt_1}{3},
\end{align*}
and the same bound holds for $-\M_n$ on $\mathbf{B}\setminus(\mathbf{F}_1^{-}\cup\mathbf{F}_2)$. Therefore, the one-sided event contains $\mathbf{B}\setminus(\mathbf{F}_1^{+}\cup\mathbf{F}_2)$ and the two-sided event in \eqref{eq th 2 1} contains $\mathbf{B}\setminus(\mathbf{F}_1^{+}\cup\mathbf{F}_1^{-}\cup\mathbf{F}_2)$. Their probabilities are no less than $1-\P(\mathbf{B}^c)-e^{-t_1}-e^{-t_2}$ and $1-\P(\mathbf{B}^c)-2e^{-t_1}-e^{-t_2}$, respectively. Thus, we finish the proof of Theorem \ref{th 2}.\\

\noindent \textbf{Proof of Theorem \ref{th 3}} Denote $v_0=n\sigma_n^2$ and $V_k=\frac{1}{3}([\M]_k+2\<\M\>_k)$ for $0\leq k\leq n$.\\

\noindent \textbf{Step 1} (Self-normalized super-martingale) For $\lambda\in\mathbb{R}$ and $0\leq k\leq n$, define
\begin{align*}
    R_k(\lambda)=\exp\Big(\lambda\M_k-\frac{\lambda^2}{2}V_k\Big)=\exp\Big(\sum_{i=1}^k\Big(\lambda Z_i-\frac{\lambda^2}{6}\big(Z_i^2+2\E[Z_i^2|\mathcal{A}_{i-1}]\big)\Big)\Big).
\end{align*}
By Lemma \ref{lemma A.5}, $\E[R_k(\lambda)|\mathcal{A}_{k-1}]\leq R_{k-1}(\lambda)$ almost surely. Thus, $(R_k(\lambda),\mathcal{A}_k)_{0\leq k\leq n}$ is a non-negative super-martingale and $\E[R_n(\lambda)]\leq 1$.\\

\noindent \textbf{Step 2} (Deviation of $\M_n$ at a deterministic level) Let $w=v_0+x$, $\lambda=\sqrt{2t/w}$ and
\begin{align*}
    \mathbf{F}^{\pm}=\Big\{\pm\M_n\geq\sqrt{2wt},\ V_n\leq w\Big\}.
\end{align*}
On $\mathbf{F}^{+}$,
\begin{align*}
    \lambda\M_n-\frac{\lambda^2}{2}V_n\geq\lambda\sqrt{2wt}-\frac{\lambda^2w}{2}=2t-t=t,
\end{align*}
so that $\mathbf{F}^{+}\subset\{R_n(\lambda)\geq e^{t}\}$ and Markov's inequality yields $\P(\mathbf{F}^{+})\leq e^{-t}$. Applying the same argument to the MDS $(-Z_i,\mathcal{A}_i)_{i\geq 1}$, which has the same $V_n$, yields $\P(\mathbf{F}^{-})\leq e^{-t}$.\\

\noindent \textbf{Step 3} (Moment-level coupling) Denote
\begin{align*}
    \mathbf{G}=\Big\{\big|\<\M\>_n-v_0\big|\leq x,\ \big|[\M]_n-v_0\big|\leq x\Big\}.
\end{align*}
On $\mathbf{G}$, we have $V_n\leq v_0+x=w$ and $w\leq[\M]_n+2x$. Hence, on $\mathbf{G}\setminus\mathbf{F}^{\pm}$,
\begin{align*}
    \pm\M_n<\sqrt{2wt}\leq\sqrt{2t\big([\M]_n+2x\big)}.
\end{align*}
Therefore, the one-sided event contains $\mathbf{G}\setminus\mathbf{F}^{+}$ and the two-sided event in \eqref{eq th 3 1} contains $\mathbf{G}\setminus(\mathbf{F}^{+}\cup\mathbf{F}^{-})$, whose probabilities are no less than $1-\P(\mathbf{G}^c)-e^{-t}$ and $1-\P(\mathbf{G}^c)-2e^{-t}$, respectively.\\

\noindent \textbf{Step 4} (Bounding $\P(\mathbf{G}^c)$) Since $[\M]_n-\<\M\>_n=\sum_{i=1}^n(Z_i^2-\E[Z_i^2|\mathcal{A}_{i-1}])$ is a sum of martingale differences, their orthogonality yields
\begin{align*}
    \E\big[([\M]_n-\<\M\>_n)^2\big]=\sum_{i=1}^n\E\Big[\big(Z_i^2-\E[Z_i^2|\mathcal{A}_{i-1}]\big)^2\Big]\leq\sum_{i=1}^n\E[Z_i^4].
\end{align*}
Since $|[\M]_n-v_0|\leq|[\M]_n-\<\M\>_n|+|\<\M\>_n-v_0|$, we have $\mathbf{G}^c\subset\{|[\M]_n-\<\M\>_n|>x/2\}\cup\{|\<\M\>_n-v_0|>x/2\}$. Thus, by Markov's and Chebyshev's inequalities,
\begin{align*}
    \P(\mathbf{G}^c)&\leq\P\Big(|[\M]_n-\<\M\>_n|>\frac{x}{2}\Big)+\P\Big(|\<\M\>_n-v_0|>\frac{x}{2}\Big)\\
    &\leq\frac{2}{x}\E\big|\<\M\>_n-v_0\big|+\frac{4}{x^2}\sum_{i=1}^n\E[Z_i^4].
\end{align*}
Together with Step 3, this proves \eqref{eq th 3 1}. Thus, we finish the proof of Theorem \ref{th 3}.\\

\noindent \textbf{Proof of Corollary \ref{corollary 1}}\\
\noindent \textbf{Step 1} (Reduction to Theorem \ref{th 1}) The sequence $(Z_i-\mu,\mathcal{A}_i)_{i\geq 1}$ is an MDS with $|Z_i-\mu|\leq b$, whose predictable quadratic variation $\sum_{i=1}^n\E[(Z_i-\mu)^2|\mathcal{A}_{i-1}]$ is deterministic. For this MDS, $\M_n=n(\overline{\mathbf{M}}_n-\mu)$ and $[\M]_n=\sum_{i=1}^n(Z_i-\mu)^2$. Denote
\begin{align*}
    Q_1=2t_1[\mathbb{V}]_n+b^2t_1t_2,\qquad Q_2=b\sqrt{t_1t_2}+\frac{bt_1}{3},
\end{align*}
and
\begin{align*}
    \mathbf{E}^{\pm}=\Big\{\pm n(\overline{\mathbf{M}}_n-\mu)\leq\sqrt{2t_1\sum_{i=1}^n(Z_i-\mu)^2+b^2t_1t_2}+Q_2\Big\}.
\end{align*}
By Theorem \ref{th 1} (i), $\P(\mathbf{E}^{\pm})\geq 1-e^{-t_1}-e^{-t_2}$ and $\P(\mathbf{E}^{+}\cap\mathbf{E}^{-})\geq 1-2e^{-t_1}-e^{-t_2}$.\\

\noindent \textbf{Step 2} (Exact solution) Denote by $L$ the radius in \eqref{eq corollary 1 1} and let $x=\pm(\overline{\mathbf{M}}_n-\mu)$. It suffices to show that $x\leq L$ on $\mathbf{E}^{\pm}$. Since $L\geq Q_2/(n-2t_1)\geq Q_2/n$, this is obvious when $nx\leq Q_2$. Otherwise, the identity $\sum_{i=1}^n(Z_i-\mu)^2=[\mathbb{V}]_n+n(\overline{\mathbf{M}}_n-\mu)^2$ yields, on $\mathbf{E}^{\pm}$,
\begin{align*}
    0<nx-Q_2\leq\sqrt{Q_1+2t_1nx^2}.
\end{align*}
Squaring both sides, we obtain
\begin{align*}
    n(n-2t_1)x^2-2nQ_2x+Q_2^2-Q_1\leq 0.
\end{align*}
Since $n(n-2t_1)>0$, $x$ is no larger than the largest root of the left-hand side, that is,
\begin{align*}
    x&\leq\frac{nQ_2+\sqrt{2t_1nQ_2^2+n(n-2t_1)Q_1}}{n(n-2t_1)}=\frac{1}{n-2t_1}\Big(Q_2+\sqrt{\frac{2t_1Q_2^2}{n}+\frac{(n-2t_1)Q_1}{n}}\Big)\\
    &\leq\frac{1}{n-2t_1}\Big(\sqrt{Q_1}+\Big(1+\sqrt{\frac{2t_1}{n}}\Big)Q_2\Big)=L,
\end{align*}
where the last inequality uses $\sqrt{x+y}\leq\sqrt{x}+\sqrt{y}$ and $n-2t_1\leq n$. Together with Step 1, this proves \eqref{eq corollary 1 1} and its one-sided versions. Under \eqref{eq th 1 2}, $2e^{-t_1}+e^{-t_2}=\alpha$ as in Step 4 of the proof of Theorem \ref{th 1}.\\

\noindent \textbf{Step 3} (Proof of (ii)) Denote by $t_{1,n}$ and $t_{2,n}$ the choice in \eqref{eq th 1 2}. Since $\E[Z_i]=\mu$, we have $\text{Var}(Z_i)=\E[(Z_i-\mu)^2]$. Step 5 of the proof of Theorem \ref{th 1}, applied to the MDS $(Z_i-\mu,\mathcal{A}_i)_{i\geq 1}$, yields $\frac{1}{n}\sum_{i=1}^n(Z_i-\mu)^2-\sigma_n^2=O_{\P}(n^{-1/2})$. Moreover, the orthogonality of martingale differences yields $\E[(\overline{\mathbf{M}}_n-\mu)^2]=\sigma_n^2/n\leq b^2/n$. Hence, by the identity in Step 2, $[\mathbb{V}]_n/n-\sigma_n^2=O_{\P}(n^{-1/2})$. Now,
\begin{align*}
    \sqrt{n}L_n=\frac{n}{n-2t_{1,n}}\Big(\sqrt{\frac{2t_{1,n}[\mathbb{V}]_n}{n}+\frac{b^2t_{1,n}t_{2,n}}{n}}+\Big(1+\sqrt{\frac{2t_{1,n}}{n}}\Big)\Big(b\sqrt{\frac{t_{1,n}t_{2,n}}{n}}+\frac{bt_{1,n}}{3\sqrt{n}}\Big)\Big).
\end{align*}
Since $n/(n-2t_{1,n})\to 1$, $t_{1,n}\to\log(2/\alpha)$ and $t_{2,n}/n\to 0$, the same argument as in Step 5 of the proof of Theorem \ref{th 1} yields $\sqrt{n}L_n-\sigma_n\sqrt{2\log(2/\alpha)}\xrightarrow[]{\P}0$. Since $\liminf_n\sigma_n^2>0$, this proves (ii). Thus, we finish the proof of Corollary \ref{corollary 1}.\\

\noindent \textbf{Proof of Theorem \ref{th 4}}\\
\noindent \textbf{Step 0} (Auxiliary triangular array) The events in (i), the radius $L_n$ and the ratio in (ii) are measurable functions of $(Z_1,\dots,Z_n)$, and $\sigma_n^2$ depends only on the distribution of $(Z_1,\dots,Z_n)$. By \ding{172}, it therefore suffices to prove Theorem \ref{th 4} with $(Z_1,\dots,Z_n)$ replaced by $(Z_{n1},\dots,Z_{nn})$ from the auxiliary triangular array introduced before Lemma \ref{lemma A.7}. We use the notation \eqref{eq notation block}. In particular, $\text{Var}(S_n)=n\sigma_n^2$.\\

\noindent \textbf{Step 1} (Reduction to Theorem \ref{th 2}) By Lemma \ref{lemma A.7} (ii) and (iii), $(M_{nj},\mathcal{F}_{n,jl})_{1\leq j\leq m}$ is a martingale difference sequence bounded by $\bar{b}_n$, and $S_n=\sum_{j=1}^mM_{nj}=n(\P_n(h_n)-\mu_n)$. Setting $M_{nj}=0$ and using the $\sigma$-algebra $\mathcal{F}_{nn}$ for $j>m$, we define $\<\M\>_m$ and $[\M]_m$ for this sequence as in \eqref{eq notation MDS}, with $n$ replaced by $m$. Let $v=\text{Var}(S_n)+\Gamma_n/2$. By Lemma \ref{lemma A.8}, $\P(\<\M\>_m\notin[v-\Gamma_n,v])=0$. If $v=0$, then $S_n=0$ almost surely and the claims of (i) hold trivially, because $L_n\geq 0$. Otherwise, Theorem \ref{th 2}, applied with $b=\bar{b}_n$, $\Gamma=\Gamma_n$ and $n$ replaced by $m$ to the martingale difference sequences $(M_{nj})_{j\geq 1}$ and $(-M_{nj})_{j\geq 1}$, shows that the events
\begin{align*}
    \mathbf{E}^{\pm}=\Big\{\pm S_n\leq\sqrt{2t_1\big([\M]_m+\Gamma_n\big)+\bar{b}_n^2t_1t_2}+\bar{b}_n\sqrt{t_1t_2}+\frac{\bar{b}_nt_1}{3}\Big\}
\end{align*}
satisfy $\P(\mathbf{E}^{\pm})\geq 1-e^{-t_1}-e^{-t_2}$ and $\P(\mathbf{E}^{+}\cap\mathbf{E}^{-})\geq 1-2e^{-t_1}-e^{-t_2}$.\\

\noindent \textbf{Step 2} (Boundary terms) By Lemma \ref{lemma A.7} (ii) and (iii), $|M_{nj}-T_{nj}|\leq 2\|h_n\|_{TV}\tilde{\Phi}_n$. Hence, by the triangle inequality in $\mathbb{R}^m$,
\begin{align*}
    \sqrt{[\M]_m}\leq\sqrt{\sum_{j=1}^mT_{nj}^2}+\delta_n,\qquad \delta_n=2\sqrt{m}\,\|h_n\|_{TV}\tilde{\Phi}_n.
\end{align*}
For $t>0$ and $a,\delta,\Gamma,Q\geq 0$, squaring both sides and using $\sqrt{2t(a^2+\Gamma)+Q}\geq\sqrt{2t}\,a$ shows that $\sqrt{2t((a+\delta)^2+\Gamma)+Q}\leq\sqrt{2t(a^2+\Gamma)+Q}+\sqrt{2t}\,\delta$. Applying this with $a^2=\sum_{j=1}^mT_{nj}^2$, $\delta=\delta_n$, $\Gamma=\Gamma_n$ and $Q=\bar{b}_n^2t_1t_2$, we obtain, on $\mathbf{E}^{\pm}$,
\begin{align*}
    \pm S_n\leq\sqrt{2t_1\Big(\sum_{j=1}^mT_{nj}^2+\Gamma_n\Big)+\bar{b}_n^2t_1t_2}+Q_2,\qquad Q_2=\bar{b}_n\sqrt{t_1t_2}+\frac{\bar{b}_nt_1}{3}+2\sqrt{2t_1m}\,\|h_n\|_{TV}\tilde{\Phi}_n.
\end{align*}

\noindent \textbf{Step 3} (Exact solution) Denote $Q_1=2t_1(n\hat{V}_n+\Gamma_n)+\bar{b}_n^2t_1t_2$ and $x=\pm(\P_n(h_n)-\mu_n)$. It suffices to show that $x\leq L_n$ on $\mathbf{E}^{\pm}$. Since $L_n\geq Q_2/(n-2t_1l)\geq Q_2/n$, this is obvious when $nx\leq Q_2$. Otherwise, Step 2 and Lemma \ref{lemma A.7} (iv) yield, on $\mathbf{E}^{\pm}$,
\begin{align*}
    0<nx-Q_2\leq\sqrt{Q_1+2t_1nlx^2}.
\end{align*}
Squaring both sides, we obtain
\begin{align*}
    n(n-2t_1l)x^2-2nQ_2x+Q_2^2-Q_1\leq 0.
\end{align*}
Since $n-2t_1l=l(m-2t_1)>0$, $x$ is no larger than the largest root of the left-hand side, that is,
\begin{align*}
    x&\leq\frac{nQ_2+\sqrt{2t_1lnQ_2^2+n(n-2t_1l)Q_1}}{n(n-2t_1l)}=\frac{1}{n-2t_1l}\Big(Q_2+\sqrt{\frac{2t_1lQ_2^2}{n}+\frac{(n-2t_1l)Q_1}{n}}\Big)\\
    &\leq\frac{1}{n-2t_1l}\Big(\sqrt{Q_1}+\Big(1+\sqrt{\frac{2t_1}{m}}\Big)Q_2\Big)=L_n,
\end{align*}
where the last inequality uses $\sqrt{x+y}\leq\sqrt{x}+\sqrt{y}$, $n-2t_1l\leq n$ and $l/n=1/m$. Together with Step 1, this proves \eqref{eq th 4 1} and its one-sided versions. Under \eqref{eq th 1 2} with $n$ replaced by $m$, $2e^{-t_1}+e^{-t_2}=\alpha$ as in Step 4 of the proof of Theorem \ref{th 1}. This proves (i).\\

\noindent \textbf{Step 4} (Preparation for (ii)) Denote $C=\sup_n\max\{b_n-a_n,\|h_n\|_{TV}\}$ and $\Psi=\sum_{r=1}^\infty r\,\tilde{\phi}_2(r)$. Since $\tilde{\Phi}_n\leq\Psi_n\leq\Psi$, we have $\bar{b}_n\leq C(l_n+2\Psi)$ and $\Gamma_n/(2m)\leq 14C^2\Psi+2C^2\Psi^2$. As in Step 2 (c) of the proof of Lemma \ref{lemma A.8}, $|\E[Y_{ni}Y_{nk}]|\leq C^2\tilde{\phi}_1(k-i)$ for $i<k$, so that $\sigma_n^2=\frac{1}{n}\text{Var}(S_n)\leq C^2(1+2\Psi)$. Moreover, $m=n/l_n\to\infty$, because $l_n^2=o(n)$. Denote by $t_{1,n}$ and $t_{2,n}$ the choice \eqref{eq th 1 2} with $n$ replaced by $m$, so that $t_{1,n}\to\log(2/\alpha)$, $2t_{1,n}<m$ for $n$ large enough, and $t_{2,n}\leq\log(\sqrt{n}/\alpha)$.\\

\noindent \textbf{Step 5} (Consistency of $\hat{V}_n$) Since $(M_{nj}^2-\E[M_{nj}^2|\mathcal{F}_{n,(j-1)l}])_{1\leq j\leq m}$ is a sequence of martingale differences and $|M_{nj}|\leq\bar{b}_n$, their orthogonality yields
\begin{align*}
    \E\big[([\M]_m-\<\M\>_m)^2\big]\leq\sum_{j=1}^m\E\big[M_{nj}^4\big]\leq\bar{b}_n^2\sum_{j=1}^m\E\big[M_{nj}^2\big]=\bar{b}_n^2n\sigma_n^2.
\end{align*}
Thus, $([\M]_m-\<\M\>_m)/n=O_{\P}(l_n/\sqrt{n})$. By Lemma \ref{lemma A.8}, $|\<\M\>_m/n-\sigma_n^2|\leq \Gamma_n/(2n)=O(1/l_n)$. Hence $[\M]_m/n-\sigma_n^2\xrightarrow[]{\P}0$. By the triangle inequality in $\mathbb{R}^m$ as in Step 2,
\begin{align*}
    \Big|\sqrt{\frac{1}{n}\sum_{j=1}^mT_{nj}^2}-\sqrt{\frac{[\M]_m}{n}}\Big|\leq\frac{\delta_n}{\sqrt{n}}=\frac{2\|h_n\|_{TV}\tilde{\Phi}_n}{\sqrt{l_n}}\to 0.
\end{align*}
Since $\sigma_n^2$ is bounded, $\frac{1}{n}\sum_{j=1}^mT_{nj}^2-\sigma_n^2\xrightarrow[]{\P}0$. Finally, $\E[l_n(\P_n(h_n)-\mu_n)^2]=l_n\sigma_n^2/n\to 0$, and Lemma \ref{lemma A.7} (iv) yields $\hat{V}_n-\sigma_n^2\xrightarrow[]{\P}0$.\\

\noindent \textbf{Step 6} (Proof of (ii)) By \eqref{eq th 4 2},
\begin{align*}
    \sqrt{n}L_n=\frac{m}{m-2t_{1,n}}\bigg(&\sqrt{2t_{1,n}\hat{V}_n+\frac{2t_{1,n}\Gamma_n}{n}+\frac{\bar{b}_n^2t_{1,n}t_{2,n}}{n}}\\
    &+\Big(1+\sqrt{\frac{2t_{1,n}}{m}}\Big)\Big(\bar{b}_n\sqrt{\frac{t_{1,n}t_{2,n}}{n}}+\frac{\bar{b}_nt_{1,n}}{3\sqrt{n}}+2\sqrt{\frac{2t_{1,n}}{l_n}}\,\|h_n\|_{TV}\tilde{\Phi}_n\Big)\bigg).
\end{align*}
By Step 4, $m/(m-2t_{1,n})\to 1$, $\Gamma_n/n=O(1/l_n)\to 0$, $\bar{b}_n^2t_{2,n}/n=O(l_n^2\log n/n)\to 0$ and $\tilde{\Phi}_n/\sqrt{l_n}\to 0$. Together with Step 5, $t_{1,n}\to\log(2/\alpha)$ and $|\sqrt{x}-\sqrt{y}|\leq\sqrt{|x-y|}$, the same argument as in Step 5 of the proof of Theorem \ref{th 1} yields $\sqrt{n}L_n-\sigma_n\sqrt{2\log(2/\alpha)}\xrightarrow[]{\P}0$. Since $\liminf_n\sigma_n^2>0$, this proves (ii). Thus, we finish the proof of Theorem \ref{th 4}.\\

\noindent \textbf{Proof of Theorem \ref{th 5}}\\
\noindent \textbf{Step 0} (Auxiliary triangular array) As in Step 0 of the proof of Theorem \ref{th 4}, it suffices to prove Theorem \ref{th 5} with $(Z_1,\dots,Z_n)$ replaced by $(Z_{n1},\dots,Z_{nn})$, and we use the notation \eqref{eq notation block}. In particular, $\text{Var}(S_n)=n\sigma_n^2$.\\

\noindent \textbf{Step 1} (Reduction to Theorem \ref{th 3}) By Lemma \ref{lemma A.7} (ii) and Lemma \ref{lemma A.9} (iv), $(M_{nj},\mathcal{F}_{n,jl})_{1\leq j\leq m}$ is a martingale difference sequence with finite fourth moments and $S_n=\sum_{j=1}^mM_{nj}=n(\P_n(h_n)-\mu_n)$. Setting $M_{nj}=0$ and using the $\sigma$-algebra $\mathcal{F}_{nn}$ for $j>m$, we define $\<\M\>_m$ and $[\M]_m$ as in \eqref{eq notation MDS}, with $n$ replaced by $m$, so that $\E[\<\M\>_m]=\text{Var}(S_n)$. Theorem \ref{th 3}, applied with $n$ replaced by $m$ to $(M_{nj})_{j\geq 1}$ and $(-M_{nj})_{j\geq 1}$, together with Lemma \ref{lemma A.9} (ii) and (iv), shows that the events $\mathbf{E}^{\pm}=\{\pm S_n\leq\sqrt{2t([\M]_m+2x)}\}$ satisfy
\begin{align*}
    \P(\mathbf{E}^{\pm})\geq 1-e^{-t}-\frac{104\,m\,c_n\Psi_n^\tau}{x}-\frac{32\big(nl^3\kappa_n^4+16\,m\,\kappa_n^{4-4\gamma}\|h_n\|_{\text{Lip}}^{4\gamma}(\Phi_n^\tau)^4\big)}{x^2},
\end{align*}
and $\P(\mathbf{E}^{+}\cap\mathbf{E}^{-})$ satisfies the same bound with $e^{-t}$ replaced by $2e^{-t}$.\\

\noindent \textbf{Step 2} (Boundary terms) Let $\mathbf{G}=\{\sum_{j=1}^m(\varepsilon_{n,jl}-\varepsilon_{n,(j-1)l})^2\leq x\}$. By Markov's inequality and Lemma \ref{lemma A.9} (iii), $\P(\mathbf{G}^c)\leq 8mc_n\Psi_n^\tau/x$. On $\mathbf{G}$, Lemma \ref{lemma A.7} (ii) and the triangle inequality in $\mathbb{R}^m$ yield $\sqrt{[\M]_m}\leq\sqrt{\sum_{j=1}^mT_{nj}^2}+\sqrt{x}$. The inequality in Step 2 of the proof of Theorem \ref{th 4}, applied with $a^2=\sum_{j=1}^mT_{nj}^2$, $\delta=\sqrt{x}$, $\Gamma=2x$ and $Q=0$, shows that, on $\mathbf{E}^{\pm}\cap\mathbf{G}$,
\begin{align*}
    \pm S_n\leq\sqrt{2t\Big(\sum_{j=1}^mT_{nj}^2+2x\Big)}+\sqrt{2tx}.
\end{align*}
By Step 1, $\P(\mathbf{E}^{\pm}\cap\mathbf{G})\geq 1-e^{-t}-\eta_n(x)$ and $\P(\mathbf{E}^{+}\cap\mathbf{E}^{-}\cap\mathbf{G})\geq 1-2e^{-t}-\eta_n(x)$.\\

\noindent \textbf{Step 3} (Exact solution) Denote $y=\pm(\P_n(h_n)-\mu_n)$. By Lemma \ref{lemma A.7} (iv), on $\mathbf{E}^{\pm}\cap\mathbf{G}$,
\begin{align*}
    ny\leq\sqrt{Q_1+2tnly^2}+Q_2,\qquad Q_1=2t\big(n\hat{V}_n+2x\big),\qquad Q_2=\sqrt{2tx}.
\end{align*}
This is the inequality of Step 3 of the proof of Theorem \ref{th 4} with $t_1=t$, and the same argument yields $y\leq L_n$. Together with Step 2, this proves (i). For $t=\log(2/\alpha)$, $2e^{-t}=\alpha$.\\

\noindent \textbf{Step 4} (Proof of (ii)) Under the assumptions of (ii), $c_n$, $\Psi_n^\tau$ and $\Phi_n^\tau$ are bounded, so that $\eta_n(x_n)=O(n/(l_nx_n))+O(nl_n^3/x_n^2)+O(n/(l_nx_n^2))\to 0$. As in Step 3 (c) of the proof of Lemma \ref{lemma A.9}, $|\E[Y_{ni}Y_{nk}]|\leq c_n\tau_1(k-i)^\beta$ for $i<k$, so that $\sigma_n^2\leq\kappa_n^2+2c_n\Psi_n^\tau$ is bounded. By Lemma \ref{lemma A.9} (iv) and the orthogonality of martingale differences, $\E[([\M]_m-\<\M\>_m)^2]\leq\sum_{j=1}^m\E[M_{nj}^4]=O(nl_n^3)$, so that $([\M]_m-\<\M\>_m)/n=O_{\P}(\sqrt{l_n^3/n})\to 0$. By Lemma \ref{lemma A.9} (ii), $\E|\<\M\>_m/n-\sigma_n^2|\leq 52c_n\Psi_n^\tau/l_n\to 0$. Hence $[\M]_m/n-\sigma_n^2\xrightarrow[]{\P}0$. By the triangle inequality in $\mathbb{R}^m$ and Lemma \ref{lemma A.9} (iii),
\begin{align*}
    \E\Big|\sqrt{\frac{1}{n}\sum_{j=1}^mT_{nj}^2}-\sqrt{\frac{[\M]_m}{n}}\Big|^2\leq\frac{1}{n}\sum_{j=1}^m\E\big[(\varepsilon_{n,jl}-\varepsilon_{n,(j-1)l})^2\big]\leq\frac{8c_n\Psi_n^\tau}{l_n}\to 0.
\end{align*}
Since $\sigma_n^2$ is bounded, $\frac{1}{n}\sum_{j=1}^mT_{nj}^2-\sigma_n^2\xrightarrow[]{\P}0$. Finally, $\E[l_n(\P_n(h_n)-\mu_n)^2]=l_n\sigma_n^2/n\to 0$, and Lemma \ref{lemma A.7} (iv) yields $\hat{V}_n-\sigma_n^2\xrightarrow[]{\P}0$. Since
\begin{align*}
    \sqrt{n}L_n=\frac{m}{m-2t}\Big(\sqrt{2t\Big(\hat{V}_n+\frac{2x_n}{n}\Big)}+\Big(1+\sqrt{\frac{2t}{m}}\Big)\sqrt{\frac{2tx_n}{n}}\Big)
\end{align*}
with $t=\log(2/\alpha)$, $m\to\infty$ and $x_n/n\to 0$, we obtain $\sqrt{n}L_n-\sigma_n\sqrt{2\log(2/\alpha)}\xrightarrow[]{\P}0$, and $\liminf_n\sigma_n^2>0$ proves (ii). Thus, we finish the proof of Theorem \ref{th 5}.\\

\noindent \textbf{Proof of Theorem \ref{th 6}}\\
\noindent \textbf{Step 0} (Notation) The events in (i), the radius $L_n$ and the ratio in (ii) are measurable functions of $(Z_1,\dots,Z_n)$, and $\sigma_n^2$ depends only on the distribution of $(Z_1,\dots,Z_n)$. Replacing $(\Omega,\mathcal{F},\P)$ by its product with $[0,1]^{2m}$, equipped with the Lebesgue measure, changes neither these quantities nor the coefficients of $\mathbf{Z}$. We may therefore assume that there exist random variables $\vartheta_1,\dots,\vartheta_m,\vartheta_1',\dots,\vartheta_m'$, uniformly distributed on $[0,1]$, such that $\mathcal{F}_{-\infty}^{\infty},\sigma(\vartheta_1),\dots,\sigma(\vartheta_m),\sigma(\vartheta_1'),\dots,\sigma(\vartheta_m')$ are independent. The big and the small blocks are $I_j=\{(j-1)l+1,\dots,jl-d\}$ and $J_j=\{jl-d+1,\dots,jl\}$, $1\leq j\leq m$, and we denote $Y_i=h_n(Z_i)-\mu_n$. Since $h_n$ and $\mu_n$ take values in $[a_n,b_n]$, we have $|Y_i|\leq b_n-a_n$, $\E[U_j]=(l-d)\mu_n$, $\E[W_j]=d\mu_n$, $|U_j-(l-d)\mu_n|\leq b_U$ and $|W_j-d\mu_n|\leq b_W$. Moreover, $n(\P_n(h_n)-\mu_n)=\sum_{j=1}^m(U_j-(l-d)\mu_n)+\sum_{j=1}^m(W_j-d\mu_n)$. For $v\in\mathbb{R}^m$, denote by $\Pi v$ the vector with entries $v_j-\frac{1}{m}\sum_{k=1}^mv_k$ and by $\|v\|$ the Euclidean norm. Then $\Pi$ is the orthogonal projection onto the orthogonal complement of $(1,\dots,1)$, and $[\mathbb{V}]_m^U=\|\Pi U\|^2$ and $[\mathbb{V}]_m^W=\|\Pi W\|^2$, where $U=(U_1,\dots,U_m)$ and $W=(W_1,\dots,W_m)$.\\

\noindent \textbf{Step 1} (Coupling) Denote by $e_j=jl-d$ the last index of $I_j$. Let $U_1^*=U_1$ and, for $2\leq j\leq m$, let $\mathcal{M}_j=\mathcal{F}_{-\infty}^{e_{j-1}}\vee\sigma(\vartheta_k:2\leq k<j)$. Since $\vartheta_j$ is independent of $\sigma(U_j)\vee\mathcal{M}_j$, Lemma 5 of \citeA{dedecker2004coupling} yields a random variable $U_j^*$, measurable with respect to $\mathcal{M}_j\vee\sigma(U_j)\vee\sigma(\vartheta_j)$, independent of $\mathcal{M}_j$ and distributed as $U_j$, such that $\E|U_j-U_j^*|=\tau(\mathcal{M}_j,U_j)$.
\begin{enumerate}[(a)]
    \item Since $U_j$ is measurable with respect to $\mathcal{F}_{-\infty}^{e_j}$, $U_j^*$ is measurable with respect to $\mathcal{M}_{j+1}$ for $1\leq j<m$. Hence $U_1^*,\dots,U_{j-1}^*$ are measurable with respect to $\mathcal{M}_j$. Since $U_j^*$ is independent of $\mathcal{M}_j$, it is independent of $(U_1^*,\dots,U_{j-1}^*)$ for every $2\leq j\leq m$. By induction, $U_1^*,\dots,U_m^*$ are independent.
    \item Let $X=(Z_i)_{i\in I_j}$. Since $\sigma(\vartheta_k:2\leq k<j)$ is independent of $\mathcal{F}_{-\infty}^\infty$, a regular conditional distribution of $X$ given $\mathcal{F}_{-\infty}^{e_{j-1}}$ is also a regular conditional distribution of $X$ given $\mathcal{M}_j$. Moreover, $U_j=H(X)$ with $H(z)=\sum_{k=1}^{l-d}h_n(z_k)$, so that $f\circ H$ is $\|h_n\|_{\text{Lip}}$-Lipschitz with respect to $\|\cdot\|_1$ on $\mathbb{R}^{l-d}$ for every $f\in\Lambda_1(\mathbb{R})$. Since the first index of $I_j$ is $e_{j-1}+d+1$, Definition \ref{def tau k} yields
    \begin{align*}
        \E|U_j-U_j^*|=\tau\big(\mathcal{F}_{-\infty}^{e_{j-1}},H(X)\big)\leq\|h_n\|_{\text{Lip}}\,\tau\big(\mathcal{F}_{-\infty}^{e_{j-1}},X\big)\leq\|h_n\|_{\text{Lip}}(l-d)\,\tau_\infty(d+1).
    \end{align*}
\end{enumerate}
In the same way, with $jl$ the last index of $J_j$, $W_1^*=W_1$, $\mathcal{M}_j'=\mathcal{F}_{-\infty}^{(j-1)l}\vee\sigma(\vartheta_k':2\leq k<j)$ and $\vartheta_j'$ in place of $\vartheta_j$, we obtain independent random variables $W_1^*,\dots,W_m^*$ such that $W_j^*$ is distributed as $W_j$ and, since the first index of $J_j$ is $(j-1)l+(l-d)+1$, $\E|W_j-W_j^*|\leq\|h_n\|_{\text{Lip}}\,d\,\tau_\infty(l-d+1)$. Since $\tau_\infty$ is non-increasing and $l-d\geq d$, we obtain
\begin{align*}
    \E\Big[\sum_{j=1}^m\big(|U_j-U_j^*|+|W_j-W_j^*|\big)\Big]\leq(m-1)\,l\,\|h_n\|_{\text{Lip}}\,\tau_\infty(d)\leq\|h_n\|_{\text{Lip}}\,n\,\tau_\infty(d).
\end{align*}

\noindent \textbf{Step 2} (Five events) Let $\mathbf{G}=\{\sum_{j=1}^m(|U_j-U_j^*|+|W_j-W_j^*|)\leq x\}$. By Markov's inequality and Step 1, $\P(\mathbf{G}^c)\leq\|h_n\|_{\text{Lip}}\,n\,\tau_\infty(d)/x$. Let $\mathcal{A}_0$ be the trivial $\sigma$-algebra and $\mathcal{A}_j=\sigma(U_1^*,\dots,U_j^*)$ for $1\leq j\leq m$. By Step 1 (a), $\E[U_j^*|\mathcal{A}_{j-1}]=\E[U_j]=(l-d)\mu_n$ and $\sum_{j=1}^m\E[(U_j^*-(l-d)\mu_n)^2|\mathcal{A}_{j-1}]=\sum_{j=1}^m\text{Var}(U_j)$ is deterministic, and $|U_j^*-(l-d)\mu_n|\leq b_U$ almost surely, since $U_j^*$ is distributed as $U_j$. Denote $U^*=(U_1^*,\dots,U_m^*)$ and $[\mathbb{V}]^{U*}=\|\Pi U^*\|^2$. Corollary \ref{corollary 1} (i), applied with $n$ replaced by $m$, $\mu=(l-d)\mu_n$ and $b=b_U$, shows that the events
\begin{align*}
    \mathbf{E}_U^{\pm}=\Big\{\pm\sum_{j=1}^m\big(U_j^*-(l-d)\mu_n\big)\leq\frac{m}{m-2t_1}\Big(\sqrt{2t_1[\mathbb{V}]^{U*}+b_U^2t_1t_2}+\Big(1+\sqrt{\frac{2t_1}{m}}\Big)\Big(b_U\sqrt{t_1t_2}+\frac{b_Ut_1}{3}\Big)\Big)\Big\}
\end{align*}
satisfy $\P(\mathbf{E}_U^{\pm})\geq 1-e^{-t_1}-e^{-t_2}$ and $\P(\mathbf{E}_U^{+}\cap\mathbf{E}_U^{-})\geq 1-2e^{-t_1}-e^{-t_2}$. In the same way, with $W^*=(W_1^*,\dots,W_m^*)$ and $[\mathbb{V}]^{W*}=\|\Pi W^*\|^2$, Corollary \ref{corollary 1} (i), applied with $n$ replaced by $m$, $\mu=d\mu_n$, $b=b_W$ and $(t_3,t_4)$ in place of $(t_1,t_2)$, shows that the events
\begin{align*}
    \mathbf{E}_W^{\pm}=\Big\{\pm\sum_{j=1}^m\big(W_j^*-d\mu_n\big)\leq\frac{m}{m-2t_3}\Big(\sqrt{2t_3[\mathbb{V}]^{W*}+b_W^2t_3t_4}+\Big(1+\sqrt{\frac{2t_3}{m}}\Big)\Big(b_W\sqrt{t_3t_4}+\frac{b_Wt_3}{3}\Big)\Big)\Big\}
\end{align*}
satisfy $\P(\mathbf{E}_W^{\pm})\geq 1-e^{-t_3}-e^{-t_4}$ and $\P(\mathbf{E}_W^{+}\cap\mathbf{E}_W^{-})\geq 1-2e^{-t_3}-e^{-t_4}$.\\

\noindent \textbf{Step 3} (Proof of (i)) Since $\Pi$ is linear with operator norm one and $\|\cdot\|\leq\|\cdot\|_1$ on $\mathbb{R}^m$,
\begin{align*}
    \Big|\sqrt{[\mathbb{V}]^{U*}}-\sqrt{[\mathbb{V}]_m^U}\Big|\leq\sum_{j=1}^m|U_j-U_j^*|,\qquad\Big|\sqrt{[\mathbb{V}]^{W*}}-\sqrt{[\mathbb{V}]_m^W}\Big|\leq\sum_{j=1}^m|W_j-W_j^*|.
\end{align*}
On $\mathbf{G}$, both sums are at most $x$. The inequality in Step 2 of the proof of Theorem \ref{th 4}, applied with $\Gamma=0$, $\delta=x$ and $a^2=[\mathbb{V}]_m^U$ or $a^2=[\mathbb{V}]_m^W$, yields $\sqrt{2t_1[\mathbb{V}]^{U*}+b_U^2t_1t_2}\leq\sqrt{2t_1[\mathbb{V}]_m^U+b_U^2t_1t_2}+\sqrt{2t_1}\,x$ and $\sqrt{2t_3[\mathbb{V}]^{W*}+b_W^2t_3t_4}\leq\sqrt{2t_3[\mathbb{V}]_m^W+b_W^2t_3t_4}+\sqrt{2t_3}\,x$ on $\mathbf{G}$. Therefore, on $\mathbf{E}_U^{\pm}\cap\mathbf{E}_W^{\pm}\cap\mathbf{G}$,
\begin{align*}
    \pm n\big(\P_n(h_n)-\mu_n\big)\leq&\ \pm\sum_{j=1}^m\big(U_j^*-(l-d)\mu_n\big)\pm\sum_{j=1}^m\big(W_j^*-d\mu_n\big)\\
    &+\sum_{j=1}^m\big(|U_j-U_j^*|+|W_j-W_j^*|\big)\\
    \leq&\ \frac{m}{m-2t_1}\Big(\sqrt{2t_1[\mathbb{V}]_m^U+b_U^2t_1t_2}+\sqrt{2t_1}\,x+\Big(1+\sqrt{\frac{2t_1}{m}}\Big)\Big(b_U\sqrt{t_1t_2}+\frac{b_Ut_1}{3}\Big)\Big)\\
    &+\frac{m}{m-2t_3}\Big(\sqrt{2t_3[\mathbb{V}]_m^W+b_W^2t_3t_4}+\sqrt{2t_3}\,x+\Big(1+\sqrt{\frac{2t_3}{m}}\Big)\Big(b_W\sqrt{t_3t_4}+\frac{b_Wt_3}{3}\Big)\Big)\\
    &+x,
\end{align*}
which is $nL_n$, because $m/(n(m-2t))=1/(n-2tl)$ for $t\in\{t_1,t_3\}$. By Step 2, $\P(\mathbf{E}_U^{\pm}\cap\mathbf{E}_W^{\pm}\cap\mathbf{G})\geq 1-e^{-t_1}-e^{-t_2}-e^{-t_3}-e^{-t_4}-\|h_n\|_{\text{Lip}}\,n\,\tau_\infty(d)/x$, and the intersection of all five events has probability no less than the right-hand side of \eqref{eq th 6 1}. This proves (i).\\

\noindent \textbf{Step 4} (Preparation for (ii)) Denote $C=\sup_n\max\{b_n-a_n,\|h_n\|_{\text{Lip}}\}$, $T=\sum_{r=1}^\infty\tau_1(r)$ and $R(u)=\sum_{r\geq u}\tau_1(r)$ for $u\geq 1$. Let $i<k$. Since $\E[Y_k|\mathcal{F}_{-\infty}^i]=\int h_n\,d\P_{Z_k|\mathcal{F}_{-\infty}^i}-\int h_n\,d\P_{Z_k}$ almost surely and $h_n$ is $\|h_n\|_{\text{Lip}}$-Lipschitz, we have $\|\E[Y_k|\mathcal{F}_{-\infty}^i]\|_1\leq\|h_n\|_{\text{Lip}}\,\tau(\mathcal{F}_{-\infty}^i,Z_k)\leq\|h_n\|_{\text{Lip}}\,\tau_1(k-i)$, so that $|\E[Y_iY_k]|=|\E[Y_i\E[Y_k|\mathcal{F}_{-\infty}^i]]|\leq C^2\tau_1(k-i)$. Consequently, $\text{Var}(\sum_{i\in A}Y_i)\leq|A|\,C^2(1+2T)$ for every finite set $A$ of indices.
\begin{enumerate}[(a)]
    \item Let $\bar{\sigma}_n^2=\frac{1}{m(l-d)}\sum_{j=1}^m\text{Var}(U_j)$, so that $\bar{\sigma}_n^2\leq C^2(1+2T)$. For $i\in I_j$, let $s=e_j-i\in\{0,\dots,l-d-1\}$. Every $k\in I_{j'}$ with $j'>j$ satisfies $k-i\geq s+d+1$, so that $\sum_{j'>j}\sum_{k\in I_{j'}}|\E[Y_iY_k]|\leq C^2R(s+d+1)$. Summing over $i$ and $j$ yields $|\text{Var}(\sum_{j=1}^mU_j)-\sum_{j=1}^m\text{Var}(U_j)|\leq 2mC^2\sum_{u=1}^lR(u)$. Moreover, $\text{Var}(\sum_{j=1}^mU_j)\leq nC^2(1+2T)$ and $\text{Var}(\sum_{j=1}^mW_j)\leq mdC^2(1+2T)$. Since $\text{Var}(\sum_{i=1}^nh_n(Z_i))=\text{Var}(\sum_{j=1}^m(U_j+W_j))$, the Cauchy--Schwarz inequality yields
    \begin{align*}
        \big|\sigma_n^2-\big(1-\tfrac{d}{l}\big)\bar{\sigma}_n^2\big|\leq 2C^2\frac{1}{l}\sum_{u=1}^lR(u)+C^2(1+2T)\Big(\frac{d}{l}+2\sqrt{\frac{d}{l}}\Big).
    \end{align*}
    Since $d_n\geq 1$ and $d_n\log n/l_n\to 0$, we have $l_n\to\infty$ and $d_n/l_n\to 0$. Since $R(u)\to 0$ as $u\to\infty$, $\frac{1}{l}\sum_{u=1}^lR(u)\to 0$ as $l\to\infty$. Hence $\sigma_n^2-\bar{\sigma}_n^2\to 0$. In particular, $\sigma_n^2$ is bounded and $\liminf_n\bar{\sigma}_n^2=\liminf_n\sigma_n^2$.
    \item Since $l_n^2\log n/n\to 0$, $m=n/l_n\to\infty$ and $l_n/m=l_n^2/n\to 0$. With the choice in (ii), $t_1=\log(2/\alpha)$ does not depend on $n$, and we write $t_{2,n}=t_{3,n}=t_{4,n}=\log n$. Hence $2t_1<m$ and $2t_{3,n}<m$ for $n$ large enough.
    \item The right-hand side of \eqref{eq th 6 1} equals $1-\alpha-4/n-\|h_n\|_{\text{Lip}}\,n\,\tau_\infty(d_n)/x_n$, which is at least $1-\alpha-4/n-Cn\tau_\infty(d_n)/x_n$ and hence tends to $1-\alpha$.
    \item Since $[\mathbb{V}]_m^W\leq\sum_{j=1}^m(W_j-d\mu_n)^2$, we have $\E\big[[\mathbb{V}]_m^W\big]\leq\sum_{j=1}^m\text{Var}(W_j)\leq mdC^2(1+2T)$.
\end{enumerate}

\noindent \textbf{Step 5} (Consistency of $[\mathbb{V}]_m^U$) By Step 1 (a), $\E[\sum_{j=1}^m(U_j^*-(l-d)\mu_n)^2]=m(l-d)\bar{\sigma}_n^2$ and
\begin{align*}
    \text{Var}\Big(\sum_{j=1}^m\big(U_j^*-(l-d)\mu_n\big)^2\Big)\leq\sum_{j=1}^m\E\big[\big(U_j^*-(l-d)\mu_n\big)^4\big]\leq l^2C^2\sum_{j=1}^m\text{Var}(U_j)=l^2C^2m(l-d)\bar{\sigma}_n^2.
\end{align*}
Since $l-d\geq l/2$, the variance of $\sum_{j=1}^m(U_j^*-(l-d)\mu_n)^2/(m(l-d))$ is at most $2lC^4(1+2T)/m\to 0$. Moreover, with $\bar{U}^*=\frac{1}{m}\sum_{j=1}^mU_j^*$, $\E[m(\bar{U}^*-(l-d)\mu_n)^2]=(l-d)\bar{\sigma}_n^2$, so that $\E[m(\bar{U}^*-(l-d)\mu_n)^2]/n\leq\bar{\sigma}_n^2/m\to 0$. Since $[\mathbb{V}]^{U*}=\sum_{j=1}^m(U_j^*-(l-d)\mu_n)^2-m(\bar{U}^*-(l-d)\mu_n)^2$, $m(l-d)/n=1-d/l\to 1$ and $\bar{\sigma}_n^2$ is bounded, we obtain $[\mathbb{V}]^{U*}/n-\bar{\sigma}_n^2\xrightarrow[]{\P}0$. By Step 3 and Step 1,
\begin{align*}
    \E\Big|\sqrt{\frac{[\mathbb{V}]_m^U}{n}}-\sqrt{\frac{[\mathbb{V}]^{U*}}{n}}\Big|\leq\frac{Cn\tau_\infty(d)}{\sqrt{n}}=\frac{Cn\tau_\infty(d)}{x_n}\cdot\frac{x_n}{\sqrt{n}}\to 0.
\end{align*}
Together with $|\sqrt{y}-\sqrt{z}|\leq\sqrt{|y-z|}$ and Step 4 (a), this yields $\sqrt{[\mathbb{V}]_m^U/n}-\sigma_n\xrightarrow[]{\P}0$.\\

\noindent \textbf{Step 6} (Proof of (ii)) By \eqref{eq th 6 2} and $\sqrt{n}/(n-2tl)=m/(\sqrt{n}(m-2t))$,
\begin{align*}
    \sqrt{n}L_n=&\ \frac{m}{m-2t_1}\bigg(\sqrt{\frac{2t_1[\mathbb{V}]_m^U}{n}+\frac{b_U^2t_1t_{2,n}}{n}}+\frac{\sqrt{2t_1}\,x_n}{\sqrt{n}}\\
    &\qquad+\Big(1+\sqrt{\frac{2t_1}{m}}\Big)\Big(b_U\sqrt{\frac{t_1t_{2,n}}{n}}+\frac{b_Ut_1}{3\sqrt{n}}\Big)\bigg)\\
    &+\frac{m}{m-2t_{3,n}}\bigg(\sqrt{\frac{2t_{3,n}[\mathbb{V}]_m^W}{n}+\frac{b_W^2t_{3,n}t_{4,n}}{n}}+\frac{\sqrt{2t_{3,n}}\,x_n}{\sqrt{n}}\\
    &\qquad+\Big(1+\sqrt{\frac{2t_{3,n}}{m}}\Big)\Big(b_W\sqrt{\frac{t_{3,n}t_{4,n}}{n}}+\frac{b_Wt_{3,n}}{3\sqrt{n}}\Big)\bigg)+\frac{x_n}{\sqrt{n}}.
\end{align*}
By Step 4 (b), $m/(m-2t_1)\to 1$ and $m/(m-2t_{3,n})\to 1$. Since $b_U\leq Cl_n$ and $l_n^2\log n/n\to 0$, the terms involving $b_U$ tend to zero. Since $b_W\leq Cd_n$, $t_{3,n}=t_{4,n}=\log n$ and $d_n\log n\leq l_n$ for $n$ large enough, we have $b_W^2t_{3,n}t_{4,n}/n=O(l_n^2/n)\to 0$, so that the terms involving $b_W$ tend to zero as well. By Step 4 (d), $\E[t_{3,n}[\mathbb{V}]_m^W/n]\leq t_{3,n}(d_n/l_n)C^2(1+2T)\to 0$, since $d_n\log n/l_n\to 0$. Moreover, $\sqrt{t_{3,n}}\,x_n/\sqrt{n}\to 0$ by assumption. Together with Step 5, $t_1=\log(2/\alpha)$, the boundedness of $\sigma_n^2$ and $|\sqrt{y+z}-\sqrt{y}|\leq\sqrt{z}$, this yields $\sqrt{n}L_n-\sigma_n\sqrt{2\log(2/\alpha)}\xrightarrow[]{\P}0$. Since $\liminf_n\sigma_n^2>0$, this and Step 4 (c) prove (ii). Thus, we finish the proof of Theorem \ref{th 6}.\\

\subsection{Proof of Results in Section \ref{sec application}}
\noindent \textbf{Proof of Theorem \ref{th 7}}\\
\noindent \textbf{Step 0} (Notation) Denote $t=\log(1/\alpha)$, $T=\log n$, $\bar{x}_n=\sqrt{n}/T^{3/2}$, $b_U=(l-d)(b-a)$ and $b_W=d\,(b-a)$. Let $W_j=\sum_{i\in J_j}\ell_i$ and $[\mathbb{V}]_m^W=\sum_{j=1}^m(W_j-\frac{1}{m}\sum_{k=1}^mW_k)^2$, and let $L_n$ be the radius \eqref{eq th 6 2} with $h_n(Z_i)$ replaced by $\ell_i$, $b_n-a_n$ replaced by $b-a$, $t_1=t$ and $t_2=t_3=t_4=T$. Since $l=\lfloor n^{1/3}\rfloor$, $d=\lceil T^2\rceil$ and $m=n/l$, there exists $N_1$, depending only on $\alpha$ and $p$, such that, for every $n\geq N_1$,
\begin{align}
\label{eq proof th 7 0}
    n^{1/3}/2\leq l,\qquad p<d\leq l/2,\qquad m\geq 2,\qquad 2t<m,\qquad 2T<m,\qquad 4tl\leq n,\qquad 4Tl\leq n.
\end{align}
We assume $n\geq N_1$ throughout. Let $g(k)=(b-a)^2$ for $|k|\leq p$ and $g(k)=(p+1)(b-a)L\,\tau_\infty(|k|-p)$ for $|k|>p$, and denote $\bar{\Gamma}=\sum_{k\in\mathbb{Z}}g(k)$ and $\bar{\rho}_k=\sum_{j\in\mathbb{Z}}\min\{|j|/k,1\}\,g(j)$ for $k\geq 1$. Since $\sum_{r=1}^\infty\tau_\infty(r)<\infty$, we have $\bar{\Gamma}<\infty$ and, by dominated convergence, $\bar{\rho}_k\to 0$ as $k\to\infty$. Both depend only on $b-a$, $L$, $p$ and the sequence $\tau_\infty$.\\

\noindent \textbf{Step 1} (Covariances) Let $\gamma(k)=\mathrm{Cov}(\ell_0,\ell_k)$ for $k\in\mathbb{Z}$, so that $\gamma(-k)=\gamma(k)$ by stationarity. Let $k>p$. The loss $\ell_0$ is measurable with respect to $\mathcal{F}_{-\infty}^0$, and the first index of $(Z_{k-p},\dots,Z_k)$ is $k-p$. Since $\ell/L\in\Lambda_1(\mathbb{R}^{p+1})$, Definition \ref{def tau k} yields
\begin{align*}
    \big\|\E[\ell_k|\mathcal{F}_{-\infty}^0]-\mathcal{R}\big\|_1\leq L\,\tau\big(\mathcal{F}_{-\infty}^0,(Z_{k-p},\dots,Z_k)\big)\leq(p+1)L\,\tau_\infty(k-p).
\end{align*}
Hence $|\gamma(k)|=|\E[(\ell_0-\mathcal{R})(\E[\ell_k|\mathcal{F}_{-\infty}^0]-\mathcal{R})]|\leq(b-a)(p+1)L\,\tau_\infty(k-p)=g(k)$. For $|k|\leq p$, $|\gamma(k)|\leq(b-a)^2=g(k)$. Consequently, the series defining $\sigma^2$ converges absolutely. For $k\geq 1$, denote $S_k=\sum_{i=1}^k\ell_i$. By stationarity, $\text{Var}(S_k)=k\sum_{|j|<k}(1-|j|/k)\gamma(j)$, so that
\begin{align}
\label{eq proof th 7 1}
    \text{Var}(S_k)\leq k\bar{\Gamma},\qquad\Big|\frac{\text{Var}(S_k)}{k}-\sigma^2\Big|\leq\sum_{|j|<k}\frac{|j|}{k}|\gamma(j)|+\sum_{|j|\geq k}|\gamma(j)|\leq\bar{\rho}_k.
\end{align}

\noindent \textbf{Step 2} (Theorem \ref{th 6} for the losses) The proof of Theorem \ref{th 6} (i) applies to the losses $\ell_i$ in place of $h_n(Z_i)$, with common mean $\mathcal{R}$ and values in $[a,b]$, so that $|U_j-(l-d)\mathcal{R}|\leq b_U$ and $|W_j-d\mathcal{R}|\leq b_W$. Only the bounds in Step 1 (b) of that proof change. Since every $Z_k$ enters at most $p+1$ of the windows $(Z_{i-p},\dots,Z_i)$, the big-block sum $U_j$ is a $(p+1)L$-Lipschitz function, with respect to $\|\cdot\|_1$, of the $l-d+p$ observations $Z_k$, $(j-1)l+1-p\leq k\leq jl-d$, whereas $U_{j-1}$ depends only on the observations $Z_k$ with $k\leq e_{j-1}=(j-1)l-d$. Since $(j-1)l+1-p=e_{j-1}+d-p+1$, Definition \ref{def tau k} yields $\E|U_j-U_j^*|\leq(p+1)L(l-d+p)\,\tau_\infty(d-p+1)$. Likewise, $W_j$ is a $(p+1)L$-Lipschitz function of the $d+p$ observations $Z_k$, $jl-d+1-p\leq k\leq jl$, whereas $W_{j-1}$ depends only on the observations $Z_k$ with $k\leq(j-1)l$. Since $jl-d+1-p=(j-1)l+(l-d-p)+1$, we obtain $\E|W_j-W_j^*|\leq(p+1)L(d+p)\,\tau_\infty(l-d-p+1)$. Since $\tau_\infty$ is non-increasing and $l-d\geq d$ by \eqref{eq proof th 7 0}, the total coupling error satisfies
\begin{align*}
    \E\Big[\sum_{j=1}^m\big(|U_j-U_j^*|+|W_j-W_j^*|\big)\Big]&\leq(m-1)(p+1)L(l+2p)\,\tau_\infty(d-p)\\
    &\leq(p+1)L(n+2pm)\,\tau_\infty(d-p).
\end{align*}
Steps 2 and 3 of the proof of Theorem \ref{th 6} then apply without change. With $\mathbf{E}_U^-$, $\mathbf{E}_W^-$ and $\mathbf{G}$ the events defined there, the event $\mathbf{E}=\mathbf{E}_U^-\cap\mathbf{E}_W^-\cap\mathbf{G}$ satisfies $\P(\mathbf{E})\geq 1-e^{-t}-3e^{-T}-(p+1)L(n+2pm)\,\tau_\infty(d-p)/x$, that is,
\begin{align*}
    \P(\mathbf{E})\geq 1-\alpha-\frac{3}{n}-\frac{(p+1)L\,(n+2pm)\,\tau_\infty(d-p)}{x},
\end{align*}
and $\mathcal{R}-\hat{\mathcal{R}}_n\leq L_n$ on $\mathbf{E}$. Moreover, $U_1^*,\dots,U_m^*$ are independent by Step 1 (a) of that proof and, since $\mathbf{Z}$ is strictly stationary, they are identically distributed as $U_1$. Likewise, $W_1^*,\dots,W_m^*$ are independent and identically distributed as $W_1$. Finally, by Step 3 of that proof, $|\sqrt{[\mathbb{V}]^{U*}}-\sqrt{[\mathbb{V}]_m^U}|\leq x$ and $|\sqrt{[\mathbb{V}]^{W*}}-\sqrt{[\mathbb{V}]_m^W}|\leq x$ on $\mathbf{G}$.\\

\noindent \textbf{Step 3} (Two events for the empirical variances) Denote $s_U^2=\text{Var}(U_1)$ and $s_W^2=\text{Var}(W_1)$, so that $s_U^2=\text{Var}(S_{l-d})$ and $s_W^2=\text{Var}(S_d)$ by stationarity. The random variables $(U_j^*-(l-d)a)/b_U$, $1\leq j\leq m$, are independent, identically distributed and take values in $[0,1]$, and their sample variance in the sense of \citeA{maurer2009empirical} is $[\mathbb{V}]^{U*}/((m-1)b_U^2)$, with expectation $s_U^2/b_U^2$. Since $m\geq 2$, the first inequality of Theorem 10 of \citeA{maurer2009empirical}, applied with $\delta=1/n$, shows that the event $\mathbf{F}_U=\{\sqrt{[\mathbb{V}]^{U*}}\geq\sqrt{m-1}\,s_U-b_U\sqrt{2T}\}$ satisfies $\P(\mathbf{F}_U^c)\leq 1/n$. In the same way, the second inequality of that theorem shows that the event $\mathbf{F}_W=\{\sqrt{[\mathbb{V}]^{W*}}\leq\sqrt{m-1}\,s_W+b_W\sqrt{2T}\}$ satisfies $\P(\mathbf{F}_W^c)\leq 1/n$. By Step 2, $\P(\mathbf{E}\cap\mathbf{F}_U\cap\mathbf{F}_W)\geq 1-\alpha-\beta_n(x)$, and, on this event,
\begin{align}
\label{eq proof th 7 2}
    \sqrt{[\mathbb{V}]_m^U}\geq\sqrt{m-1}\,s_U-b_U\sqrt{2T}-x,\qquad\sqrt{[\mathbb{V}]_m^W}\leq\sqrt{m-1}\,s_W+b_W\sqrt{2T}+x.
\end{align}

\noindent \textbf{Step 4} (Deterministic bound on the radius) It remains to show that there exists $N_0\geq N_1$, depending only on $\alpha$, $b-a$, $L$, $p$, $\sigma^2$ and $\tau_\infty$, such that \eqref{eq proof th 7 2} implies $L_n\leq(1+1/\sqrt{T})\bar{L}_n(\alpha)$ for every $n\geq N_0$ and every $0<x\leq\bar{x}_n$. Then $\mathcal{R}\leq\hat{\mathcal{R}}_n+(1+1/\sqrt{T})\bar{L}_n(\alpha)$ on $\mathbf{E}\cap\mathbf{F}_U\cap\mathbf{F}_W$, which proves the first statement of Theorem \ref{th 7}.

Since $l-d\to\infty$, $d/l\to 0$ and $m\to\infty$, we may choose $N_0$ such that, for $n\geq N_0$, $\bar{\rho}_{l-d}\leq\sigma^2/8$, $(1-1/m)(1-d/l)\geq 9/14$, $b_U\sqrt{2T}\leq\sqrt{n}\,\sigma/8$ and $\bar{x}_n\leq\sqrt{n}\,\sigma/8$, where we use $b_U\leq n^{1/3}(b-a)$. By \eqref{eq proof th 7 1}, $s_U^2\geq(l-d)(\sigma^2-\bar{\rho}_{l-d})\geq\frac{7}{8}(l-d)\sigma^2$, so that $(m-1)s_U^2\geq\frac{7}{8}n(1-1/m)(1-d/l)\sigma^2\geq\frac{9}{16}n\sigma^2$. Together with \eqref{eq proof th 7 2} and $x\leq\bar{x}_n$, this yields $\sqrt{[\mathbb{V}]_m^U}\geq\frac{3}{4}\sqrt{n}\,\sigma-\frac{1}{4}\sqrt{n}\,\sigma=\frac{1}{2}\sqrt{n}\,\sigma$. Moreover, \eqref{eq proof th 7 1} yields $s_W^2\leq d\bar{\Gamma}$, so that \eqref{eq proof th 7 2} yields $\sqrt{[\mathbb{V}]_m^W}\leq\sqrt{md\bar{\Gamma}}+b_W\sqrt{2T}+\bar{x}_n$.

We now bound $L_n/\bar{L}_n(\alpha)$, where $\bar{L}_n(\alpha)=\sqrt{2t[\mathbb{V}]_m^U}/n$ and $\sqrt{2t[\mathbb{V}]_m^U}\geq\sigma\sqrt{tn/2}$. By \eqref{eq proof th 7 0}, $n/(n-2tl)\leq 1+4tl/n\leq 2$, $n/(n-2Tl)\leq 2$, $1+\sqrt{2t/m}\leq 2$ and $1+\sqrt{2T/m}\leq 2$. Since $\sqrt{1+u}\leq 1+u/2$ for $u\geq 0$, the first term of $L_n$ satisfies
\begin{align*}
    \frac{n}{n-2tl}\cdot\frac{\sqrt{2t[\mathbb{V}]_m^U+b_U^2tT}}{\sqrt{2t[\mathbb{V}]_m^U}}\leq\Big(1+\frac{4tl}{n}\Big)\Big(1+\frac{b_U^2T}{4[\mathbb{V}]_m^U}\Big)\leq 1+\frac{4tl}{n}+\frac{2b_U^2T}{n\sigma^2}.
\end{align*}
The remaining terms of $nL_n$ are at most $2\sqrt{2t}\,x+4b_U(\sqrt{tT}+t/3)$ for the big blocks, at most $2\sqrt{2T}(\sqrt{[\mathbb{V}]_m^W}+x)+2b_WT+\frac{16}{3}b_WT$ for the small blocks, and $x$. Inserting the bound on $\sqrt{[\mathbb{V}]_m^W}$ and $x\leq\bar{x}_n$, and dividing by $\sigma\sqrt{tn/2}$, we obtain
\begin{align*}
    \frac{L_n}{\bar{L}_n(\alpha)}-1\leq&\ \frac{4tl}{n}+\frac{2b_U^2T}{n\sigma^2}+\frac{\sqrt{2}\big(4b_U(\sqrt{tT}+t)+2\sqrt{2T}\sqrt{md\bar{\Gamma}}+12\,b_WT\big)}{\sigma\sqrt{tn}}\\
    &+\frac{\sqrt{2}\,(2\sqrt{2t}+4\sqrt{2T}+1)\,\bar{x}_n}{\sigma\sqrt{tn}}.
\end{align*}
Since $n^{1/3}/2\leq l\leq n^{1/3}$, $md=nd/l$, $d\leq T^2+1$, $b_U\leq n^{1/3}(b-a)$ and $b_W\leq(T^2+1)(b-a)$, the first line is at most a constant times $n^{-1/6}T^{3/2}$, where the constant depends only on $\alpha$, $b-a$, $\sigma^2$ and $\bar{\Gamma}$. Since $\bar{x}_n/\sqrt{n}=T^{-3/2}$, the second line equals $(4/T^{3/2}+8/(\sqrt{t}\,T)+\sqrt{2}/(\sqrt{t}\,T^{3/2}))/\sigma$. Hence $L_n/\bar{L}_n(\alpha)-1$ is at most a constant times $1/T$, with a constant depending only on $\alpha$, $b-a$, $\sigma^2$ and $\bar{\Gamma}$, and thus at most $1/\sqrt{T}$ for $n$ large enough. Enlarging $N_0$ accordingly completes the proof of the first statement of Theorem \ref{th 7}.\\

\noindent \textbf{Step 5} (Sharpness of the leading term) Let $c=(l-d)\mathcal{R}$ and $\bar{U}=\frac{1}{m}\sum_{j=1}^mU_j$, so that $[\mathbb{V}]_m^U=\sum_{j=1}^m(U_j-c)^2-m(\bar{U}-c)^2$ and $|U_j-c|\leq b_U$. For a finite set $A$ of indices, Step 1 yields $\text{Var}(\sum_{i\in A}\ell_i)=\sum_{i,k\in A}\gamma(k-i)\leq|A|\,\bar{\Gamma}$.
\begin{enumerate}[(a)]
    \item By stationarity, $\E[\sum_{j=1}^m(U_j-c)^2]=m\,\text{Var}(S_{l-d})$. Since $m(l-d)/n=1-d/l\to 1$ and $\bar{\rho}_{l-d}\to 0$, \eqref{eq proof th 7 1} yields $\E[\sum_{j=1}^m(U_j-c)^2]/n\to\sigma^2$. Moreover, $\E[m(\bar{U}-c)^2]=\text{Var}(\sum_{j=1}^mU_j)/m\leq(l-d)\bar{\Gamma}$, so that $\E[m(\bar{U}-c)^2]/n\leq l\bar{\Gamma}/n\to 0$.
    \item Since $(U_j-c)^2\leq b_U^2$, $\text{Var}((U_j-c)^2)\leq b_U^2\,\text{Var}(U_j)\leq b_U^2(l-d)\bar{\Gamma}\leq(b-a)^2\bar{\Gamma}\,l^3$. Let $j<k$. The variable $(U_j-c)^2$ is measurable with respect to $\mathcal{F}_{-\infty}^{e_j}$ with $e_j=jl-d$. Since $|u^2-v^2|\leq 2b_U|u-v|$ for $|u|,|v|\leq b_U$, Step 2 shows that $(U_k-c)^2$ is a $2b_U(p+1)L$-Lipschitz function of the $l-d+p$ observations $Z_s$, $(k-1)l+1-p\leq s\leq kl-d$, and $(k-1)l+1-p=e_j+(k-j-1)l+d-p+1$. Definition \ref{def tau k} therefore yields
    \begin{align*}
        \big|\text{Cov}\big((U_j-c)^2,(U_k-c)^2\big)\big|&\leq b_U^2\cdot 2b_U(p+1)L(l-d+p)\,\tau_\infty\big((k-j-1)l+d-p+1\big)\\
        &\leq K_1l^4\,\tau_\infty\big((k-j-1)l+d-p+1\big),
    \end{align*}
    with $K_1=4(p+1)L(b-a)^3$, where we use $b_U\leq l(b-a)$ and $l-d+p\leq 2l$. Since $\tau_\infty$ is non-increasing, $l\,\tau_\infty(rl+q)\leq\sum_{s=(r-1)l+q+1}^{rl+q}\tau_\infty(s)$ for $r\geq 1$, so that $\sum_{r\geq 0}\tau_\infty(rl+d-p+1)\leq\tau_\infty(d-p)+\frac{1}{l}\sum_{s\geq 1}\tau_\infty(s)$. Summing over the pairs $j<k$, of which at most $m$ satisfy $k-j-1=r$ for each $r\geq 0$, we obtain
    \begin{align*}
        \frac{1}{n^2}\text{Var}\Big(\sum_{j=1}^m(U_j-c)^2\Big)\leq\frac{(b-a)^2\bar{\Gamma}\,ml^3}{n^2}+\frac{2K_1ml^4}{n^2}\Big(\tau_\infty(d-p)+\frac{1}{l}\sum_{s\geq 1}\tau_\infty(s)\Big),
    \end{align*}
    and, since $ml=n$ and $l^3\leq n$, the right-hand side is at most $(b-a)^2\bar{\Gamma}\,l^2/n+2K_1(\tau_\infty(d-p)+l^{-1}\sum_{s\geq 1}\tau_\infty(s))$, which tends to zero because $l^2/n\to 0$, $d\to\infty$ and $l\to\infty$.
\end{enumerate}
By (a), (b) and Markov's inequality, $[\mathbb{V}]_m^U/n\xrightarrow[]{\P}\sigma^2$. Since $\bar{L}_n(\alpha)=\sqrt{2\log(1/\alpha)[\mathbb{V}]_m^U}/n$, this proves $\sqrt{n}\,\bar{L}_n(\alpha)\xrightarrow[]{\P}\sigma\sqrt{2\log(1/\alpha)}$. Thus, we finish the proof of Theorem \ref{th 7}.\\

\noindent \textbf{Proof of Corollary \ref{corollary 2}}\\
Let $N_0$ be as in the proof of Theorem \ref{th 7}, let $n\geq N_0$ and $0<x\leq\bar{x}_n$, and use the notation of that proof. By Steps 2 and 3 of that proof, $\mathcal{R}-\hat{\mathcal{R}}_n\leq L_n$ on $\mathbf{E}$ and $\P(\mathbf{E}\cap\mathbf{F}_U\cap\mathbf{F}_W)\geq 1-\alpha-\beta_n(x)$. It therefore suffices to show that, after enlarging $N_0$, $L_n\leq(1+1/\sqrt{T})\tilde{L}_n(\alpha)$ on $\mathbf{E}\cap\mathbf{F}_U\cap\mathbf{F}_W$.

Denote $Y=(Y_1,\dots,Y_m)$, $U=(U_1,\dots,U_m)$ and $W=(W_1,\dots,W_m)$, so that $Y=U+W$ and $[\mathbb{V}]_m^Y=\|\Pi Y\|^2$ in the notation of Step 0 of the proof of Theorem \ref{th 6}. Since $\Pi$ is linear, the triangle inequality yields
\begin{align}
\label{eq proof cor 2 1}
    \Big|\sqrt{[\mathbb{V}]_m^Y}-\sqrt{[\mathbb{V}]_m^U}\Big|\leq\sqrt{[\mathbb{V}]_m^W}.
\end{align}
By Step 4 of the proof of Theorem \ref{th 7}, the following holds on $\mathbf{E}\cap\mathbf{F}_U\cap\mathbf{F}_W$: $\sqrt{[\mathbb{V}]_m^U}\geq\frac{1}{2}\sqrt{n}\,\sigma$, $\sqrt{[\mathbb{V}]_m^W}\leq w_n$ with $w_n=\sqrt{md\bar{\Gamma}}+b_W\sqrt{2T}+\bar{x}_n$, and $L_n/\bar{L}_n(\alpha)-1\leq K/T$ for a constant $K$ depending only on $\alpha$, $b-a$, $\sigma^2$ and $\bar{\Gamma}$. Since $md=nd/l$, $d\leq T^2+1$, $l\geq n^{1/3}/2$ by \eqref{eq proof th 7 0} and $b_W\leq(T^2+1)(b-a)$, we have
\begin{align*}
    \frac{w_n}{\sqrt{n}}\leq\sqrt{2(T^2+1)\bar{\Gamma}}\,n^{-1/6}+\frac{(T^2+1)(b-a)\sqrt{2T}}{\sqrt{n}}+T^{-3/2}.
\end{align*}
After enlarging $N_0$, the right-hand side is at most $\sigma/4$, and it is at most $K'T^{-3/2}$ for a constant $K'$ depending only on $b-a$ and $\bar{\Gamma}$. By \eqref{eq proof cor 2 1}, $\sqrt{[\mathbb{V}]_m^Y}\geq\frac{1}{2}\sqrt{n}\,\sigma-w_n\geq\frac{1}{4}\sqrt{n}\,\sigma>0$, so that
\begin{align*}
    \frac{\bar{L}_n(\alpha)}{\tilde{L}_n(\alpha)}=\frac{\sqrt{[\mathbb{V}]_m^U}}{\sqrt{[\mathbb{V}]_m^Y}}\leq 1+\frac{\sqrt{[\mathbb{V}]_m^W}}{\sqrt{[\mathbb{V}]_m^Y}}\leq 1+\frac{4w_n}{\sqrt{n}\,\sigma}\leq 1+\frac{4K'}{\sigma T^{3/2}}.
\end{align*}
Hence $L_n/\tilde{L}_n(\alpha)=(L_n/\bar{L}_n(\alpha))(\bar{L}_n(\alpha)/\tilde{L}_n(\alpha))\leq(1+K/T)(1+4K'/(\sigma T^{3/2}))$, which is at most $1+1/\sqrt{T}$ after enlarging $N_0$. Finally, since $[\mathbb{V}]_m^W\leq\sum_{j=1}^m(W_j-d\mathcal{R})^2$, stationarity and \eqref{eq proof th 7 1} yield $\E[[\mathbb{V}]_m^W]/n\leq md\bar{\Gamma}/n=d\bar{\Gamma}/l\to 0$. Together with \eqref{eq proof cor 2 1} and Step 5 of the proof of Theorem \ref{th 7}, this yields $[\mathbb{V}]_m^Y/n\xrightarrow[]{\P}\sigma^2$, and hence $\sqrt{n}\,\tilde{L}_n(\alpha)\xrightarrow[]{\P}\sigma\sqrt{2\log(1/\alpha)}$. Thus, we finish the proof of Corollary \ref{corollary 2}.\\

\section{Relation Between the Coefficients $\tilde{\phi}_2$ and $\tau_2$}
\label{app B}
\setcounter{equation}{0}
\renewcommand\theequation{B.\arabic{equation}}
\setcounter{proposition}{0}
\renewcommand\theproposition{B.\arabic{proposition}}
\renewcommand\theHproposition{B.\arabic{proposition}}

Throughout this appendix, $\mathcal{M}\subset\mathcal{F}$ is a $\sigma$-algebra, and $X$ and $Y$ are real-valued random variables. The coefficients $\tilde{\phi}(\mathcal{M},\cdot)$, $\tilde{\phi}_1$, $\tilde{\phi}_2$, $\tau(\mathcal{M},\cdot)$, $\tau_1$ and $\tau_2$ are given in Definitions \ref{def phi tilde k} and \ref{def tau k}.

\begin{proposition}
\label{prop B.1}
\begin{enumerate}[(i)]
    \item (Monotone transformations) Let $I\subset\mathbb{R}$ be an interval such that $\P(X\in I)=\P(Y\in I)=1$, and let $G:I\to\mathbb{R}$ be strictly increasing. Then $\tilde{\phi}(\mathcal{M},G(X))\leq\tilde{\phi}(\mathcal{M},X)$ and $\tilde{\phi}(\mathcal{M},(G(X),G(Y)))\leq\tilde{\phi}(\mathcal{M},(X,Y))$. Consequently, if $\P(Z_t\in I)=1$ for every $t$, the coefficients $\tilde{\phi}_1$ and $\tilde{\phi}_2$ of $(G(Z_t))_{t\in\mathbb{Z}}$ are no larger than those of $\mathbf{Z}$.
    \item (Bounded variables) Suppose that $X$ and $Y$ take values in an interval of length $K$. Then
    \begin{align*}
        \tau(\mathcal{M},X)\leq K\tilde{\phi}(\mathcal{M},X),\qquad \tau(\mathcal{M},(X,Y))\leq 6K\sqrt{\tilde{\phi}(\mathcal{M},(X,Y))}.
    \end{align*}
    Consequently, if $\mathbf{Z}$ takes values in an interval of length $K$, then $\tau_1(r)\leq K\tilde{\phi}_1(r)$ and $\tau_2(r)\leq 3K\sqrt{\tilde{\phi}_2(r)}$ for every $r\geq 1$.
    \item (No linear bound) There is no constant $C$ such that $\tau(\mathcal{M},(X,Y))\leq C\,\tilde{\phi}(\mathcal{M},(X,Y))$ for every $\sigma$-algebra $\mathcal{M}$ and every pair $(X,Y)$ with values in $[0,1]^2$. More precisely, for every $\delta=1/(2N)$, $N\geq 1$, there exist $\mathcal{M}$ and such a pair with $\tilde{\phi}(\mathcal{M},(X,Y))\leq\delta^2$ and $\tau(\mathcal{M},(X,Y))\geq\delta/6$. In particular, the exponent $\frac{1}{2}$ in (ii) cannot be improved.
    \item (Example \ref{example 1}) The process $\mathbf{Z}$ of Example \ref{example 1} satisfies $\tilde{\phi}_1(r)\leq 2^{-r}$ and $\tilde{\phi}_2(r)\leq 2^{1-r}$ for every $r\geq 1$. The process $Z_t'=\tan(\pi(Z_t-\frac{1}{2}))$, $t\in\mathbb{Z}$, satisfies the same bounds, whereas $Z_t'$ has the standard Cauchy distribution, so that $\E|Z_t'|=\infty$ and $\tau(\mathcal{M},Z_t')$ is not defined.
    \item (Gaussian autoregression) Let $Z_t=aZ_{t-1}+\varepsilon_t$ with $0<|a|<1$ and i.i.d.\ standard normal innovations $(\varepsilon_t)_{t\in\mathbb{Z}}$, and let $\mathbf{Z}$ be its stationary solution. Then $\tilde{\phi}_1(r)=1$ for every $r\geq 1$.
\end{enumerate}
\end{proposition}

By (iv) and (v), neither of the two classes of processes contains the other. By (ii) and (iii), even for bounded processes, $\tilde{\phi}_2$ controls $\tau_2$ only through a square root, while $\tau_2$ does not control $\tilde{\phi}_2$ at all, as (v) shows for unbounded processes. That the Gaussian autoregression in (v) is $\tau$-mixing, with $\tau_2(r)\leq 2\E|Z_0||a|^r$, follows from Example 2 in Section 7.2 of \citeA{dedecker2005new}, applied with $g(x)=ax$.\\

\noindent \textbf{Proof of Proposition \ref{prop B.1}}\\
The univariate statements in (i), (ii) and (iv) are known: the first claim in (i) is a variant of Lemma 3 of \citeA{dedecker2005new}, the first inequality in (ii) follows from the representation of $\tau(\mathcal{M},X)$ in Section 1 of \citeA{dedecker2005new}, and the bound $\tilde{\phi}_1(r)\leq 2^{-r}$ in (iv) is shown in Example 1 of Section 4.1 of \citeA{dedecker2005new}. The remaining statements concern the bivariate coefficient $\tilde{\phi}(\mathcal{M},(X,Y))$ of Definition \ref{def phi tilde k}. For self-containedness, we nevertheless give complete proofs of all statements.\\
\noindent \textbf{Step 0} (Preliminaries) Let $\P_{(X,Y)|\mathcal{M}}$ be a regular conditional distribution of $(X,Y)$ given $\mathcal{M}$, and denote by $\lambda_\omega=\P_{(X,Y)|\mathcal{M}}(\omega,\cdot)-\P_{(X,Y)}$ the resulting signed measure on $\mathbb{R}^2$, so that $\lambda_\omega(\mathbb{R}^2)=0$. As in Step 1 of the proof of Lemma \ref{lemma A.6}, there exists a null set $N$ such that, for every $\omega\notin N$,
\begin{align}
\label{eq proof of prop B.1 0}
    \big|\lambda_\omega\big(I_1\times I_2\big)\big|\leq\tilde{\phi}(\mathcal{M},(X,Y))
\end{align}
for all sets $I_1,I_2$ of the form $(-\infty,s]$, $(-\infty,s)$ or $\mathbb{R}$, $s\in\mathbb{R}$: this holds first for closed half-lines with rational end points, then for all closed half-lines by continuity from above, and then for open half-lines and for $\mathbb{R}$ by continuity from below. The same statement holds for the marginal $\lambda_\omega^X=\P_{X|\mathcal{M}}(\omega,\cdot)-\P_X$ and $\tilde{\phi}(\mathcal{M},X)$.\\

\noindent \textbf{Step 1} (Proof of (i)) Since $G$ is strictly increasing on the interval $I$, for every $y\in\mathbb{R}$ the set $\{u\in I:G(u)\leq y\}$ is the intersection of $I$ with a set of the form $(-\infty,s]$, $(-\infty,s)$, $\emptyset$ or $\mathbb{R}$. Since $\P(X\in I)=1$ and $\P_{X|\mathcal{M}}(\omega,I)=1$ for almost every $\omega$, we obtain, almost surely, $|\P(G(X)\leq y|\mathcal{M})-\P(G(X)\leq y)|=|\lambda^X_\omega(I_1)|$ for such a set $I_1$, which is no larger than $\tilde{\phi}(\mathcal{M},X)$ by Step 0. Taking the supremum over $y$ proves the first claim. The bivariate claim follows in the same way from \eqref{eq proof of prop B.1 0}, because $\{(u,v)\in I^2:G(u)\leq y_1,G(v)\leq y_2\}$ is the intersection of $I^2$ with a product $I_1\times I_2$ as in Step 0. The claim on $\tilde{\phi}_1$ and $\tilde{\phi}_2$ follows from the definitions, since $\sigma(G(Z_t):t\leq p)\subset\mathcal{F}_{-\infty}^p$ and, by Step 3 of the proof of Lemma \ref{lemma A.6}, the coefficients do not increase when $\mathcal{F}_{-\infty}^p$ is replaced by a smaller $\sigma$-algebra.\\

\noindent \textbf{Step 2} (Proof of (ii), univariate) Let $X$ take values in $[c-K,c]$ and let $f\in\Lambda_1(\mathbb{R})$. Since $f$ is absolutely continuous with $|f'|\leq 1$ almost everywhere, Fubini's theorem yields, for every probability measure $\mu$ on $[c-K,c]$,
\begin{align*}
    \int f\,d\mu=f(c)-\int_{c-K}^cf'(s)\,\mu\big((-\infty,s]\big)\,ds.
\end{align*}
Applying this to $\P_{X|\mathcal{M}}(\omega,\cdot)$ and $\P_X$ and using Step 0, we obtain $|\int f\,d\lambda^X_\omega|\leq K\tilde{\phi}(\mathcal{M},X)$ for every $\omega\notin N$ and every $f\in\Lambda_1(\mathbb{R})$. Hence $\tau(\mathcal{M},X)\leq K\tilde{\phi}(\mathcal{M},X)$.\\

\noindent \textbf{Step 3} (Proof of (ii), bivariate) Let $(X,Y)$ take values in $[c-K,c]^2$, denote $\varphi=\tilde{\phi}(\mathcal{M},(X,Y))$, and let $f\in\Lambda_1(\mathbb{R}^2)$. If $\varphi=0$, then $\lambda_\omega=0$ for $\omega\notin N$ by \eqref{eq proof of prop B.1 0}, and there is nothing to prove. Otherwise, let $\kappa(u)=(1-|u|)_+$, $\kappa_\varepsilon(u)=\varepsilon^{-1}\kappa(u/\varepsilon)$ for $\varepsilon>0$, and
\begin{align*}
    f_\varepsilon(x,y)=\int\!\!\int f(x-u,y-v)\,\kappa_\varepsilon(u)\kappa_\varepsilon(v)\,du\,dv.
\end{align*}
Since $f$ is $1$-Lipschitz with respect to $\|\cdot\|_1$ and $\int|u|\kappa_\varepsilon(u)\,du=\varepsilon/3$, we have $|f_\varepsilon-f|\leq 2\varepsilon/3$. Moreover, $f$ is $1$-Lipschitz in each coordinate, so that its partial derivatives exist almost everywhere and are bounded by one. Hence $\partial_1f_\varepsilon=(\partial_1 f)*(\kappa_\varepsilon\otimes\kappa_\varepsilon)$ and $\partial_2\partial_1f_\varepsilon=(\partial_1 f)*(\kappa_\varepsilon\otimes\kappa_\varepsilon')$ are continuous, $|\partial_1f_\varepsilon|\leq 1$, $|\partial_2f_\varepsilon|\leq 1$ and $|\partial_2\partial_1f_\varepsilon|\leq\int|\kappa_\varepsilon'|=2/\varepsilon$. For $(x,y)\in[c-K,c]^2$,
\begin{align*}
    f_\varepsilon(x,y)=f_\varepsilon(c,c)-\int_x^c\partial_1f_\varepsilon(s,c)\,ds-\int_y^c\partial_2f_\varepsilon(c,t)\,dt+\int_x^c\!\!\int_y^c\partial_2\partial_1f_\varepsilon(s,t)\,dt\,ds.
\end{align*}
Integrating with respect to $\lambda_\omega$, which is concentrated on $[c-K,c]^2$ for almost every $\omega$ and satisfies $\lambda_\omega(\mathbb{R}^2)=0$, Fubini's theorem yields
\begin{align*}
    \int f_\varepsilon\,d\lambda_\omega=&-\int_{c-K}^c\partial_1f_\varepsilon(s,c)\,\lambda_\omega\big((-\infty,s]\times\mathbb{R}\big)\,ds-\int_{c-K}^c\partial_2f_\varepsilon(c,t)\,\lambda_\omega\big(\mathbb{R}\times(-\infty,t]\big)\,dt\\
    &+\int_{c-K}^c\!\!\int_{c-K}^c\partial_2\partial_1f_\varepsilon(s,t)\,\lambda_\omega\big((-\infty,s]\times(-\infty,t]\big)\,dt\,ds.
\end{align*}
By \eqref{eq proof of prop B.1 0}, for almost every $\omega$,
\begin{align*}
    \Big|\int f\,d\lambda_\omega\Big|\leq\Big|\int (f-f_\varepsilon)\,d\lambda_\omega\Big|+\Big|\int f_\varepsilon\,d\lambda_\omega\Big|\leq\frac{4\varepsilon}{3}+2K\varphi+\frac{2K^2\varphi}{\varepsilon}.
\end{align*}
The choice $\varepsilon=K\sqrt{3\varphi/2}$ yields $|\int f\,d\lambda_\omega|\leq K(\sqrt{32/3}\sqrt{\varphi}+2\varphi)\leq 6K\sqrt{\varphi}$, because $\varphi\leq 1$. Since this holds for every $f\in\Lambda_1(\mathbb{R}^2)$ and almost every $\omega$, $\tau(\mathcal{M},(X,Y))\leq 6K\sqrt{\varphi}$.

For the consequence, Definitions \ref{def phi tilde k} and \ref{def tau k} and the bounds just proved yield $\tau_1(r)\leq K\tilde{\phi}_1(r)$ and $\frac{1}{2}\tau(\mathcal{F}_{-\infty}^p,(Z_i,Z_j))\leq 3K\sqrt{\tilde{\phi}_2(r)}$ for $p+r\leq i<j$. Since $\tilde{\phi}_1(r)\leq\tilde{\phi}_2(r)\leq 1$, we also have $K\tilde{\phi}_1(r)\leq 3K\sqrt{\tilde{\phi}_2(r)}$, which proves $\tau_2(r)\leq 3K\sqrt{\tilde{\phi}_2(r)}$.\\

\noindent \textbf{Step 4} (Proof of (iii)) Let $N\geq 1$, $\delta=1/(2N)$ and $s(u)=(-1)^{\lfloor u/\delta\rfloor}$ for $u\in\mathbb{R}$. Let $\xi$ be a random variable with $\P(\xi=1)=\P(\xi=-1)=\frac{1}{2}$, let $\mathcal{M}=\sigma(\xi)$, and, conditionally on $\xi$, let $(X,Y)$ have the density $1+\xi s(x)s(y)$ on $[0,1]^2$, which is non-negative. Unconditionally, $(X,Y)$ is uniformly distributed on $[0,1]^2$. Denote $S(x)=\int_0^xs(u)\,du$, so that $0\leq S(x)\leq\delta$ for every $x\in[0,1]$. Since
\begin{align*}
    \P(X\leq x,Y\leq y|\mathcal{M})-\P(X\leq x,Y\leq y)=\xi S(x)S(y),\qquad x,y\in[0,1],
\end{align*}
and the left-hand side vanishes if $x$ or $y$ lies outside $[0,1]$, we have $\tilde{\phi}(\mathcal{M},(X,Y))\leq\delta^2$. On the other hand, let $m(x,y)=\min\{\text{dist}(x,\delta\mathbb{Z}),\text{dist}(y,\delta\mathbb{Z})\}$ and $q(x,y)=s(x)s(y)m(x,y)$. On each closed square $[i\delta,(i+1)\delta]\times[j\delta,(j+1)\delta]$, $i,j\in\mathbb{Z}$, the function $q$ coincides with $\pm m$, which is $1$-Lipschitz with respect to $\|\cdot\|_1$, and $q$ vanishes on the boundaries of these squares. Hence $q$ is continuous, and splitting any segment at the boundaries of the squares that it crosses shows that $q\in\Lambda_1(\mathbb{R}^2)$. Since $s(x)^2s(y)^2=1$ and the average of $m$ over each square equals $\delta/6$,
\begin{align*}
    \int q\,d\P_{(X,Y)|\mathcal{M}}-\int q\,d\P_{(X,Y)}=\xi\int_0^1\!\!\int_0^1 m(x,y)\,dx\,dy=\xi\frac{\delta}{6},
\end{align*}
so that $\tau(\mathcal{M},(X,Y))\geq\delta/6$. Here we used that, for $U$ and $V$ independent and uniformly distributed on $[0,1]$, $\P(\min\{U,1-U,V,1-V\}>t)=(1-2t)^2$ for $t\in[0,\frac{1}{2}]$, so that $\E[\min\{U,1-U,V,1-V\}]=\int_0^{1/2}(1-2t)^2dt=\frac{1}{6}$. Therefore $\tau(\mathcal{M},(X,Y))/\tilde{\phi}(\mathcal{M},(X,Y))\geq 1/(6\delta)\to\infty$ as $N\to\infty$, and, for every $\beta>\frac{1}{2}$, $\tau(\mathcal{M},(X,Y))/\tilde{\phi}(\mathcal{M},(X,Y))^\beta\geq\delta^{1-2\beta}/6\to\infty$.\\

\noindent \textbf{Step 5} (Proof of (iv)) Let $p\in\mathbb{Z}$ and $j\geq p+r$. We have $Z_j=2^{-(j-p)}Z_p+U$ with $U=\sum_{k=0}^{j-p-1}2^{-(k+1)}\varepsilon_{j-k}$, which is independent of $\mathcal{F}_{-\infty}^p$ and uniformly distributed on the grid $\{kh:0\leq k<2^{j-p}\}$ with mesh $h=2^{-(j-p)}\leq 2^{-r}$. Hence, almost surely, $\P(Z_j\leq z|\mathcal{F}_{-\infty}^p)=\P(c+U\leq z)|_{c=hZ_p}$. For $c\in[0,h]$ and $z\in[0,1]$, the number of grid points $kh$ with $c+kh\leq z$ is $0$ if $z<c$ and $\min\{\lfloor (z-c)/h\rfloor+1,2^{j-p}\}$ otherwise, so that $|\P(c+U\leq z)-z|\leq h$. For $z\notin[0,1]$, both probabilities coincide. Since $Z_j$ is uniformly distributed on $[0,1]$, this proves $\tilde{\phi}(\mathcal{F}_{-\infty}^p,Z_j)\leq 2^{-r}$ and hence $\tilde{\phi}_1(r)\leq 2^{-r}$.

Let now $p+r\leq i<j$. We have $Z_j=2^{-(j-i)}Z_i+V$ with $V=\sum_{k=0}^{j-i-1}2^{-(k+1)}\varepsilon_{j-k}$, which is independent of $\mathcal{F}_{-\infty}^i\supset\mathcal{F}_{-\infty}^p$. For $x,y\in\mathbb{R}$, let $g(u)=1\{u\leq x\}\P(V\leq y-2^{-(j-i)}u)$. Then $\P(Z_i\leq x,Z_j\leq y|\mathcal{F}_{-\infty}^p)=\E[g(Z_i)|\mathcal{F}_{-\infty}^p]$ almost surely and $\P(Z_i\leq x,Z_j\leq y)=\E[g(Z_i)]$. The function $u\mapsto\P(V\leq y-2^{-(j-i)}u)$ is non-increasing with values in $[0,1]$, so that $g$ is of bounded variation with $\|g\|_{TV}\leq 2$. Lemma \ref{lemma A.6} (i) and the first part of this step yield
\begin{align*}
    \big|\P(Z_i\leq x,Z_j\leq y|\mathcal{F}_{-\infty}^p)-\P(Z_i\leq x,Z_j\leq y)\big|\leq 2\tilde{\phi}(\mathcal{F}_{-\infty}^p,Z_i)\leq 2^{1-r}\ \ \text{a.s.}
\end{align*}
Together with $\tilde{\phi}_1(r)\leq 2^{-r}$, this proves $\tilde{\phi}_2(r)\leq 2^{1-r}$.

Finally, $G(u)=\tan(\pi(u-\frac{1}{2}))$ is strictly increasing on $I=(0,1)$ and $\P(Z_t\in I)=1$, so that (i) yields the same bounds for $(Z_t')_{t\in\mathbb{Z}}$. Since $Z_t$ is uniformly distributed on $[0,1]$, $Z_t'$ has the standard Cauchy distribution, and $\E|Z_t'|=\infty$.\\

\noindent \textbf{Step 6} (Proof of (v)) The stationary solution is $Z_t=\sum_{k\geq 0}a^k\varepsilon_{t-k}$, so that $Z_t$ is normal with mean zero and variance $s_0^2=1/(1-a^2)$. For $p\in\mathbb{Z}$ and $r\geq 1$, we have $Z_{p+r}=a^rZ_p+W$ with $W=\sum_{k=0}^{r-1}a^k\varepsilon_{p+r-k}$, which is independent of $\mathcal{F}_{-\infty}^p$ and normal with mean zero and variance $s_r^2=(1-a^{2r})/(1-a^2)>0$. Hence, with $\Phi$ the standard normal distribution function, $\P(Z_{p+r}\leq z|\mathcal{F}_{-\infty}^p)=\Phi((z-a^rZ_p)/s_r)$ almost surely. For fixed $z$, the map $u\mapsto|\Phi((z-u)/s_r)-\Phi(z/s_0)|$ is continuous, and $a^rZ_p$ is normal with positive variance, so that its support is $\mathbb{R}$. Therefore
\begin{align*}
    \big\|\P(Z_{p+r}\leq z|\mathcal{F}_{-\infty}^p)-\P(Z_{p+r}\leq z)\big\|_\infty&=\sup_{u\in\mathbb{R}}\big|\Phi((z-u)/s_r)-\Phi(z/s_0)\big|\\
    &=\max\{\Phi(z/s_0),1-\Phi(z/s_0)\},
\end{align*}
which tends to one as $z\to\infty$. Hence $\tilde{\phi}(\mathcal{F}_{-\infty}^p,Z_{p+r})=1$, and $\tilde{\phi}_1(r)=1$ for every $r\geq 1$. Thus, we finish the proof of Proposition \ref{prop B.1}.

\end{document}